\documentclass[12pt,a4paper]{article}

\usepackage{amsmath}    
\usepackage{amssymb}
\usepackage{graphicx}   
\usepackage{verbatim}   
\usepackage{color}      

\usepackage{hyperref}   

\usepackage{enumerate}
\usepackage{mathrsfs}
\usepackage{empheq}
\usepackage{bbold}
\usepackage{multirow}
\usepackage[english]{babel} 

\usepackage[final]{showlabels}

\usepackage{amsthm}
\newtheorem{definition}{Definition}[section]

\newtheorem{corollary}{Corollary}[section]
\newtheorem{theorem}{Theorem}[section]

\newtheorem{remark}{Remark}[section]

\newtheorem{scheme}{Scheme}

\newcommand{\del}{\partial}
\renewcommand{\theta}{\vartheta}
\renewcommand{\phi}{\varphi}
\newcommand{\vecc}[2]{\left ( \begin{array}{c}#1\\#2\\ \end{array}\right )}
\newcommand{\veccc}[3]{\left ( \begin{array}{c}#1\\#2\\#3\\ \end{array}\right )}

\newcommand{\dd}{\mathrm{d}}

\newcommand{\id}{\mathbb{1}}

\renewcommand{\div}{\mathrm{div\,}}
\renewcommand{\vec}{\mathbf}

\newcommand{\const}{\mathrm{const}}

\newcommand{\ii}{\mathbb{i}}
\newcommand{\ee}{\mathrm{e}}

\newcommand{\curl}{\text{curl\,}}
\newcommand{\sgn}{\mathrm{sgn}}

\renewcommand{\title}{All-speed numerical methods for the Euler equations via a sequential explicit time integration}

\newcommand{\authorOne}{Wasilij Barsukow\footnote{Bordeaux Institute of Mathematics, Bordeaux University and CNRS/UMR5251, Talence, 33405 France}}

\begin{document}

\begin{center} \Large
\title

\vspace{1cm}

\date{}
\normalsize

\authorOne
\end{center}

\begin{abstract}

This paper presents a new strategy to deal with the excessive diffusion that standard finite volume methods for compressible Euler equations display in the limit of low Mach number. The strategy can be understood as using centered discretizations for the acoustic part of the Euler equations and stabilizing them with a leap-frog-type (``sequential explicit'') time integration, a fully explicit method. This time integration takes inspiration from time-explicit staggered grid numerical methods. In this way, advantages of staggered methods carry over to collocated methods. The paper provides a number of new collocated schemes for linear acoustic/Maxwell equations that are inspired by the Yee scheme. They are then extended to an all-speed method for the full Euler equations on Cartesian grids. By taking the opposite view and taking inspiration from collocated methods, the paper also suggests a new way of staggering the variables which increases the stability as compared to the traditional Yee scheme.



\end{abstract}

\section{Introduction}

Numerical methods for compressible flow for a long time have been focusing on stability under explicit time-integration (upwinding) and shock-capturing. The very idea of obtaining the numerical flux through the introduction of Riemann problems at cell interfaces exemplifies the focus on supersonic phenomena. Besides refining the grid, the standard way of obtaining more accurate methods is higher order of approximation (e.g. higher degree polynomials). Both strategies come at the cost of increased computational time and memory.

As was shown in e.g. \cite{ebin77,klainerman81,metivier01}, in the limit of low Mach number, the solutions to the Euler equations fulfill the incompressible Euler equations for well-prepared data. For the numerical solution of low Mach number flow, \textbf{two problems} arise. First, the explicit time step becomes small because it is computed with respect to the (fast) speed of sound, while relevant time scales involve the (slow) speed of the fluid: the problem becomes stiff. But even if one accepts to wait long, the upwind-based space discretization adds too much diffusion, which is the second problem. Arguments based on asymptotic analysis (e.g. in \cite{guillard99}) demonstrate how standard methods fail to comply with the low Mach number limit at finite discretization by introducing inadequate pressure fluctuations. In practice this means that they require excessive grid refinement to resolve low Mach number phenomena. It is this latter problem that the present work addresses. It thus pursues the question how the numerical diffusion can be reduced in the limit of low Mach number while retaining enough of it for an explicit integration in time.

If it is a priori known that the flow has a very low Mach number, and if there is no interest in resolving compressible phenomena (such as sound waves or a density/pressure stratification), then one might directly consider using an incompressible, or weakly-compressible code. Otherwise, an \textbf{all-speed} numerical method is necessary, which can deal with all regimes simultaneously.

Implicit time integration can help bridging the disparity of the advective and acoustic time scales in the low Mach number regime and allows time steps to be based on a CFL condition involving just the advective time scale. However, even in flows with shocks, there are low Mach number regions (e.g. turbulent wakes or unstable slip lines) that would be worth resolving. In all-speed regimes, where low Mach number regions coexist with shock waves or if one is interested in resolving sound waves, the CFL condition of an explicit method does not pose a restriction: out of accuracy considerations the time step would be chosen based on the acoustic time scale anyway. In this work, only \textbf{explicit time stepping} is considered.

Several ways of constructing low Mach number, or even all-speed numerical methods are already well-established in the literature. Formal asymptotic analysis of finite volume schemes for the compressible Euler equations suggests that certain terms in the numerical diffusion (which arises via upwinding, or through the usage of a Riemann solver) are not compliant with the low Mach number limit (see e.g. \cite{dellacherie10,barsukow16} for more details). As early as in \cite{turkel87,weiss95} it has been found that the accuracy of standard compressible methods in the regime of low Mach number could be drastically increased by modifying, or even removing these terms. Thus, there is potential to dramatically increase the accuracy in the low Mach number regime -- at no extra computational cost.

Numerical methods which modify those terms, and make them vanish as the Mach number approaches zero, are referred to as \textbf{low Mach fixes}. Through these fixes certain parts of the equations are discretized asymptotically centrally, and a large number of fixes has been found, which still possess sufficient numerical diffusion to be amenable to explicit time integration: e.g. \cite{li08,thornber08,dellacherie10,rieper11,li13,chalons16,birken16,barsukow16,dellacherie16} (sometimes with a more severe CFL condition, as shown in e.g. \cite{birken05,barsukow16}). So far virtually all explicit all-speed methods employed a ``low Mach fix'' strategy, which -- if giving rise to a stable method -- involves free parameters and lacks a first-principles derivation.

In \cite{barsukow18thesis,barsukow20cgk} another way of constructing explicit all-speed schemes was presented. The key observation is that the low Mach number limit is only non-trivial in multiple spatial dimensions. One-dimensional numerical methods need to be modified through a low Mach fix only because they are subsequently applied to multiple dimensions in a dimensionally split fashion. In \cite{barsukow18thesis,barsukow20cgk}, one-dimensional numerical methods were instead \textbf{extended to multiple dimensions in a} very particular, \textbf{all-speed way}, leaving the one-dimensional scheme as it is. The choice of this multi-dimensional extension is inspired by vorticity-preserving numerical methods for \textbf{linear acoustics}, which were shown to be low Mach number compliant in \cite{barsukow17a}, and the methods from \cite{barsukow18thesis,barsukow20cgk} can be understood as one way of extending those results from linear acoustics to the fully nonlinear Euler equations. This approach results in enhanced stability and does not require free parameters. The particular multi-dimensional extension needs to be chosen carefully; the use of an \emph{exact}, truly multi-dimensional Riemann solver, for example, was shown in \cite{barsukow17} \emph{not} to be low Mach compliant.

Although the focus of the present work is on explicit time integration, it is still insightful to review those approaches that resort to implicit time integration. Although they primarily target the first problem of stiffness in time, the natural choice of central discretization in space, i.e. the absence of upwinding, simultaneously relieves them from the second problem. There exist fully implicit treatments such as \cite{viallet11,miczek15,abbate19}, however they still typically choose a time step based on the advective time scale for accuracy. Therefore, splittings (\textbf{IMEX/semi-implicit methods}) seem more efficient where terms associated with acoustics are solved implicitly, while those terms whose upwinding is not harmful for the low Mach number limit (advection) are solved explicitly in time (e.g. \cite{degond07,degond11,cordier12,haack12,dimarco17,bispen17,boscarino19,boscheri20,thomann20,boscheri21,boscheri21a}). These methods use central derivatives in space.

Time-implicit discretizations of the compressible Euler equations on \textbf{staggered grids} are inspired by the MAC method \cite{harlow65} for incompressible flow, which was extended to include compressibility effects in \cite{harlow71,casulli84,karki89,shyy92,bijl96,roller00,wenneker02,munz03}. In \cite{park05}, the conservative form of the equations allowed to compute shocks, and examples of the usage of staggered grids to achieve the all-speed property with time-implicit methods are \cite{dumbser16,dumbser19}. 

Staggered finite differences are essentially central, i.e. they do not include upwinding and are generally not stable under explicit time integration. It has been mentioned above that the equations of linear acoustics serve as an excellent test bed for questions of involution preservation and low Mach number compliance. In two spatial dimensions, they are equivalent to Maxwell's equations. Interestingly, for these latter the \textbf{Yee method} \cite{yee66} employs staggered grids while being essentially explicit. The time integration is of leap-frog type, and appears as a natural time discretization of the (second order) wave equation as early as in \cite{courant28,lax56}, for example. If this time integration can be used for the Euler equations, one would obtain an explicitly-stable staggered/central difference method, which is all-speed because it does not contain the usual type of upwinding. To show a way how this can be done is the aim of the paper.

Starting from Yee's method for Maxwell's equations, in the first part of this paper some of the conceptual links, analogies and differences between staggered and collocated time-explicit methods are elucidated. In particular, an improvement of the staggered-grid Yee scheme is suggested (Scheme \ref{scheme:yeeextended}), which takes inspiration from multi-dimensional vorticity-preserving collocated numerical methods for linear acoustics. In the second part of the paper, a path towards time-explicit all-speed numerical methods for the Euler equations is shown. The new method presented here uses an \textbf{explicit stabilization of central derivatives} of leap-frog type inspired by Yee's method, well-known in computational electromagnetism (see e.g. \cite{pinto14}) and in connection with Hamiltonian systems. In this paper its usage is extended to collocated numerical methods and for problems without a Hamiltonian or an energy to be conserved: The special kind of time integration is used merely to stabilize the acoustic part of the Euler equations, for which one would like to get rid of upwinding. The result is a new \textbf{collocated all-speed numerical method for the Euler equations} (Scheme \ref{scheme:euler}), which is easy to implement as it is similar to Riemann-solver-based Finite Volume methods. Table \ref{tab:overview} shows an overview of the numerical methods proposed in this paper.

\begin{table}
 \centering
 \begin{tabular}{c|c|c||c|c|c}
  \textbf{Method} & \textbf{grid}& \textbf{CFL$_\text{max}$} & \textbf{Maxwell} & \textbf{Acoustics} & \textbf{Euler} \\\hline
  original Yee (\ref{scheme:yeeoriginal})  &sta & $1/\sqrt{2}$ &  \eqref{eq:yeeorig1}--\eqref{eq:yeeorig3} & \eqref{eq:yeeorigacoustics1}--\eqref{eq:yeeorigacoustics3} & -- \\
  Yee (\ref{scheme:yeecollocated}) & col & $1/\sqrt{2}$& \eqref{eq:yeecolloc1}--\eqref{eq:yeecolloc3} & & --\\
  explicit Yee (\ref{scheme:yeecollocatedexplicit}) &col & $1/\sqrt{2}$ &\eqref{eq:yeecollocexpl1}--\eqref{eq:yeecollocexpl3} & & --  \\
  \hline
  Yee extended (\ref{scheme:yeecollocatedextended}) & col & 1 &\eqref{eq:yeecollocextended1}--\eqref{eq:yeecollocextended3} & \eqref{eq:yeecollocatedextacoustics1}--\eqref{eq:yeecollocatedextacoustics3} & -- \\ 
  Yee extented (\ref{scheme:yeeextended}) & sta & 1&\eqref{eq:yeeextended1}--\eqref{eq:yeeextended3} & & --  \\
  Yee extended 3D (\ref{scheme:yeeextended3d}) & sta & 1&\eqref{eq:yeeextended3d1}--\eqref{eq:yeeextended3d3} & & --  \\
  central (\ref{scheme:central}) &col & $\sqrt{2}$& \eqref{eq:centralmaxwell1}--\eqref{eq:centralmaxwell3} & & --  \\
  central extended (\ref{scheme:centralextended}) &col & 2 & \eqref{eq:centralextended1}--\eqref{eq:centralextended3} & \eqref{eq:centralextendedacoustics1}--\eqref{eq:centralextendedacoustics3} & -- \\
  \hline
  central extended $+$  & \multirow{2}{*}{col}&  \multirow{2}{*}{1} &  \multirow{2}{*}{--} &  \multirow{2}{*}{--} &  \multirow{2}{*}{\eqref{eq:updateeulermomentum}--\eqref{eq:updateeulerene}} \\
   pressureless (\ref{scheme:euler}) & &&&&
 \end{tabular}
 \caption{Overview of the numerical methods (capital letters) presented in this paper (besides \ref{scheme:yeeoriginal}, it seems that only \ref{scheme:central} has been suggested before as the FVTD method \cite{remaki99,piperno02}). They all are inspired by the original Yee scheme, but this work places a focus on collocated methods (see second column). Some of the suggestions of this work improve the stability constraint of the original Yee scheme (third column, the values are valid for two spatial dimensions). Method \ref{scheme:euler} is for the full Euler equations, while the others can be used for Maxwell's and the acoustic equations. The word ``extended'' refers to a multi-dimensional enlargement of the stencil. The numerical methods for different systems are referenced via the equation numbers.}
 \label{tab:overview}
\end{table}

\begin{table}
 \centering
 \begin{tabular}{c||c}
  \textbf{Method} &  \textbf{Fourier transform}\\\hline
  original Yee (\ref{scheme:yeeoriginal})  & \eqref{eq:yeecolloc1fourier}--\eqref{eq:yeecolloc3fourier}\\
  Yee (\ref{scheme:yeecollocated}) & \eqref{eq:yeecolloc1fourier}--\eqref{eq:yeecolloc3fourier}\\
  explicit Yee &   \\
  \hline
  Yee extended (\ref{scheme:yeecollocatedextended}) & \eqref{eq:yeeextended1fourier}--\eqref{eq:yeeextended3fourier}\\ 
  Yee extended (\ref{scheme:yeeextended}) & \eqref{eq:yeeextended3d1fourier}--\eqref{eq:yeeextended3d3fourier} \\
  central (\ref{scheme:central}) & \\
  central extended (\ref{scheme:centralextended}) & \eqref{eq:centralextended1fourier}--\eqref{eq:centralextended3fourier}\\
  \hline
  central extended $+$  & \multirow{2}{*}{--}\\
   pressureless (\ref{scheme:euler}) & 
 \end{tabular}
 \caption{Overview of the numerical methods and their discrete Fourier transforms, used for both stability and structure preservation analysis.}
 \label{tab:overviewfourier}
\end{table}

The paper is organized as follows: Section \ref{sec:yee} is an analysis of the Yee scheme, focusing particularly on interpretations as a collocated method. Section \ref{sec:yeeextensions} presents extensions of the Yee scheme which show enhanced stability for the Maxwell equations. Section \ref{sec:acoustics} shows the close relation of these results to numerical methods for linear acoustics, which are then extended to include advection in Section \ref{sec:euler}. This yields an all-speed numerical method for the full Euler equations.

\section{Review of leap-frog-type time integration methods}

The transfer of ideas from staggered-grid to collocated methods is possible because the low Mach number properties are rather related to the specific time integration customarily used for staggered-grid methods, than to the staggering itself. In its simplest form this particular time integration reads 
\begin{align}
 \frac{a^{n+1} - a^n}{\Delta t} &= f(b^n) &      a,b &\colon \mathbb R^+_0 \to \mathbb R \label{eq:triangintro1}\\
 \frac{b^{n+1} - b^n}{\Delta t} &= g(a^{n+1}) & f,g &\colon \mathbb R \to \mathbb R \text{ given} \label{eq:triangintro2}
\end{align}
which is a discretization of the system
\begin{align}
 \del_t \vecc{a}{b} &= \vecc{f(b)}{g(a)} \label{eq:examplesystem}
\end{align}
In the application to PDEs, of course, $f$ and $g$ will be replaced by appropriate spatial discretizations of differential operators.

In a way, it is a leap-frog method, but -- as is reviewed next -- in the context of PDEs other methods are also called leap-frog, such that here a different name is used to reduce confusion: \emph{sequential explicit}. Although the right-hand side of \eqref{eq:triangintro1}--\eqref{eq:triangintro2} formally is implicit, due to the special (``off-diagonal'') structure of the Jacobian of \eqref{eq:examplesystem} it can be implemented as an explicit method. The name ``sequential explicit'' derives from sequentially solving the equations and immediately using the new values of the variables as soon as they become available.

\subsection{Relation to semi-implicit methods}

The ability to rewrite sequential explicit methods as explicit methods is the main difference to ``semi-implicit'' methods (e.g. \cite{guerra86,guerra86a,gustafsson87}) which employ operator splitting and solve one of the operators fully implicitly. There, for example, $\del_t q = P_0 q + P_1 q$ is solved by
\begin{align}
 \frac{q^{n+1} - q^{n-1}}{2 \Delta t} = P_0 q^{n+1} + P_1 q^n
\end{align}
where $q$ is the vector of conserved quantities, and $P_0$ and $P_1$ are first order differential operators in space (or their discretizations), with $P_0$ stiff.
The same is true for the Crank-Nicolson scheme \cite{crank47}, also sometimes called leap-frog (e.g. in \cite{guerra86}). 

A time integration similar to \eqref{eq:triangintro1}--\eqref{eq:triangintro2} was used in \cite{dimarco17,boscarino19}:
\begin{align}
 \frac{a^{n+1} - a^n}{\Delta t} &= f(b^{n+1}) &      a,b &\colon \mathbb R^+_0 \to \mathbb R  \label{eq:leapfrogimplicit1}\\
 \frac{b^{n+1} - b^n}{\Delta t} &= g_1(a^{n+1}) + g_2(b^n)  & f,g_1,g_2 &\colon \mathbb R \to \mathbb R \text{ given} \label{eq:leapfrogimplicit2}
\end{align}
Observe that this is a truly implicit time integration, which cannot be rewritten as an explicit one. However, it can be reduced to just one implicit equation for $a^{n+1}$, such that the other can be updated explicitly. By inserting the second equation into the first one obtains
\begin{align}
 \frac{a^{n+1} - a^n}{\Delta t} &= f\Big(b^n + \Delta t (g_1(a^{n+1}) + g_2(b^n)) \Big) \label{eq:leapfrogimplicitinserted}
\end{align}
The discretization \eqref{eq:leapfrogimplicit1}--\eqref{eq:leapfrogimplicit2} therefore comes at the advantage that only a single implicit equation \eqref{eq:leapfrogimplicitinserted} needs to be solved, and \eqref{eq:leapfrogimplicit2} then amounts to an explicit update. The difference to the approach of using the time integration \eqref{eq:triangintro1}--\eqref{eq:triangintro2} is that the latter can be rewritten as an explicit update for \emph{all} the equations.

\subsection{Relation to symplectic and energy-conserving methods}

In the context of ODEs/dynamical systems, the time-stepping scheme \eqref{eq:triangintro1}--\eqref{eq:triangintro2}, in particular when $f$ is linear, is widely used and also carries the name ``symplectic Euler'' and is related to the Verlet method. This is because upon defining the Hamiltonian $H(a,b) = G(a) - F(b)$ with $F,G$ primitives of $f$ and $g$ (i.e. $F'=f$, $G'=g$), \eqref{eq:examplesystem} can be rewritten as
\begin{align}
 \del_t a &= - \frac{\del H}{\del b}\\
 \del_t b &= \frac{\del H}{\del a}
\end{align}
i.e. as a Hamiltonian system. Observe that the Hamiltonian is conserved:
\begin{align}
 \del_t H(a,b) = \frac{\del H}{\del a} \del_t a + \frac{\del H}{\del b} \del_t b = - \frac{\del H}{\del a} \frac{\del H}{\del b} + \frac{\del H}{\del b} \frac{\del H}{\del a} = 0
\end{align}
While symplectic Euler in general is an implicit method, for separable Hamiltonians such as the one defined above it is implementable explicitly, i.e. as \eqref{eq:triangintro1}--\eqref{eq:triangintro2}. There also exist other, equivalent ways of writing the symplectic Euler method. However, the author is unaware of its usage for compressible flow problems, or the low Mach number limit, where no Hamiltonian structure is available. 

For linear $f$ and $g$ (i.e. $f'=\const$, $g'=\const$), $H(a,b) = g' \frac{a^2}{2} - f' \frac{b^2}{2}$, and one can easily show by explicit computation that the following discrete Hamiltonian is conserved in time:
\begin{align}
 H_\text{discrete} = g' \frac{a^n a^{n+1}}{2} - f' \frac{(b^n)^2}{2}
\end{align}

The preservation of an (albeit modified) Hamiltonian (i.e. of energy) is important for long-term simulations of e.g. the Maxwell equations, where stabilization via upwinding is very quickly degrading the numerical results. Besides time-implicit methods, leap-frog-type time discretizations as in Yee's scheme together with appropriate spatial discretization (centered fluxes/derivatives) have been used to achieve energy-conservation (e.g.  \cite{remaki99,rodrigue01,fezoui05}) for Maxwell's equations. In e.g. \cite{cheng14a,cheng14b,cheng14}, the Maxwell system is coupled with further equations (e.g. the Vlasov equation), with the sequential explicit time integration used for Maxwell's equations, while e.g. the Vlasov equation is updated explicitly. It has been emphasized that choosing a central discretization in space for the former yields energy conservation (and upwind flux does not), while for the latter, central derivatives yield instability and therefore upwind flux is required.

However, for the eventual application to low Mach flows, it is rather not the aspect of symplecticity/energy-conservation that is most interesting, but the fact that sequential explicit numerical methods are dissipation-free. This property conceptually does not rely on existence of a Hamiltonian structure.

\subsection{Non-dissipativity}

A well-known property of the sequential explicit method \eqref{eq:triangintro1}--\eqref{eq:triangintro2} for linear $f,g$ is that it is non-dissipative. First of all, one notes that for $f'g' < 0$, the system \eqref{eq:examplesystem} has oscillatory solutions. In the linear case, \eqref{eq:triangintro1}--\eqref{eq:triangintro2} can be rewritten as
\begin{align}
 \vecc{a}{b}^{n+1} &= \left( \begin{array}{cc} 1 & \Delta t f' \\ \Delta t g' & 1 + \Delta t^2 f'g' \end{array} \right) \vecc{a}{b}^n =: A \vecc{a}{b}^n
\end{align}
If the time-dependence of $(a,b)^\text{T}$ is of the form $\exp(-\ii \omega t) = \exp(-\ii \omega \Delta t n)$ for some $\omega$, then $\exp(-\ii \omega \Delta t)$ must be an eigenvalue of $A$. Moreover, if $\omega$ is always real (i.e. if all eigenvalues of $A$ are unit complex numbers), then $(a,b)^\text{T}$ merely oscillates in time, while an imaginary part of $\omega$ will cause $(a,b)^\text{T}$ to decay/grow. A growth would be referred to as instability, while a decay would be associated to numerical dissipation. In absence of both, a numerical method is called non-dissipative.

\begin{theorem}
 Assume that $f,g$ in \eqref{eq:triangintro1}--\eqref{eq:triangintro2} are linear, that $f'g' <0$ and that $\Delta t  < \frac{2}{\sqrt{-f'g'}}$. Then \eqref{eq:triangintro1}--\eqref{eq:triangintro2} is non-dissipative.
\end{theorem}
\begin{proof}
The eigenvalues $\lambda$ of $A$ fulfill
\begin{align}
 0 &= (1 - \lambda)(1 + \Delta t^2 f'g' - \lambda) - \Delta t^2  f'g' = 1 - \lambda (2 + \Delta t^2 f'g') + \lambda^2
\end{align}
Define $y := - \frac{\Delta t^2 f'g'}{4}$ and by assumption $y \in (0,1)$, such that
\begin{align}
 \lambda &= 1 + \frac{\Delta t^2 f'g'}{2} \pm \sqrt{\Delta t^2 f'g' \left(  1 + \frac{\Delta t^2 f'g'}{4}  \right ) } 
= \left( \sqrt{1 -y }   \pm \ii \sqrt{y}     \right )^2 \\
|\lambda|^2 &= \left( 1 -y  +  y     \right )^2 = 1
\end{align}

\end{proof}

For Maxwell's equations (leaving the spatial derivatives continuous for the moment) one finds something similar\footnote{Boldface symbols are reserved for vectors with as many components as there are spatial dimensions. Indices never denote derivatives.}:
 \begin{align}
  \frac{\vec B^{n+1} - \vec B^n}{\Delta t} &= - \curl \vec E^n \\
  \frac{\vec E^{n+1} - \vec E^n}{\Delta t} &= \curl \vec B^{n+1} 
 \end{align}
 Considering a Fourier mode in space $\vecc{\vec B^n(x)}{\vec E^n(x)} = \vecc{\hat{\vec B}^n}{\hat {\vec E}^n} \exp(\ii \vec k \cdot \vec x)$ one rewrites
 \begin{align}
  \vecc{\hat{\vec B}^{n+1}}{\hat {\vec E}^{n+1}} &= \left(\begin{array}{cc} \id & - \Delta t \ii K \\ \Delta \ii K & \id + \Delta t^2 K^2  \end{array}\right)\vecc{\hat{\vec B}^n}{\hat {\vec E}^n} \label{eq:maxwellsystemfourierspace}
 \end{align}
 with $K = \left( \begin{array}{ccc}\phantom{a}\\\vec k \times \vec e_x &\vec k \times \vec e_y&\vec k \times \vec e_z\\\phantom{m}\end{array}\right) = \left(\begin{array}{ccc} 0 & -k_z & k_y \\ k_z &0&-k_x \\ -k_y & k_x & 0  \end{array}\right)$. Using the same analysis as before one verifies that, as long as $\Delta t |\vec k| < 2$, all eigenvalues $\lambda$ of the matrix in \eqref{eq:maxwellsystemfourierspace} fulfill $|\lambda| = 1$.

This property of non-dissipativity has been the focus of \cite{thomas93,roe98,kim04}. They compare, for linear advection $\del_t q + c \del_x q = 0$, the standard leap-frog
\begin{align}
 \frac{q_j^{n+1} -  q_j^{n-1}}{2 \Delta t} + c \frac{q_{j+1}^n - q_{j-1}^n}{2 \Delta x} &= 0 \label{eq:leapfrogadvection}
\end{align}
to what they call upwind-leap-frog ($c> 0$)
\begin{align}
 \frac{q_j^{n+1} -  q_j^{n} + q_{j-1}^n - q_{j-1}^{n-1}}{2 \Delta t} + c \frac{q_{j}^n - q_{j-1}^n}{ \Delta x} &= 0 \label{eq:upwindleapfrogadvection}
\end{align}
Observe that, if $q_j^n = \exp(-\ii \omega \Delta t n + \ii k \Delta x j)$, then \eqref{eq:leapfrogadvection} implies
\begin{align}
 \frac{\lambda -  1/\lambda}{2 \Delta t} + c \frac{\ii \sin (k \Delta x)}{\Delta x} &= 0 
\end{align}
\begin{align}
 \lambda = - \frac{c\Delta t}{\Delta x} \ii \sin (k \Delta x) \pm \sqrt{1 -\left( \frac{c\Delta t}{\Delta x} \right )^2 \sin^2 (k \Delta x)}
\end{align}
having defined again $\lambda := \exp(-\ii \omega \Delta t)$. Assuming $0 < \frac{c\Delta t}{\Delta x} < 1$, 
\begin{align}
 |\lambda|^2 &= 1
\end{align}
Thus, indeed, \eqref{eq:leapfrogadvection} is non-dissipative and the numerical error manifests itself as erroneous dispersion only. As can be seen from the proof, this property depends on the spatial discretization chosen. Analogous results are valid for \eqref{eq:upwindleapfrogadvection}. Whereas for systems in one spatial dimension these methods can easily be used for each characteristic, in multiple spatial dimensions, in order to apply them to bicharacteristic relations the authors of \cite{thomas93,roe98} find themselves forced to double the number of pressure variables. The continuum of bicharacteristics is replaced by a finite subset which is made compatible with a particular arrangement of the degrees of freedom in the grid. The resulting methods therefore are very different from the ones suggested here. In particular, the methods proposed in this work do not rely on a characteristic or bicharacteristic decomposition.

\section{A collocated interpretation of the Yee scheme} \label{sec:yee}

Consider the Maxwell equations 
 \begin{align}
  \del_t \vec B &= - \curl \vec E \\
  \del_t \vec E &= \curl \vec B \\
  &\phantom{m} \nonumber\\
  \div \vec E &= 0 \label{eq:maxwellinv1}\\
  \div \vec B &= 0 \label{eq:maxwellinv2}
 \end{align}
 and in particular their form in the transverse-magnetic case in two spatial dimensions:
 \begin{align}
  \del_t B^z &= -\Big( \del_x E^y - \del_y E^x \Big ) \label{eq:maxwell1}\\
  \del_t E^x &= \del_y B^z \\
  \del_t E^y &= - \del_x B^z \label{eq:maxwell3} 
 \end{align}

 The constraints \eqref{eq:maxwellinv1}--\eqref{eq:maxwellinv2} are fulfilled for all times, if they are filfilled at initial time, i.e. they are so-called \emph{involutions}. In two dimensions they amount to $\del_x E^x + \del_y E^y = 0$ and $\del_z B^z = 0$.

\subsection{Review of the original Yee scheme}

The original formulation of the scheme from \cite{yee66} reads

\begin{figure}
 \centering
 \includegraphics[width=0.2\textwidth]{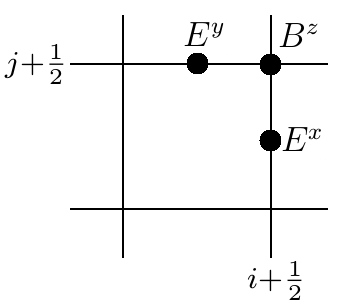} \qquad\qquad\qquad \includegraphics[width=0.2\textwidth]{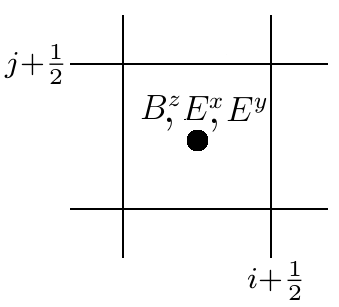}
 \qquad\qquad\qquad \includegraphics[width=0.2\textwidth]{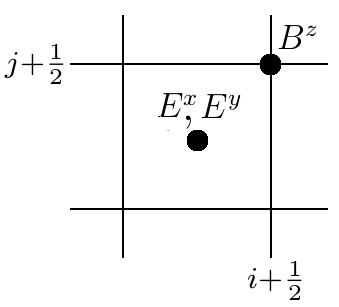}
 \caption{Degrees of freedom of Scheme \ref{scheme:yeeoriginal} (\emph{left}), of its collocated version (Scheme \ref{scheme:yeecollocated}, \emph{center}), and of Scheme \ref{scheme:yeeextended} (\emph{right}).}
 \label{fig:yee}
 \end{figure}

\setcounter{scheme}{24}
\begin{scheme}[Yee, 1966] \label{scheme:yeeoriginal}
\begin{align}
 &\frac{(B^z)^{n+\frac12}_{i+\frac12,j+\frac12} - (B^z)^{n-\frac12}_{i+\frac12,j+\frac12}}{\Delta t} = \label{eq:yeeorig1} \\\nonumber&\phantom{mmm}-\left( \frac{(E^y)^n_{i+1,j+\frac12}  - (E^y)^n_{i,j+\frac12}}{\Delta x} - \frac{(E^x)^n_{i+\frac12,j+1}  - (E^x)^n_{i+\frac12,j}}{\Delta y}  \right ) \\
 &\frac{(E^x)^{n+1}_{i+\frac12,j} - (E^x)^{n}_{i+\frac12,j}}{\Delta t} = \frac{(B^z)^{n+\frac12}_{i+\frac12,j+\frac12} - (B^z)^{n+\frac12}_{i+\frac12,j-\frac12}}{\Delta y}\\
 &\frac{(E^y)^{n+1}_{i,j+\frac12} - (E^y)^{n}_{i,j+\frac12}}{\Delta t} = -\frac{(B^z)^{n+\frac12}_{i+\frac12,j+\frac12} - (B^z)^{n+\frac12}_{i-\frac12,j+\frac12}}{\Delta x} \label{eq:yeeorig3}
\end{align}
\end{scheme} \setcounter{scheme}{0}

Here, the magnetic field $B^z$ is associated to locations $(i+\frac12,j+\frac12)$ (nodes), while the components of the electric field are associated with $(i+\frac12,j)$, $(i,j+\frac12)$ (edges) (see Figure \ref{fig:yee}). On Cartesian grids, therefore, there is one magnetic field $B^z$ and one of each components $E^x, E^y$ of the electrical field \emph{per cell}.

As has been derived in \cite{taflove75}, the stability condition for this scheme reads
\begin{align}
 \Delta t < \frac{1}{\sqrt{\frac{1}{\Delta x^2} + \frac{1}{\Delta y^2}}}
\end{align}
i.e. for equidistant two-dimensional grids $\Delta y = \Delta x$ the maximum CFL number is
\begin{align}
 \text{CFL}_\text{max} =  \frac{\Delta t_\text{max}}{\Delta x} = \frac{1}{\sqrt{2 }}
\end{align}

\subsection{A collocated interpretation of the Yee scheme}

Counted per cell, there is one of each variables $B^z, E^x, E^y$ stored as a degree of freedom. The association of some of them with an edge, and others not, is a matter of interpretation. This interpretation becomes important when the discrete value is compared to the exact solution. A discrete variable might, for example, be a higher order approximation to the exact value at one location than at another. Besides such comparison, however, the algorithm is performing the same algebraic manipulations irrespective of how the variables are named. Changing their interpretation (and accepting a possible reduction of the formal order of accuracy of the method) can then give what shall be called a ``collocated interpretation'' of the same scheme. The renaming amounts to moving the discrete degrees of freedom without changing the way they are updated (see Figure \ref{fig:yee}). Once the collocated version is understood, higher order of accuracy can be restored (see Section \ref{ssec:central}), yielding then a genuinely new algorithm and not just a reinterpretation.

Renaming 
\begin{align}
(B^z)^{n+\frac12}_{i+\frac12,j+\frac12} &\mapsto (B^z)^{n+\frac12}_{ij} & E^x_{i+\frac12,j} &\mapsto E^x_{ij} & E^y_{i,j+\frac12} &\mapsto E^y_{ij} \label{eq:yeerenaming}
\end{align}
yields the scheme
\begin{align}
 \frac{(B^z)^{n+\frac12}_{ij} - (B^z)^{n-\frac12}_{ij}}{\Delta t} &= -\left( \frac{(E^y)^n_{i+1,j}  - (E^y)^n_{ij}}{\Delta x} - \frac{(E^x)^n_{i,j+1}  - (E^x)^n_{ij}}{\Delta y}  \right ) \\
 \frac{(E^x)^{n+1}_{ij} - (E^x)^{n}_{ij}}{\Delta t} &= \frac{(B^z)^{n+\frac12}_{ij} - (B^z)^{n+\frac12}_{i,j-1}}{\Delta y}\\
 \frac{(E^y)^{n+1}_{ij} - (E^y)^{n}_{ij}}{\Delta t} &= -\frac{(B^z)^{n+\frac12}_{ij} - (B^z)^{n+\frac12}_{i-1,j}}{\Delta x}  
\end{align}

In fact, the same renaming can be performed with respect to the time index, yielding

\begin{scheme}[collocated Yee] \label{scheme:yeecollocated}
\begin{align}
 \frac{(B^z)^{n+1}_{ij} - (B^z)^{n}_{ij}}{\Delta t} &= -\left( \frac{(E^y)^n_{i+1,j}  - (E^y)^n_{ij}}{\Delta x} - \frac{(E^x)^n_{i,j+1}  - (E^x)^n_{ij}}{\Delta y}  \right ) \label{eq:yeecolloc1}\\
 \frac{(E^x)^{n+1}_{ij} - (E^x)^{n}_{ij}}{\Delta t} &= \frac{(B^z)^{n+1}_{ij} - (B^z)^{n+1}_{i,j-1}}{\Delta y} \label{eq:yeecolloc2}\\
 \frac{(E^y)^{n+1}_{ij} - (E^y)^{n}_{ij}}{\Delta t} &= -\frac{(B^z)^{n+1}_{ij} - (B^z)^{n+1}_{i-1,j}}{\Delta x}  \label{eq:yeecolloc3}
\end{align}
\end{scheme}

\begin{remark} \label{rem:otherdirectionshiftyee}
Note that renaming $(B^z)^{n+\frac12}_{i+\frac12,j+\frac12} \mapsto (B^z)^{n+\frac12}_{ij}$, $E^x_{i+\frac12,j+1} \mapsto E^x_{ij}$ and $E^y_{i+1,j+\frac12} \mapsto E^y_{ij}$ yields
\begin{align}
 \frac{(B^z)^{n+1}_{ij} - (B^z)^{n}_{ij}}{\Delta t} &= -\left( \frac{(E^y)^n_{ij}  - (E^y)^n_{i-1,j}}{\Delta x} - \frac{(E^x)^n_{ij}  - (E^x)^n_{i,j-1}}{\Delta y}  \right )\\
 \frac{(E^x)^{n+1}_{ij} - (E^x)^{n}_{ij}}{\Delta t} &= \frac{(B^z)^{n+1}_{i,j+1} - (B^z)^{n+1}_{ij}}{\Delta y}\\
 \frac{(E^y)^{n+1}_{ij} - (E^y)^{n}_{ij}}{\Delta t} &= -\frac{(B^z)^{n+1}_{i+1,j} - (B^z)^{n+1}_{ij}}{\Delta x}
\end{align}
instead of \eqref{eq:yeecolloc1}--\eqref{eq:yeecolloc3}, such that having forward finite differences in the first equation, and backward differences in the second does not have any fundamental meaning, and could be the other way around.
\end{remark}

\subsection{Reinterpretation as a fully explicit method}

 The numerical method \eqref{eq:yeecolloc1}--\eqref{eq:yeecolloc3} is not really implicit in time, even if the right hand side involves values at time step $n+1$. Inserting \eqref{eq:yeecolloc1} in \eqref{eq:yeecolloc2}--\eqref{eq:yeecolloc3} yields its fully explicit form:
 
 \begin{scheme} \label{scheme:yeecollocatedexplicit}
 \begin{align}
 \frac{(B^z)^{n+1}_{ij} - (B^z)^{n}_{ij}}{\Delta t} &= -\left( \frac{(E^y)^n_{i+1,j}  - (E^y)^n_{ij}}{\Delta x} - \frac{(E^x)^n_{i,j+1}  - (E^x)^n_{ij}}{\Delta y}  \right )\label{eq:yeecollocexpl1} \\\nonumber \phantom{m}\\
 \frac{(E^x)^{n+1}_{ij} - (E^x)^{n}_{ij}}{\Delta t} &= \frac{
 (B^z)^{n}_{ij} - (B^z)^{n}_{i,j-1}}{\Delta y} \label{eq:yeecollocexpl2}
 \\&\!\!\!\!\!\!\!\!\!\!\!\!\!\!\!\!\!\!\!\!\!\!\!\!\nonumber- \frac{\Delta t}{\Delta y} \left( \frac{(E^y)^n_{i+1,j} - (E^y)^n_{ij} - (E^y)^n_{i+1,j-1} + (E^y)^n_{i,j-1}}{\Delta x} \right . \\ \nonumber &\phantom{mmmmmmmm} \left . - \frac{(E^x)^n_{i,j+1}  -2 (E^x)^n_{ij}  + (E^x)^n_{i,j-1}}{\Delta y}  \right )\\ \nonumber\phantom{m}\\
 \frac{(E^y)^{n+1}_{ij} - (E^y)^{n}_{ij}}{\Delta t} &= -\frac{(B^z)^{n}_{ij} - (B^z)^{n}_{i-1,j}}{\Delta x} \label{eq:yeecollocexpl3}
  \\&\!\!\!\!\!\!\!\!\!\!\!\!\!\!\!\!\!\!\!\!\!\!\!\!\nonumber+ \frac{\Delta t}{\Delta x} \left( \frac{(E^y)^n_{i+1,j}  -2 (E^y)^n_{ij}  + (E^y)^n_{i-1,j}}{\Delta x} 
  \right . \\ \nonumber &\phantom{mmm} \left . - \frac{(E^x)^n_{i,j+1}  -(E^x)^n_{ij}  - (E^x)^n_{i-1,j+1} + (E^x)^n_{i-1,j}}{\Delta y} \right )
\end{align}
\end{scheme}

Observe that the right-hand side of \eqref{eq:yeecollocexpl2}--\eqref{eq:yeecollocexpl3} now contains \emph{second} derivatives in space. This is reminiscent of the second derivatives appearing in the dimensionally split upwind method for \eqref{eq:maxwell1}--\eqref{eq:maxwell3}
 \begin{align}
 \frac{(B^z)^{n+1}_{ij} - (B^z)^{n}_{ij}}{\Delta t} &= \label{eq:collocupwind1}  \\ \nonumber & \!\!\!\!\!\!\!\!\!\!\!\!\!\!\!\!\!\!\!\!\!\!\!\!\!\!\!\!\!\!\!\!\!\!\! -\left( \frac{(E^y)^n_{i+1,j}  - (E^y)^n_{i-1,j}}{2\Delta x} - \frac{(E^x)^n_{i,j+1}  - (E^x)^n_{i,j-1}}{2\Delta y}  \right ) 
 \\&\nonumber \!\!\!\!\!\!\!\!\!\!\!\!\!\!\!\!\!\!\!\!\!\!\!\!\!\!\!\!\!\!\!\!\!\!\!+ \frac12 \frac{(B^z)^n_{i+1,j} - 2 (B^z)^n_{ij} + (B^z)^n_{i-1,j}}{\Delta  x} + \frac12 \frac{(B^z)^n_{i,j+1} - 2 (B^z)^n_{ij} + (B^z)^n_{i,j-1}}{\Delta y} \\ \nonumber \phantom{m}\\
 \frac{(E^x)^{n+1}_{ij} - (E^x)^{n}_{ij}}{\Delta t} &= \frac{
 (B^z)^{n}_{i,j+1} - (B^z)^{n}_{i,j-1}}{2\Delta y} \\\nonumber & + \frac12  \frac{(E^x)^n_{i,j+1}  -2 (E^x)^n_{ij}  + (E^x)^n_{i,j-1}}{\Delta y} \\ \nonumber \phantom{m}\\
 \frac{(E^y)^{n+1}_{ij} - (E^y)^{n}_{ij}}{\Delta t} &= -\frac{(B^z)^{n}_{i+1,j} - (B^z)^{n}_{i-1,j}}{2\Delta x}\label{eq:collocupwind3}
 \\\nonumber& + \frac12  \frac{(E^y)^n_{i+1,j}  -2 (E^y)^n_{ij}  + (E^y)^n_{i-1,j}}{\Delta x} 
\end{align}

In the context of the Maxwell equations, an important question is whether the numerical method is involution preserving, i.e. whether a discretization of the involution is kept stationary. It has been shown in \cite{barsukow17a} that a linear involution preserving numerical method is also stationarity preserving, i.e. its stationary states are a discretization of \emph{all} the stationary states of the PDE. A special property of linear acoustics, discussed later, is that its low Mach number limit is equivalent to the long time limit. Low Mach number compliance for linear acoustics is then equivalent to the property of stationarity preservation, because von Neumann stable numerical methods dissipate away any Fourier mode that is not stationary. This also is described in more detail in \cite{barsukow17a}. Thus, when obtaining an involution preserving method for the Maxwell equations, stationarity preservation comes for free. Once the methods are applied to linear acoustics, this then implies their low Mach number compliance.

It is known that the upwind method is not involution preserving. It is thus instructive to also compare \eqref{eq:yeecollocexpl1}--\eqref{eq:yeecollocexpl3} to the involution-preserving, truly multi-dimensional method from \cite{barsukow17a} (following \cite{morton01,jeltsch06,mishra09preprint}): \newpage
 \begin{align}
 \frac{(B^z)^{n+1}_{ij} - (B^z)^{n}_{ij}}{\Delta t} &= \label{eq:statpresmaxwell1}
 \\\nonumber &\!\!\!\!\!\!\!\!\!\!\!\!\!\!\!\!\!\!\!\!\!\!\!\!-\left( \frac{\langle(E^y)^n_{i+1}\rangle_j  - \langle(E^y)^n_{i-1}\rangle_j}{2\Delta x} - \frac{\langle(E^x)^n_{\cdot,j+1}\rangle_i  - \langle(E^x)^n_{\cdot,j-1}\rangle_i}{2\Delta y}  \right )  
 \\&\nonumber \!\!\!\!\!\!\!\!\!\!\!\!\!\!\!\!\!\!\!\!\!\!\!\!+ \frac12 \frac{\langle(B^z)^n_{i+1}\rangle_j - 2 \langle(B^z)^n_{i}\rangle_j + \langle(B^z)^n_{i-1}\rangle_j}{\Delta  x} \\&\nonumber \!\!\!\!\!\!\!\!\!\!\!\!\!\!\!\!\!\!\!\!\!\!\!\!+ \frac12 \frac{\langle(B^z)^n_{\cdot,j+1}\rangle_i - 2 \langle(B^z)^n_{\cdot,j}\rangle_i + \langle(B^z)^n_{\cdot,j-1}\rangle_i}{\Delta y} \\ \nonumber \phantom{m}\\
 \frac{(E^x)^{n+1}_{ij} - (E^x)^{n}_{ij}}{\Delta t} &= \frac{\langle(B^z)^{n}_{\cdot,j+1}\rangle_i - \langle(B^z)^{n}_{\cdot,j-1}\rangle_i}{2\Delta y} 
 \\&\nonumber \!\!\!\!\!\!\!\!\!\!\!\!\!\!\!\!\!\!\!\!\!\!\!\!+ \frac12 \left( -\frac{(E^y)^n_{i+1,j} - (E^y)^n_{ij} - (E^y)^n_{i+1,j-1} + (E^y)^n_{i,j-1}}{\Delta x} \right. \\&\nonumber \!\!\!\!\!\!\!\! \left. + \frac{\langle(E^x)^n_{\cdot,j+1}\rangle_i  -2 \langle(E^x)^n_{\cdot,j}\rangle_i  + \langle(E^x)^n_{\cdot,j-1}\rangle_i}{\Delta y} \right )
 \\ \nonumber \phantom{m}\\
 \frac{(E^y)^{n+1}_{ij} - (E^y)^{n}_{ij}}{\Delta t} &= -\frac{\langle (B^z)^{n}_{i+1} \rangle_j - \langle(B^z)^{n}_{i-1} \rangle_j}{2\Delta x} \label{eq:statpresmaxwell3}
 \\&\nonumber \!\!\!\!\!\!\!\!\!\!\!\!\!\!\!\!\!\!\!\!\!\!\!\!+ \frac12 \left( \frac{\langle(E^y)^n_{i+1}\rangle_j  -2\langle (E^y)^n_{i}\rangle_j  + \langle(E^y)^n_{i-1}\rangle_j}{\Delta x} \right. \\&\nonumber \!\!\!\!\!\!\!\! \left.- \frac{(E^x)^n_{i,j+1}  -(E^x)^n_{ij}  - (E^x)^n_{i-1,j+1} + (E^x)^n_{i-1,j}}{\Delta y} \right )
\end{align}
where 
\begin{align}
\langle q_i \rangle_j &:= \frac{q_{i,j+1} + 2 q_{ij} + q_{i,j-1}}{4} & \langle q_{\cdot,j} \rangle_i &:= \frac{q_{i+1,j} + 2 q_{ij} + q_{i-1,j}}{4} \label{eq:averagingoperator}
\end{align}
are averaging operators.

The first terms on the right-hand side of \eqref{eq:yeecollocexpl1}--\eqref{eq:yeecollocexpl3} are a one-sided discretization of the PDE, while the other two schemes employ central differences there. Although in one spatial dimension, Maxwell's equations \eqref{eq:maxwell1}--\eqref{eq:maxwell3} can be diagonalized and thus decoupled into a left-going and a right-going advection equations, there is no physical reason behind taking $E^x$ and $E^y$ from the right / top and $B^z$ from the left and from below. As has been seen in Remark \ref{rem:otherdirectionshiftyee}, this is purely an artifact of the renaming, when the initially staggered method was reinterpreted as a collocated one. This asymmetry is addressed further in Section \ref{ssec:central}.
 
As is well known, scheme \eqref{eq:yeecolloc1}--\eqref{eq:yeecolloc3} does not have diffusion in the sense that the absolute values of all the Fourier modes are stationary (all the eigenvalues in the von Neumann analysis have absolute value 1). This then is also true for \eqref{eq:yeecollocexpl1}--\eqref{eq:yeecollocexpl3}. Thus, contrary to \eqref{eq:collocupwind1}--\eqref{eq:collocupwind3}, the second derivatives appearing in \eqref{eq:yeecollocexpl1}--\eqref{eq:yeecollocexpl3} are not a diffusion.

\subsection{Sequential explicit time integration and structure preservation}

For involution preservation/stationarity preservation, one of the eigenvalues needs to be exactly 1, not just in absolute value. For \eqref{eq:yeecollocexpl1}--\eqref{eq:yeecollocexpl3}, this is indeed the case, a consequence of the following Theorem,

\begin{theorem} \label{thm:statpresfortriangimpl} 
 The numerical method
 \begin{align}
 \frac{(B^z)_{ij}^{n+1} - (B^z)_{ij}^{n}}{\Delta t} &= -\left( \Big(\mathscr D_x (E^y)^n\Big)_{ij} - \Big(\mathscr D_y (E^x)^n\Big)_{ij}   \right )  \\
 \frac{(E^x)_{ij}^{n+1} - (E^x)_{ij}^{n}}{\Delta t} &= \Big(\mathscr D'_y (B^z)^{n+1}\Big)_{ij}  \\
 \frac{(E^y)_{ij}^{n+1} - (E^y)_{ij}^{n}}{\Delta t} &= -\Big(\mathscr D'_x (B^z)^{n+1}\Big)_{ij}   
\end{align}
 for the Maxwell equations \eqref{eq:maxwell1}--\eqref{eq:maxwell3} is stationarity preserving / involution preserving for any choice of linear finite difference operators $\mathscr D_x, \mathscr D'_x, \mathscr D_y, \mathscr D'_y$. 
\end{theorem}

This result asserts involution preservation for numerical methods that employ sequential explicit time integration, and is proven in the Appendix (Section \ref{app:stability}). The terms appearing in \eqref{eq:yeecollocexpl1}--\eqref{eq:yeecollocexpl3} and in the involution preserving method \eqref{eq:statpresmaxwell1}--\eqref{eq:statpresmaxwell3} are very similar, and the most striking difference to the upwind scheme (which is not involution preserving) is the presence of terms discretizing cross-derivatives $\del_x \del_y E^x$ and $\del_x \del_y E^y$. In \eqref{eq:yeecollocexpl2}--\eqref{eq:yeecollocexpl3} they are generated by the leap-frog manner of time integration \eqref{eq:triangintro1}--\eqref{eq:triangintro2}. What is different is the $\Delta t$-dependent prefactor and the absence of second derivatives in the first equation \eqref{eq:yeecollocexpl1}.

For linear systems, involution preservation and stationarity preservation (and for acoustics, stationarity preservation and low Mach number compliance) are all the same (see \cite{barsukow17a}). In view of Theorem \ref{thm:statpresfortriangimpl} which implies involution preservation for a large class of sequential explicit methods, the conclusion therefore is that staggering the grid is not the decisive ingredient for low Mach number compliance, but rather this special time integration. At the discrete stationary state (stationarity preserving), the diffusion of scheme \eqref{eq:statpresmaxwell1}--\eqref{eq:statpresmaxwell3} vanishes, and so does the right-hand side of \eqref{eq:yeecolloc1}--\eqref{eq:yeecolloc2}. As will be seen below, it is possible to modify the Yee scheme such that it actually maintains the same discrete stationary state (and the same discrete involution) as \eqref{eq:statpresmaxwell1}--\eqref{eq:statpresmaxwell3}. The schemes differ, however, with respect to the behaviour of non-stationary solutions: while \eqref{eq:statpresmaxwell1}--\eqref{eq:statpresmaxwell3} adds diffusion, \eqref{eq:yeecolloc1}--\eqref{eq:yeecolloc2} is purely dispersive/non-dissipative.

Having identified the sequential explicit time integration as the relevant strategy to achieve involution preservation (for Maxwell's equations), or low Mach number compliance (for linear acoustics, and later for Euler), it is worth commenting on a few more of its properties.

Consider the linear version of \eqref{eq:triangintro1}--\eqref{eq:triangintro2}
\begin{align}
 \frac{a^{n+1} - a^n}{\Delta t} &= F b^n     \label{eq:triang1}\\
 \frac{b^{n+1} - b^n}{\Delta t} &= G a^{n+1} \label{eq:triang2}
\end{align}
which is a discretization of the system
\begin{align}
 \del_t \vecc{a}{b} &= \vecc{F b}{G a}
\end{align}
with $F, G \in \mathbb C$. Then, it is well-known (e.g. \cite{courant28,lax56}) that writing
\begin{align}
 a^{n+1} &= a^n + \Delta t F b^n = a^n + \Delta t F (b^{n-1} + \Delta t G a^{n})\\
 &= 2 a^n - a^{n-1} + \Delta t^2 F G a^{n} \label{eq:rewriteyeeaslaplace}
\end{align}
makes \eqref{eq:triang1}--\eqref{eq:triang2} equivalent to the most standard three-time-level discretization
\begin{align}
 \frac{a^{n+1} - 2 a^n + a^{n-1}}{\Delta t^2} &= FG a^n\\
 \frac{b^{n+1} - 2b^{n} + b^{n-1}}{\Delta t^2} &= FG b^{n}
\end{align}
of the decoupled second-order system
\begin{align}
 \del_t^2 \vecc{a}{b} & = \vecc{FG a}{FG b} 
\end{align}
In the context of Maxwell's and acoustic equations, where $G$ and $F$ are spatial derivatives/discretizations, the decoupled second order PDEs are (discretized) wave equations. In order to obtain in $FD$ the standard discretizations of the Laplacian, $F$ needs to be a forward, and $G$ a backward derivative, or vice versa. This can be seen in Scheme \ref{scheme:yeecollocated}.

\section{Modifications of the Yee scheme with improved stability} \label{sec:yeeextensions}

\subsection{A sequential explicit method with different staggering} \label{ssec:yeeextended}

Having established the way in which the collocated version of the Yee scheme achieves structure preservation, and in particular having identified the close relation to the stationarity preserving schemes of \cite{barsukow17a}, there are two natural questions that arise. Is there a way how the success of the Yee scheme for the Maxwell equations (and involution preservation) can be transferred to the acoustic equations (and stationarity preservation), thus ultimately yielding a low Mach number compliant scheme for the Euler equations? And are there aspects of stationarity preserving methods from \cite{barsukow17a} that can be used to improve the Yee scheme, be it its collocated version, or the original staggered one?

The answer to both questions is yes. The discussion of the former question forms the content of subsequent chapters, while the latter question is subject of the present section. Replacing the derivative
\begin{align}
 \frac{(E^y)^n_{i+1,j}  - (E^y)^n_{ij}}{\Delta x}
\end{align}
by the average
\begin{align}
 \frac12 \left( \frac{(E^y)^n_{i+1,j+1}  - (E^y)^n_{i,j+1}}{2\Delta x} + \frac{ (E^y)^n_{i+1,j}  - (E^y)^n_{ij}}{2\Delta x}   \right )
\end{align}
lets a vertex-based curl discretization appear. This is inspired by the involution preserving schemes in e.g. \cite{jeltsch06,mishra09preprint,roe15,barsukow17a} where derivatives are also ``extended'' to truly multi-dimensional finite differences through averaging in perpendicular directions.

\begin{scheme}[collocated Yee extended] \label{scheme:yeecollocatedextended}
\begin{align}
 \frac{(B^z)^{n+1}_{ij} - (B^z)^{n}_{ij}}{\Delta t} &= -\left( \frac{(E^y)^n_{i+1,j+1}  - (E^y)^n_{i,j+1} + (E^y)^n_{i+1,j}  - (E^y)^n_{ij}}{2\Delta x} \right . \label{eq:yeecollocextended1} \\\nonumber &\phantom{mmmm}- \left. \frac{(E^x)^n_{i+1,j+1}  - (E^x)^n_{i+1,j} + (E^x)^n_{i,j+1}  - (E^x)^n_{ij}}{2\Delta y}   \right ) \\
 \frac{(E^x)^{n+1}_{ij} - (E^x)^{n}_{ij}}{\Delta t} &= \frac{(B^z)^{n+1}_{ij} - (B^z)^{n+1}_{i,j-1} + (B^z)^{n+1}_{i-1,j} - (B^z)^{n+1}_{i-1,j-1}}{2 \Delta y} \\
 \frac{(E^y)^{n+1}_{ij} - (E^y)^{n}_{ij}}{\Delta t} &= - \frac{(B^z)^{n+1}_{ij}  - (B^z)^{n+1}_{i-1,j} + (B^z)^{n+1}_{i,j-1} - (B^z)^{n+1}_{i-1,j-1}}{2 \Delta x}  \label{eq:yeecollocextended3}
\end{align}
\end{scheme}

As is shown in Section \ref{ssec:fourieryeeextension}, in fact, the acoustic variant of Scheme \ref{scheme:yeecollocatedextended} keeps stationary precisely the same discrete divergence as e.g. the method in \cite{barsukow17a}. This extension also bears resemblance to ideas appearing in \cite{koh06}.

\begin{figure}
 \centering
 \includegraphics[width=0.3\textwidth]{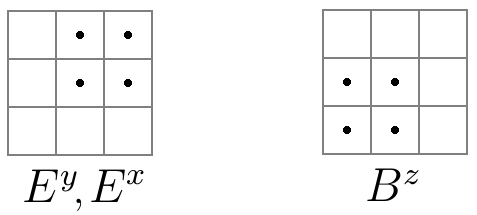}
 \caption{Stencils appearing in Scheme \ref{scheme:yeecollocatedextended}. \emph{Left}: Stencils of the discretizations of the curl of the electric field in the evolution of $B^z$. \emph{Right}: Stencils of the discretizations of the derivatives of the magnetic field in the evolution equations of $E^x, E^y$.}
 \label{fig:stencilsyeeextendedcoll}
\end{figure}

\begin{figure}
 \centering
 \includegraphics[width=0.7\textwidth]{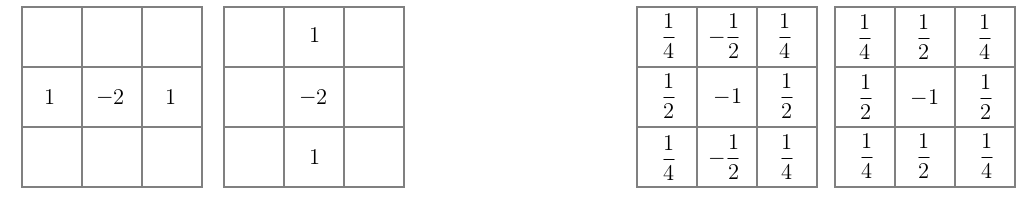}
 \caption{Stencils and weights of Laplacian discretizations. \emph{Left}: Standard Laplacian. \emph{Right}: Extended Laplacian appearing in a second-order-reformulation of Scheme \ref{scheme:yeecollocatedextended}.}
 \label{fig:stencillaplacian}
\end{figure}

Scheme \ref{scheme:yeecollocatedextended} can be read as a splitting of the discretization of the wave equations
\begin{align}
 \del_t^2 B^z &= \del_x^2 B^z + \del_y^2 B^z\\
 \del_t^2 E^x &= - \del_x \del_y E^y + \del_y^2 E^x\\
 \del_t^2 E^y &= \del_x^2 E^y - \del_x \del_y E^x
\end{align}
by analogy with \eqref{eq:rewriteyeeaslaplace}. In particular, Scheme \ref{scheme:yeecollocatedextended} is equivalent to {\footnotesize
\begin{align}
 \frac{(B_z)^{n+1}_{ij} - 2 (B_z)^{n}_{ij} + (B_z)^{n-1}_{ij}}{\Delta t^2 } &=  \frac14 \frac{(B_z)^{n}_{i+1,j+1} - 2 (B_z)^{n}_{i,j+1} + (B_z)^{n}_{i-1,j+1}}{\Delta x^2} \\\nonumber & \!\!\!\!\!\!\!\!\!\!\!\!\!\!\!\!\!\!\!\!\!\!\!\!\!\!\!\!\!\!\!\!\!\!\!\!\!\!\!\!\!\!\!\!\!\!\!\!\!\!\!\!\!\!\!\!\!\!\!\!\!\!\!
 + \frac12 \frac{(B_z)^{n}_{i+1,j} - 2 (B_z)^{n}_{ij} + (B_z)^{n}_{i-1,j}}{\Delta x^2} + \frac14 \frac{(B_z)^{n}_{i+1,j-1} - 2 (B_z)^{n}_{ij-1} + (B_z)^{n}_{i-1,j-1} }{\Delta x^2} \\
 &\nonumber \!\!\!\!\!\!\!\!\!\!\!\!\!\!\!\!\!\!\!\!\!\!\!\!\!\!\!\!\!\!\!\!\!\!\!\!\!\!\!\!\!\!\!\!\!\!\!\!\!\!\!\!\!\!\!\!\!\!\!\!\!\!\!+
 \frac14 \frac{  (B^z)^n_{i+1,j+1} - 2 (B^z)^n_{i+1,j} + (B^z)^n_{i+1,j-1} }{\Delta y^2} + \frac12 \frac{(B^z)^n_{i,j+1} - 2 (B^z)^n_{ij} + (B^z)^n_{i,j-1} }{\Delta y^2} \\\nonumber &\!\!\!\!\!\!\!\!\!\!\!\!\!\!\!\!\!\!\!\!\!\!\!\!\!\!\!\!\!\!\!\!\!\!\!\!\!\!\!\!\!\!\!\!\!\!\!\!\!\!\!\!\!\!\!\!\!\!\!\!\!\!\!
 + \frac14 \frac{ (B^z)^n_{i-1,j+1} - 2 (B^z)^n_{i-1,j} + (B^z)^n_{i-1,j-1}  }{\Delta y^2} \\
 \frac{(E^x)^{n+1}_{ij} - 2 (E^x)^n_{ij} + (E^x)^{n-1}_{ij} }{\Delta t^2} &= \\&\nonumber \!\!\!\!\!\!\!\!\!\!\!\!\!\!\!\!\!\!\!\!\!\!\!\!\!\!\!\!\!\!\!\!\!\!\!\!\!\!\!\!\!\!\!\!\!\!\!\!\!\!\!\!\!\!\!\!\!\!\!\!\!\!\! - \frac{ (E^y)^n_{i+1,j+1} - (E^y)^n_{i-1,j+1} - (E^y)^n_{i+1,j-1} + (E^y)^n_{i-1,j-1} }{4 \Delta x \Delta y} \\
 \nonumber &\!\!\!\!\!\!\!\!\!\!\!\!\!\!\!\!\!\!\!\!\!\!\!\!\!\!\!\!\!\!\!\!\!\!\!\!\!\!\!\!\!\!\!\!\!\!\!\!\!\!\!\!\!\!\!\!\!\!\!\!\!\!\!
 + \frac14 \frac{ (E^x)^n_{i+1,j+1} - 2 (E^x)^n_{i+1,j} + (E^x)^n_{i+1,j-1} }{\Delta y^2} + \frac12 \frac{ (E^x)^n_{i,j+1} - 2 (E^x)^n_{ij} + (E^x)^n_{i,j-1} }{\Delta y^2}\\
 \nonumber &\!\!\!\!\!\!\!\!\!\!\!\!\!\!\!\!\!\!\!\!\!\!\!\!\!\!\!\!\!\!\!\!\!\!\!\!\!\!\!\!\!\!\!\!\!\!\!\!\!\!\!\!\!\!\!\!\!\!\!\!\!\!\!
 + \frac14 \frac{ (E^x)^n_{i-1,j+1} - 2 (E^x)^n_{i-1,j} + (E^x)^n_{i-1,j-1} }{\Delta y^2}
\end{align}}
and similarly for the $E^y$-equation. Observe the appearance of a standard approximation of a second derivative in time and a natural multi-dimensionally extended Laplacian in space (see Figure \ref{fig:stencillaplacian}).

By Theorem \ref{thm:statpresfortriangimpl}, Scheme \ref{scheme:yeecollocatedextended} is involution preserving and by Corollary \ref{cor:yeeextendedstability} it shows improved stability up to $\mathrm{CFL} = 1$. 

Undoing the renaming \eqref{eq:yeerenaming} would yield
\begin{align}
 \frac{(B^z)^{n+1}_{i+\frac12,j+\frac12} - (B^z)^{n}_{i+\frac12,j+\frac12}}{\Delta t} &= \\&\nonumber \!\!\!\!\!\!\!\!\!\!\!\!\!\!\!\!\!\!\!\!\!\!\!\!\!\!\!\!\!\!\!\!\!\!\!\!\!\!\!\! -\left( \frac{(E^y)^n_{i+1,j+\frac32}  - (E^y)^n_{i,j+\frac32} + (E^y)^n_{i+1,j+\frac12}  - (E^y)^n_{i,j+\frac12}}{2\Delta x} \right . \\ \nonumber &\!\!\!\!\!\!\!\!\!\!\!\!\!\!\!\!\!\!\!\!\!\!\!\!\! - \left . \frac{(E^x)^n_{i+\frac32,j+1}  - (E^x)^n_{i+\frac32,j} + (E^x)^n_{i+\frac12,j+1}  - (E^x)^n_{i+\frac12,j}}{2\Delta y}   \right ) \\
 \frac{(E^x)^{n+1}_{i+\frac12,j} - (E^x)^{n}_{i+\frac12,j}}{\Delta t} &= \frac{(B^z)^{n+1}_{i+\frac12,j+\frac12} - (B^z)^{n+1}_{i+\frac12,j-\frac12} + (B^z)^{n+1}_{i-\frac12,j+\frac12} - (B^z)^{n+1}_{i-\frac12,j-\frac12}}{2 \Delta y} \\
 \frac{(E^y)^{n+1}_{i,j+\frac12} - (E^y)^{n}_{i,j+\frac12}}{\Delta t} &= - \frac{(B^z)^{n+1}_{i+\frac12,j+\frac12}  - (B^z)^{n+1}_{i-\frac12,j+\frac12} + (B^z)^{n+1}_{i+\frac12,j-\frac12} - (B^z)^{n+1}_{i-\frac12,j-\frac12}}{2 \Delta x} 
\end{align}
This is asymmetric, though. Undoing only the renaming of the magnetic field, however, restores symmetry: \newpage
\begin{align}
 &\frac{(B^z)^{n+1}_{i+\frac12,j+\frac12} - (B^z)^{n}_{i+\frac12,j+\frac12}}{\Delta t} = \\&\nonumber \phantom{mmmmmmmmmm} -\left( \frac{(E^y)^n_{i+1,j+1}  - (E^y)^n_{i,j+1} + (E^y)^n_{i+1,j}  - (E^y)^n_{ij}}{2\Delta x} \right . \\ \nonumber &\phantom{mmmmmmmmmmmm} - \left. \frac{(E^x)^n_{i+1,j+1}  - (E^x)^n_{i+1,j} + (E^x)^n_{i,j+1}  - (E^x)^n_{ij}}{2\Delta y}   \right ) \\
 &\frac{(E^x)^{n+1}_{ij} - (E^x)^{n}_{ij}}{\Delta t} = \frac{(B^z)^{n+1}_{i+\frac12,j+\frac12} - (B^z)^{n+1}_{i+\frac12,j-\frac12} + (B^z)^{n+1}_{i-\frac12,j+\frac12} - (B^z)^{n+1}_{i-\frac12,j-\frac12}}{2 \Delta y} \\
 &\frac{(E^y)^{n+1}_{ij} - (E^y)^{n}_{ij}}{\Delta t} = - \frac{(B^z)^{n+1}_{i+\frac12,j+\frac12}  - (B^z)^{n+1}_{i-\frac12,j+\frac12} + (B^z)^{n+1}_{i+\frac12,j-\frac12} - (B^z)^{n+1}_{i-\frac12,j-\frac12}}{2 \Delta x} 
\end{align}
This is an indication that the correct staggered-grid interpretation of the new scheme \ref{scheme:yeecollocatedextended} should use the magnetic field stored at the nodes $(i+\frac12,j+\frac12)$ (as in the original Yee scheme) and both components of the electric field stored at cell centers $(i,j)$. The scheme thus only employs the face-vertex dual meshes:

\begin{scheme}[Yee extended] \label{scheme:yeeextended}
\begin{align}
 &\frac{(B^z)^{n+1}_{i+\frac12,j+\frac12} - (B^z)^{n}_{i+\frac12,j+\frac12}}{\Delta t} = \label{eq:yeeextended1}\\\nonumber& \phantom{mmmmmmmmm} -\left( \frac{(E^y)^n_{i+1,j+1}  - (E^y)^n_{i,j+1} + (E^y)^n_{i+1,j}  - (E^y)^n_{ij}}{2\Delta x} \right .  \\\nonumber &\phantom{mmmmmmmmmmm}- \left.  \frac{(E^x)^n_{i+1,j+1}  - (E^x)^n_{i+1,j} + (E^x)^n_{i,j+1}  - (E^x)^n_{ij}}{2\Delta y}   \right ) \\
 &\frac{(E^x)^{n+1}_{ij} - (E^x)^{n}_{ij}}{\Delta t} = \frac{(B^z)^{n+1}_{i+\frac12,j+\frac12} - (B^z)^{n+1}_{i+\frac12,j-\frac12} + (B^z)^{n+1}_{i-\frac12,j+\frac12} - (B^z)^{n+1}_{i-\frac12,j-\frac12}}{2 \Delta y} \\
 &\frac{(E^y)^{n+1}_{ij} - (E^y)^{n}_{ij}}{\Delta t} = \label{eq:yeeextended3} - \frac{(B^z)^{n+1}_{i+\frac12,j+\frac12}  - (B^z)^{n+1}_{i-\frac12,j+\frac12} + (B^z)^{n+1}_{i+\frac12,j-\frac12} - (B^z)^{n+1}_{i-\frac12,j-\frac12}}{2 \Delta x}  
\end{align}
\end{scheme}

The scheme now is ``symmetric'': the vertex-based $B^z$-field is updated using a vertex-based curl and the cell-based $E$-fields are updated using cell-based gradients. As this scheme differs from Scheme \ref{scheme:yeecollocatedextended} only by renaming, it is also stable for $\text{CFL} < 1$.

The original Yee scheme in \cite{yee66} is formulated for the three-dimensional Maxwell equations, and so the question arises whether the proposed extension can also be used in three dimensions and whether it retains its increased stability. This is the case, and details are given in Section \ref{app:maxwell3d} of the Appendix.

\subsection{A collocated sequential explicit method with central derivatives} \label{ssec:central}

The asymmetry of the finite difference approximation that arise in the collocated method \eqref{eq:yeecolloc1}--\eqref{eq:yeecolloc3} after renaming variables in the staggered-grid method \eqref{eq:yeeorig1}--\eqref{eq:yeeorig3} can be addressed by replacing them by central differences, thus also restoring the higher order of accuracy of the spatial discretization:

\begin{scheme}[Central sequential explicit] \label{scheme:central}
\begin{align}
 \frac{(B^z)^{n+1}_{ij} - (B^z)^{n}_{ij}}{\Delta t} &= -\left( \frac{(E^y)^n_{i+1,j}  - (E^y)^n_{i-1,j}}{2\Delta x} - \frac{(E^x)^n_{i,j+1}  - (E^x)^n_{i,j-1}}{2\Delta y}  \right ) \label{eq:centralmaxwell1}\\
 \frac{(E^x)^{n+1}_{ij} - (E^x)^{n}_{ij}}{\Delta t} &= \frac{(B^z)^{n+1}_{i,j+1} - (B^z)^{n+1}_{i,j-1}}{2\Delta y} \\
 \frac{(E^y)^{n+1}_{ij} - (E^y)^{n}_{ij}}{\Delta t} &= -\frac{(B^z)^{n+1}_{i+1,j} - (B^z)^{n+1}_{i-1,j}}{2\Delta x} \label{eq:centralmaxwell3}
\end{align}
\end{scheme}

This is the FVTD method from \cite{remaki99,piperno02}. It is very easy to implement, as it amounts to a particular time-stepping of just the central scheme. This method is also involution preserving (see Theorem \ref{thm:statpresfortriangimpl}) and stable with a maximum CFL of $\sqrt{2}$ (see Corollary \ref{cor:centralstability}).

Consider also a multi-dimensionally extended central collocated method:

\begin{scheme}[Central sequential explicit extended] \label{scheme:centralextended}
\begin{align}
 \frac{(B^z)^{n+1}_{ij} - (B^z)^{n}_{ij}}{\Delta t} &= \label{eq:centralextended1}\\\nonumber &\!\!\!\!\!\!\!\!\!\!\!\!-\left( \frac{\langle(E^y)^n_{i+1}\rangle_j  - \langle(E^y)^n_{i-1}\rangle_j}{2\Delta x} - \frac{\langle(E^x)^n_{\cdot,j+1}\rangle_i  - \langle(E^x)^n_{\cdot,j-1}\rangle_i}{2\Delta y}  \right ) \\
 \frac{(E^x)^{n+1}_{ij} - (E^x)^{n}_{ij}}{\Delta t} &= \frac{\langle(B^z)^{n+1}_{\cdot,j+1}\rangle_i - \langle(B^z)^{n+1}_{\cdot,j-1}\rangle_i}{2\Delta y} \\
 \frac{(E^y)^{n+1}_{ij} - (E^y)^{n}_{ij}}{\Delta t} &= -\frac{\langle(B^z)^{n+1}_{i+1}\rangle_j - \langle(B^z)^{n+1}_{i-1}\rangle_j}{2\Delta x}  \label{eq:centralextended3}
\end{align} 
with notation defined by \eqref{eq:averagingoperator}.
\end{scheme}

This scheme is stable until $\text{CFL} = 2$ (see Corollary \ref{cor:centralextendedstability}) and involution preserving by Theorem \ref{thm:statpresfortriangimpl}. The increased CFL number is not surprising, as the method uses a five-point stencil in each direction. This can be seen by inserting the values of $(B^z)^{n+1}$ into the equations for $E^x$, $E^y$.

\section{The low Mach number limit of sequential explicit collocated methods for linear acoustics} \label{sec:acoustics}

The equations of linear acoustics 
\begin{align}
 \del_t \vec v + \nabla p &= 0 & \vec v &\colon \mathbb R^+_0 \times \mathbb R^d \to \mathbb R^d\\
 \del_t p + c^2 \nabla \cdot \vec v &= 0 & p &\colon \mathbb R^+_0 \times \mathbb R^d \to \mathbb R
\end{align}
are a linearization of the Euler equations and govern the evolution of small perturbations (sound waves) on top of a background of constant density, velocity and pressure. For a derivation, see e.g. \cite{barsukow17}. With an evolution based on characteristic cones rather than characteristics, in multiple spatial dimensions they are as important as linear advection is for one-dimensional problems. The quest for adequate numerical methods for linear acoustics can be understood as paving the way towards truly multi-dimensional discretizations for more general problems (\cite{eymann13,roe17,barsukow19activeflux}). As a linear system, they allow for an explicitly known exact solution (\cite{barsukow17}) which employs characteristic cones and spherical means. This is reminiscent of the scalar wave equation, but linear acoustics is more complicated, because the evolution of $\vec v$ \emph{cannot} be rewritten as a component-wise scalar wave equation. This results in e.g. singularity formation (\cite{amadori2015,barsukow17}) for discontinuous initial data in multiple spatial dimensions. 

Linear acoustics is also an important system to gain understanding of the behaviour of numerical methods in the limit of low Mach number. To this end, consider a family of equations parametrized by $\epsilon \in \mathbb R^+$:
\begin{align}
 \del_t \vec v_\epsilon + \frac{\nabla p_\epsilon}{\epsilon^2} &= 0 \label{eq:acousticslowmach1}\\
 \del_t p_\epsilon + c^2 \nabla \cdot \vec v_\epsilon &= 0  \label{eq:acousticslowmach2}
\end{align}
This system is inspired by the low Mach number scaling of the Euler equations, which is introduced below (Section \ref{sec:euler}). To simplify notation, the $\epsilon$-subscript is dropped from now on. As $\epsilon \to 0$, equations \eqref{eq:acousticslowmach1}--\eqref{eq:acousticslowmach2} formally become
\begin{align}
 \nabla p &= 0 & \nabla \cdot \vec v &= 0
\end{align}
which bears a lot of similarity to the low Mach number limit of the Euler equations. For stable discretizations of linear acoustics, the low Mach number limit is equivalent to the limit of long time. It has been shown in \cite{barsukow17a} that the failure of Finite Volume methods for \eqref{eq:acousticslowmach1}--\eqref{eq:acousticslowmach2} to be low Mach number compliant has its origin in the inability of the discrete stationary states of the numerical method to discretize all the stationary states of the PDE. For more details, see \cite{barsukow17a,barsukow18hypproceeding} and a similar discussion in \cite{barsukow20cgk}. By the results of \cite{barsukow17a}, for linear numerical methods involution preservation implies stationarity preservation.

The two-dimensional acoustic equations ($\vec v = (u,v)^\text{T}$) can be obtained from the two-dimensional Maxwell equations \eqref{eq:maxwell1}--\eqref{eq:maxwell3} essentially through renaming:
\begin{align}
 \veccc{u}{v}{p} = \left(  \begin{array}{ccc} 0 & 0 & -1 \\ 0 & 1 & 0 \\ -1 & 0 & 0  \end{array} \right ) \veccc{B^z}{E^x}{E^y}
\end{align}

Because of this close analogy, there is not much left to be shown when transferring a scheme for Maxwell equations to the acoustic case. The involution $\nabla \cdot \vec E = 0$ for the acoustic equations is replaced by the stationary vorticity $\del_t (\nabla \times \vec v) = 0$. The close connection between the different extensions of the Yee scheme and the stationarity preserving method from \cite{barsukow17a} has been also pointed out already. For reference, here some of the versions of the schemes previously stated for the Maxwell equations are given again for linear acoustics, which serves as a preparation for an extension to the full Euler equations.

The original Yee scheme \ref{scheme:yeeoriginal}, applied to linear acoustics reads
\begin{align}
 \frac{p^{n+\frac12}_{i+\frac12,j+\frac12} - p^{n-\frac12}_{i+\frac12,j+\frac12}}{\Delta t} &=  - c^2\left(\frac{u^n_{i+1,j+\frac12}  - u^n_{i,j+\frac12}}{\Delta x} + \frac{v^n_{i+\frac12,j+1}  - v^n_{i+\frac12,j}}{\Delta y} \right ) \label{eq:yeeorigacoustics1} \\
 \frac{u^{n+1}_{i,j+\frac12} - u^{n}_{i,j+\frac12}}{\Delta t} &= -\frac{p^{n+\frac12}_{i+\frac12,j+\frac12} - p^{n+\frac12}_{i-\frac12,j+\frac12}}{\Delta x} \\
 \frac{v^{n+1}_{i+\frac12,j} - v^{n}_{i+\frac12,j}}{\Delta t}& = -\frac{p^{n+\frac12}_{i+\frac12,j+\frac12} - p^{n+\frac12}_{i+\frac12,j-\frac12}}{\Delta y} \label{eq:yeeorigacoustics3}
\end{align}

Its multi-dimensional extension (Scheme \ref{scheme:yeecollocatedextended}) is:
\begin{align}
 \frac{p^{n+1}_{ij} - p^{n}_{ij}}{\Delta t} &= -c^2\left(\frac{u^n_{i+1,j+1}  - u^n_{i,j+1} +u^n_{i+1,j}  - u^n_{ij}}{2\Delta x} \right . \label{eq:yeecollocatedextacoustics1} \\\nonumber &\phantom{mmmmmm} \left. + \frac{v^n_{i+1,j+1}  - v^n_{i+1,j} + v^n_{i,j+1}  - v^n_{ij}}{2\Delta y}   \right )  \\
 \frac{u^{n+1}_{ij} - u^{n}_{ij}}{\Delta t} &= - \frac{p^{n+1}_{ij}  - p^{n+1}_{i-1,j} + p^{n+1}_{i,j-1} - p^{n+1}_{i-1,j-1}}{2 \Delta x} \\
\frac{v^{n+1}_{ij} - v^{n}_{ij}}{\Delta t} &= -\frac{p^{n+1}_{ij} - p^{n+1}_{i,j-1} + p^{n+1}_{i-1,j} - p^{n+1}_{i-1,j-1}}{2 \Delta y}  \label{eq:yeecollocatedextacoustics3}
 \end{align}

and the extended central sequential explicit scheme (Scheme \ref{scheme:centralextended}) becomes

\begin{align}
 \frac{p^{n+1}_{ij} - p^{n}_{ij}}{\Delta t} &= -c^2\left( \frac{\langle u^n_{i+1}\rangle_j  - \langle u^n_{i-1}\rangle_j}{2\Delta x} + \frac{\langle v^n_{\cdot,j+1}\rangle_i  - \langle v^n_{\cdot,j-1}\rangle_i}{2\Delta y}  \right ) \label{eq:centralextendedacoustics1} \\
 \frac{v^{n+1}_{ij} - v^{n}_{ij}}{\Delta t} &= -\frac{\langle  p^{n+1}_{\cdot,j+1}\rangle_i - \langle p^{n+1}_{\cdot,j-1}\rangle_i}{2\Delta y} \\
 \frac{u^{n+1}_{ij} - u^{n}_{ij}}{\Delta t} &= -\frac{\langle p^{n+1}_{i+1}\rangle_j - \langle p^{n+1}_{i-1}\rangle_j}{2\Delta x}  \label{eq:centralextendedacoustics3}
\end{align}

Recall that the two latter schemes have a maximum CFL number $\frac{c \Delta t}{\Delta x}$ of 1 and 2, respectively.

\section{Sequential explicit collocated methods for the full Euler equations} \label{sec:euler}

\subsection{Introduction}

The Euler equations are
\begin{align}
 \del_t \rho + \nabla \cdot (\rho \vec v) &= 0 & \rho &\colon \mathbb R^+_0 \times \mathbb R^d \to \mathbb R^+ \label{eq:euler1}\\
 \del_t (\rho \vec v) + \nabla \cdot (\rho \vec v \otimes \vec v + p \id) &= 0 & \vec v &\colon \mathbb R^+_0 \times \mathbb R^d \to \mathbb R^d\\
 \del_t e + \nabla \cdot (\vec v (e+p)) &= 0 & p &\colon \mathbb R^+_0 \times \mathbb R^d \to \mathbb R^+ \label{eq:euler3}
\end{align}
and, whenever necessary, an equation of state of the ideal gas
\begin{align}
 e &= \frac{p}{\gamma-1} + \frac12 \rho |\vec v|^2 \label{eq:eulereos}
\end{align}
is assumed. 

It is customary to make the low Mach number limit explicit by switching to the \emph{rescaled Euler equations}, i.e. the following $\epsilon$-dependent family of equations:
\begin{align}
 \del_t \rho + \nabla \cdot (\rho \vec v) &= 0 \label{eq:eulerresc1}\\
 \del_t (\rho \vec v) + \nabla \cdot \left(\rho \vec v \otimes \vec v + \frac{p}{\epsilon^2} \id\right) &= 0 \\
 \del_t e + \nabla \cdot (\vec v (e+p)) &= 0 \label{eq:eulerresc3}
\end{align}
\begin{align}
 e &= \frac{p}{\gamma-1} + \frac12 \epsilon^2 \rho |\vec v|^2 \label{eq:eulereosresc}
\end{align}
Here, the same notation is used for rescaled and the original quantities; for details see \cite{klein95,barsukow16,barsukow20cgk}. Note that the local Mach number $M = |v|/\sqrt{\gamma p / \rho}$ scales as $M \sim \epsilon$. From now on, system \eqref{eq:eulerresc1}--\eqref{eq:eulereosresc} is used for the theoretical analysis, as it allows to immediately see the scaling of the individual terms. The implementation uses \eqref{eq:euler1}-\eqref{eq:eulereos} and does not depend on the rescaling process. The limit of \eqref{eq:eulerresc1}--\eqref{eq:eulereosresc} has been studied theoretically e.g. in \cite{metivier01}, while earlier works such as \cite{ebin77,klainerman81} focus on the isentropic case. Therein it is shown, that as the Mach number vanishes, the solutions tend to those of the incompressible Euler equations
\begin{align}
 \nabla \cdot \vec v &\in \mathcal O(\epsilon) & \nabla p &\in \mathcal O(\epsilon^2)
\end{align}
if the initial data are well-prepared.

Assuming differentiability, the equations can be split into the following three terms:
\begin{align}
 \del_t \rho &&&+ (u \del_x + v \del_y) \rho &&+ \rho (\del_x u + \del_y v) &&= 0 \label{eq:euler1splitform}\\
 \del_t (\rho u) &&&+ (u \del_x + v \del_y) (\rho u) &&+ \rho u (\del_x u + \del_y v) &&+ \frac{\del_x p}{\epsilon^2} = 0\\
 \del_t (\rho v) &&&+ (u \del_x + v \del_y) (\rho v) &&+ \rho v (\del_x u + \del_y v) &&+ \frac{\del_y p}{\epsilon^2} = 0\\
 \del_t e &&&+ (u \del_x + v \del_y) e &&+ e (\del_x u + \del_y v) &&+ \del_x (up) + \del_y (vp)= 0 \label{eq:euler4splitform}\\
 \nonumber &&&\text{\,\,\,\,\,\,advection}&&\text{\,\,\,\,\,\,compression}&&\text{nonlinear acoustics}
\end{align}
where for simplicity, the equations are given in two spatial dimensions, and $\vec v = (u,v)$.

By analogy with the acoustic methods described in Section \ref{sec:acoustics}, an all-speed scheme is constructed by first updating the momentum in time, and then using its updated value in the computation of the fluxes of the scalar quantities $\rho$ and $e$\footnote{It also seems possible to do it the other way around.}. A difference here is that the advective operator 
\begin{align}
 \del_t \rho + (u \del_x + v \del_y) \rho &= 0 \\
 \del_t (\rho u) + (u \del_x + v \del_y) (\rho u)&= 0\\
 \del_t (\rho v) + (u \del_x + v \del_y) (\rho v) &= 0\\
 \del_t e + (u \del_x + v \del_y) e &= 0
\end{align}
necessarily involves the quantity that is being updated (i.e. it is ``diagonal''), and so cannot be discretized in any other way than explicitly. This means that it necessitates upwinding, i.e. the inclusion of diffusive terms. For the low Mach number limit this is not a problem, as the difficulties of standard explicit Finite Volume methods all originate in the numerical diffusion associated only to the acoustic operator. The different discretization of acoustics and advection is a topic brought up e.g. in \cite{roe17}.

\subsection{Treatment of compressive terms and Lagrange-Projection methods}

Finally, the question arises whether the compressive terms in \eqref{eq:euler1splitform}--\eqref{eq:euler4splitform} require special treatment. Hereby it is self-evident that the numerical method needs to be consistent with the PDE. The question is whether a central discretization of these terms is enough, or whether additional terms are necessary, or at least advantageous. A short overview of existing results therefore is due:

\begin{enumerate}[1.]
\item Taking advection and compression together one obtains the \textbf{pressureless Euler equations}, a system with all the eigenvalues of the Jacobian equal to the velocity. On the one hand, a Rusanov method for the pressureless Euler equations would therefore involve the same diffusion as if the compressive terms were not present. The Riemann solver between states $(\rho_\text{L},u_\text{L})$ and $(\rho_\text{R}, u_\text{R})$ derived through relaxation (\cite{berthon06}), on the other hand, obtains the following numerical flux
\begin{align}
  f^x &= \vecc{\rho^* u^* }{\rho^* (u^*)^2}     \label{eq:presslessrelaxflux}
\end{align}
with
\begin{align}
 u^* &:= \frac{u_\text L + u_\text R}{2} &
 \rho^* &:= \begin{cases} \displaystyle \frac{\rho_\text L}{1 + \frac{\rho_\text L}{2a} (u_\text R - u_\text L)} & u^* > 0\\ \displaystyle \frac{\rho_\text R}{1 + \frac{\rho_\text R}{2a} (u_\text R - u_\text L)} & u^* \leq 0 \end{cases} \label{eq:presslessrelax}
\end{align}
($a >0$ is the relaxation speed with dimensions ``density $\times$ speed'' that needs to be chosen large enough.) Here, there is a denominator which seems to suggest special treatment of the compression (in one spatial dimension, the compression term amounts to the derivative of $u$). As has been remarked in \cite{barsukow20cgk}, also the relaxation solvers from \cite{bouchut04,bouchut09,chalons10,girardin14} for the full Euler equations contain denominators of that kind. 

\item Consider the much simpler problem of \textbf{conservative advection with spatially non-constant velocity}:
\begin{align}
 \del_t q + \nabla \cdot  (\vec U(\vec x) q) &= 0 & q &\colon \mathbb R^+_0 \times \mathbb R^d \to \mathbb R & &\vec U(\vec x) \text{ given} \label{eq:nonconstadvection}
\end{align}
This equation can also be split into advection $+$ compression: $\nabla \cdot  (\vec U q) = \vec U \cdot \nabla q + q (\nabla \cdot \vec U)$.

In one spatial dimension, \cite{leveque02}, Section 9.2, suggests (for positive $U(x)$)
\begin{align}
 q_i^{n+1} &= q_i^n - \frac{\Delta t}{\Delta x} \Big( U(x_i) q_i^n  - U(x_{i-1}) q_{i-1}^n \Big ) \label{eq:levequeconadv}
\end{align}
This method is derived in \cite{leveque02} as a Riemann solver, replacing $U(x)$ by a piecewise constant approximant.

A locally-linearized solver (Roe-type) would employ (with $\bar U$ a suitable average of $U(x_i)$ and $U(x_{i+1})$)
\begin{align}
 f_{i+\frac12} = \frac{U(x_i) q_i^n + U(x_{i+1}) q_{i+1}^n}{2} - |\bar U| \frac{q_{i+1}^n - q_i^n}{2} \label{eq:roeadvectionnonconst}
\end{align}
and thus, for positive U, and assuming that the average fulfills the Roe condition
\begin{align}
 f(q_{i+1}^n) - f(q_{i}^n) = U(x_{i+1}) q_{i+1}^n - U(x_i) q_i^n  &= \bar U (q_{i+1}^n - q_i^n)
\end{align}
one obtains the same method as \eqref{eq:levequeconadv}
\begin{align}
 f_{i+\frac12} = U(x_i) q_i^n 
\end{align}

From \eqref{eq:roeadvectionnonconst} it is clear that the numerical diffusion is coming only from advection.

\item Employing edge-based velocities, another method suggested for \eqref{eq:nonconstadvection} in \cite{leveque02} (Section 9.5, and again, for positive $U$) is, on the one hand,
\begin{align}
 q_i^{n+1} &= q_i^n - \frac{\Delta t}{\Delta x} \left( U(x_{i+\frac12}) q_i^n  - U(x_{i-\frac12}) q_{i-1}^n \right ) \label{eq:levequeedgebasedconsadv}
\end{align}
and amounts to pure upwind with respect to the edge velocity.

On the other hand, a Lagrange-Projection method would obtain the following discretization:
\begin{align}
 f_{i+\frac12} &= \frac{U(x_{i+\frac12})}{2} \left( \frac{q_i^n}{L_i} + \frac{q_{i+1}^n}{L_{i+1}} \right ) - \frac{|U(x_{i+\frac12})|}{2} \left(  \frac{q_{i+1}^n}{L_{i+1}} - \frac{q_i^n}{L_i} \right ) \label{eq:agrangeporjectioninaction}\\
  L_{i} &:= 1 + \Delta t \frac{U(x_{i+\frac12}) - U(x_{i-\frac12})}{\Delta x}
\end{align}
which for positive $U$ amounts to
\begin{align}
 q_i^{n+1} &= q_i^n - \frac{\Delta t}{\Delta x} \left( U(x_{i+\frac12}) \frac{q_i^n}{1 + \Delta t \frac{U(x_{i+\frac12}) - U(x_{i-\frac12})}{\Delta x}} \right.  \label{eq:lagrprojconsadv} \\\nonumber &\phantom{mmmmmmmmmm} \left. - U(x_{i-\frac12}) \frac{q_{i-1}^n}{1 + \Delta t \frac{U(x_{i-\frac12}) - U(x_{i-\frac32})}{\Delta x}}  \right ) 
\end{align}

The details on the derivation are given in the Appendix \ref{app:lagproj}. This flux (and in particular the denominator) is reminiscent of \eqref{eq:presslessrelax}.

\item In \cite{barsukow20cgk} a simple way of explaining the appearance of the denominators involving the divergence has been given from the point of view of operator splitting: Rewrite \eqref{eq:nonconstadvection} as
\begin{align}
 \del_t q + \vec U(x) \cdot \nabla q + q (\nabla \cdot \vec U(\vec x)) &= 0
\end{align}
Apply now the \emph{backward} Euler method to the ODE-type last term:
\begin{align}
 \frac{q_i^{n+1} - q_i^n}{\Delta t} + \vec U(x) \cdot \nabla q + q^{n+1} \nabla \cdot \vec U(\vec x) &= 0
\end{align}
i.e.
\begin{align}
 q_i^{n+1}  =  \frac{q_i^n - \Delta t \vec U(x) \cdot \nabla q }{1  + \Delta t (\nabla \cdot \vec U(\vec x)) } \label{eq:odecompressionresult}
\end{align}
The advection term $\vec U(x) \cdot \nabla q$ can now be discretized using the upwind method, and the divergence in the denominator by a central derivative.

\end{enumerate}

Thus, generally speaking, Lagrange-Projection methods and Riemann solvers derived by relaxation tend to include an additional treatment of compressive terms. Again, it is worth emphasizing that the additional treatment goes beyond consistency. For example, both the method \eqref{eq:levequeedgebasedconsadv} and \eqref{eq:lagrprojconsadv} are consistent discretizations, but the latter takes the compressive terms into account in a more sophisticated way.

As discussed in the Appendix (Section \ref{app:stabilityeuler}), for the purpose of this work -- i.e. in order to combine a sequential explicit scheme for acoustics with a discretization of pressureless Euler -- the more careful treatment of the compression brings practical advantages for the stability of the resulting method. In particular, it is found experimentally, that neither the flux inspired by \eqref{eq:roeadvectionnonconst}
\begin{align}
 \frac{f(q_\text L) + f(q_\text R)}{2} - \frac12 |u^*| (q_\text R - q_\text L)
\end{align}
with $u^* = \frac{u_\text L + u_\text R}{2}$, nor using $u^* = \frac{u_\text L + u_\text R}{2} - \frac{p_\text R - p_\text L}{2 \rho_L c_L}$, nor the flux
\begin{align}
 \frac{f(q_\text L) + f(q_\text R)}{2} - \frac12 (|u_\text R| q_\text R - |u_\text L| q_\text L)
\end{align}
nor the multi-dimensional flux (inspired by stationarity preserving schemes in \cite{barsukow17a})
\begin{align}
 \frac{f(q_\text L) + f(q_\text R)}{2} - \frac12 \sgn \,u^* \cdot (\text{discretization of } u^* \del_x q + v^* \del_y q)
\end{align}
yield sufficiently stable schemes (see Section \ref{app:stabilityeuler} for more details). Including numerical diffusion associated with compressive terms, however, has helped overcoming these difficulties.

As can be seen from the above overview, there are different ways of including the divergence in the denominator. It can either be preceded by a factor containing $\Delta t$ (as in Lagrange-Projection methods, or naturally appearing in \eqref{eq:odecompressionresult}), or be multiplied with terms involving the sound speed. It is clear that by virtue of the CFL-condition these terms are in the end not very different in value; here preference is given to using a $\Delta t$ factor in the denominator\footnote{In numerical codes with a rigid structure, the presence of a $\Delta t$-factor in the flux might be difficult to implement, in which case one will have to opt for a prefactor involving the speed of sound, as in \eqref{eq:presslessrelax}.}. 

Finally, comparing the Lagrange-Projection method \eqref{eq:lagrprojconsadv} and the relaxation solver \eqref{eq:presslessrelax} one realizes another conceptual difference: whether the divergence in the denominator is cell-based or edge-based. Cell-based velocities naturally give an edge-based divergence, and vice versa. The relaxation solver provides a flux through the interface, and thus it is natural to use a discrete divergence centered on the edge. Moreover, as the velocities are also dependent variables, their edge values are not immediately available. Finally, when each conserved variable is to be divided by the respective cell-based term $L_i$ (as in \eqref{eq:agrangeporjectioninaction}), it is not immediately clear how the central flux for full Euler is to be constructed, or at least such a construction seems cumbersome. Here, therefore, preference is given to constructing a central $+$ advection-upwind flux, and to divide it by a term $1  + \Delta t (\nabla \cdot \vec U(\vec x))$, which then involves a discrete edge-centered divergence (as the flux is edge-centered, too), i.e. schematically:
\begin{align*}
 f^x_{i+\frac12,j} = \frac{\boxed{\text{flux average}}_{i+\frac12,j} - \boxed{\text{upwind difference}}_{i+\frac12,j}}{1  + \Delta t \boxed{\text{divergence}}_{i+\frac12,j}}
\end{align*}

\subsection{All-speed sequential explicit numerical method} \label{ssec:allspeedeulermethod}

The question remaining at this point is whether there are constraints on the choice of the finite difference operators involved. Such a constraint does exist. In fact, the construction of the flux, as outlined above, will involve \emph{two} discrete divergences: the one in the denominator and the one in the energy equation
\begin{align}
 \del_t e + \nabla \cdot (\vec v (e+p)) &= 0
\end{align}
which remains when the highest orders in the $\epsilon$-expansion of $e$ and $p$ are found constant in the limit of low Mach number $\epsilon \to 0$. The interaction of these two divergences shall be exemplified by taking a naive discretization of the divergence centered at an edge $i+\frac12$, for example
\begin{align}
 \frac{u_{i+1,j} - u_{ij}}{\Delta x} + \frac{v_{i+1,j+1} + v_{i,j+1} - v_{i+1,j-1} - v_{i,j-1}}{4 \Delta y}
\end{align}
which shall be used in the denominator. The strategy outlined above would yield the following numerical flux:
\begin{align}
  f_{i+\frac12,j}^x &= \frac{\displaystyle\frac{f^x(q_{i+1,j}) + f^x(q_{ij}) }{2} - \frac12 |u^*_{i+\frac12,j}| (q_{i+1,j} - q_{ij})}{\displaystyle1 + \Delta t \left(\frac{u_{i+1,j} - u_{ij}}{\Delta x} + \frac{v_{i+1,j+1} + v_{i,j+1} - v_{i+1,j-1} - v_{i,j-1}}{4 \Delta y}\right)} \label{eq:eulerfluxx}
\end{align}
with $u^*_{i+\frac12,j} = \frac{ u_{i+1,j} + u_{ij} }{2}$, for example. Observe the presence of a symmetric flux, endowed with the upwind diffusion associated only to the advection operator; the denominator is due to the compression operator. The flux has as many components, as there are conserved quantities; e.g. the density-flux shall be denoted by $( f^\rho)_{i+\frac12,j}^x$. Note that the time indices do not appear in these formulae yet.

This method, however, is not low Mach compliant. Observe that in the numerical flux \eqref{eq:eulerfluxx}, the upwinding and the denominator are $\mathcal O(1)$ in $\epsilon$. Thus, the $\mathcal O(\epsilon^{-2})$ and $\mathcal O(\epsilon^{-1})$ equations are entirely due to the physical fluxes $\frac{f^x(q_{i+1,j}) + f^x(q_{ij}) }{2}$. Therein, all terms are $\mathcal O(1)$, except for the pressure appearing as $\frac{p}{\epsilon^2}$ in the momentum flux. Because of the denominator, however, one cannot immediately conclude that the pressure is constant to $\mathcal O(\epsilon^2)$. This would be the case if the divergence appearing in the denominator would vanish -- not an entirely hopeless situation, as of course, at continuous level, the divergence \emph{does} vanish in the low Mach number limit. 

But in the discrete setting, one cannot guarantee that two different discretizations of the divergence vanish simultaneously -- unless they are very special. The two discretizations of the divergence in question are the edge-centered discretization in the denominator and the cell-centered divergence in the energy equation. In order to be able to prove low Mach compliance of the numerical method, one needs to make sure that both the divergence arising from the central discretization of the energy equation and the divergence in the denominator are, if not the same, then at least vanishing simultaneously. 

Such questions concerning simultaneously vanishing discretizations of the divergence have already been dealt with in \cite{barsukow17a,barsukow20cgk} and the results can be used here. In fact, they have already appeared in parts upon the definition of the multi-dimensionally extended central scheme \ref{scheme:centralextended}. Essentially, to vanish simultaneously, one of the divergences has to be an average of the other. However, the simple five-point divergence (in 2-d) cannot be written as an average of an edge-centered discrete divergence without losing symmetry (see \cite{barsukow17a} for a proof). The special discretization of the divergence required at this point has to be truly multi-dimensional, i.e. it needs to involve 9 cells (in 2-d). Due to the thus increased number of terms in the multi-dimensionally extended finite difference operators before stating the numerical flux the following notation is introduced.

\begin{definition}
 In one spatial dimension, the \textbf{single-bracket operators} are a jump and a sum (both located at an edge $i+\frac12$) and are denoted by
\begin{align}
 [q]_{i+\frac12} &:= q_{i+1} - q_i & \{q\}_{i+\frac12} &:= q_{i+1} + q_i\\
 [q]_{i\pm1} &:= q_{i+1} - q_{i-1}
\end{align}
One can consider the brackets as operators on sequences, e.g. $[\cdot] \colon \mathbb R^{\mathbb N} \to \mathbb R^{\mathbb N}$, the subscript then referring to one of the elements of $[q]$, with purely notational natural shifts by $\frac12$. The \textbf{double-bracket operators} naturally arise by applying the jump/sum twice. The result then is again located at a cell:
\begin{align}
 [[q]]_{i\pm\frac12} &:= [q]_{i+\frac12} - [q]_{i-\frac12} & \{ \{ q \}\}_{i\pm\frac12}&:= \{q\}_{i+\frac12} + \{q\}_{i-\frac12}\\
 &= q_{i+1} - 2 q_i + q_{i-1} & &= q_{i+1} + 2 q_i + q_{i-1}
\end{align}
The only nontrivial identity is
\begin{align}
 \{[q]\}_{i\pm\frac12} = [q]_{i+\frac12} + [q]_{i-\frac12} = [q]_{i\pm1}
\end{align}
The average $\langle \cdot \rangle_j$ introduced in \eqref{eq:averagingoperator} is precisely $\{\{ \cdot \}\}_{j\pm\frac12}$.

For Cartesian grids in multiple spatial dimensions, the finite difference operators are applied to the different indices individually. The indices are ordered: first come those associated to the $x$-direction ($i$), then the ones of the $y$-direction ($j$). One first defines the finite differences that are essentially one-dimensional, e.g.:
\begin{align}
 [q]_{i,j\pm\frac12} \equiv [q_i]_{j\pm\frac12} &:= q_{i,j+1} - q_{ij} & [q]_{i\pm\frac12,j} &:= q_{i+1,j} - q_{ij}
\end{align}
and analogously for other single-bracket or double-bracket operators. The combination of finite differences in different directions, which is truly multi-dimensional, is then defined by expanding from outside inwards:
\begin{align}
 \{\{[q]_{i+\frac12}\}\}_{j+\frac12} &= [q]_{i+\frac12,j+1} + 2 [q]_{i+\frac12,j} + [q]_{i+\frac12,j-1} \\&=  q_{i+1,j+1} - q_{i,j+1} + 2q_{i+1,j} - 2q_{ij} + q_{i+1,j-1} - q_{i,j-1}\\
 [[q]_{i\pm1}]_{j\pm1} &= [q]_{i\pm1,j+1} - [q]_{i\pm1,j-1} \\&=  q_{i+1,j+1} - q_{i-1,j+1} - q_{i+1, j-1} + q_{i-1,j-1}
\end{align}
As one easily proves, the order of expansion (starting with the outer brackets, or the inner ones) is, in fact, immaterial. It is sometimes useful to reverse the order of indices in the following way, e.g.:
\begin{align}
 \left[ [q]_{\cdot,j+\frac12} + c \{q\}_{\cdot,j+\frac12} \right]_{i\pm1} = [[q]_{i\pm1}]_{j+\frac12} + c \{[q]_{i\pm1}\}_{j+\frac12}
\end{align}
Observe that a dot is inserted to avoid ambiguity.

\end{definition}
The number of square brackets in an expression gives the order of the differential operator approximated. Note also that the curly brackets are a sum which needs to be divided by 2 as many times as there are pairs of curly brackets in order to become an average.

In view of what has been said above about the necessity to provide a dicretization of divergence that has a chance to vanish simultaneously with another finite difference divergence, instead of \eqref{eq:eulerfluxx}, the following extended flux shall be considered

\begin{align}
  f_{i+\frac12,j}^x &= \frac{\frac{ \{\{ \,\{ f(q) \}_{i+\frac12} \}\}_{j\pm\frac12}}{8} - \frac12 |u^*_{i+\frac12,j}| [q]_{i+\frac12,j}}{1 + \Delta t \left(\frac{\{\{[u]_{i+\frac12}\}\}_{j\pm\frac12}}{4\Delta x} + \frac{[\{ v \}_{i+\frac12}]_{j\pm1}}{4 \Delta y}\right)} \label{eq:eulerfluxxmultid}
\end{align}
Here, the advective terms have been not extended to multi-d, but they could be. The rationale behind the choice of these operators and the precise moment when certain discretizations vanish simultaneously with others is detailed in the proof of Theorem \ref{thm:lowmacheuler}.

The flux in $y$-direction reads
\begin{align}
  f_{i,j+\frac12}^y &= \frac{\frac{ \{\, \{ \{ f(q) \}\}_{i\pm\frac12} \}_{j+\frac12}}{8} - \frac12 |v^*_{i,j+\frac12}| [q]_{i,j+\frac12}}{1 + \Delta t \left(\frac{ \{ [u]_{i\pm1} \}_{j+\frac12}}{4\Delta x} + \frac{[\{\{ v \}\}_{i\pm\frac12}]_{j+\frac12}}{4 \Delta y}\right)} \label{eq:eulerfluxymultid}
\end{align}

\begin{figure}
 \centering
 \includegraphics[width=0.4\textwidth]{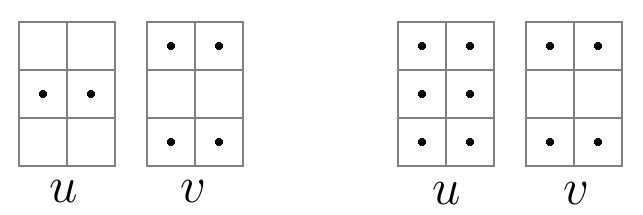}
 \caption{Stencils of the edge-centered divergence discretizations accounting for the compressive terms. \emph{Left}: Stencil of the divergence appearing in the denominator of \eqref{eq:eulerfluxx}. \emph{Right}: Stencil of the divergence appearing in the denominator of \eqref{eq:eulerfluxxmultid}.}
 \label{fig:stencilsdivergence}
\end{figure}

\newcommand{\rhs}{\mathrm{RHS}}

In order to describe the sequential explicit time integration it is necessary to make explicit the dependence on the individual conserved quantities. To this end, write $q_{ij}^n = (\rho_{ij}^n, \vec m^n_{ij}, e^n_{ij})$, having denoted the momentum $\rho \vec v $ by $\vec m$. Denote by $\{ q_{\cdot,\cdot} \}$ the set of all values on the grid and define the right-hand-side
\begin{align}
 \rhs_{ij}( \rho_{\cdot,\cdot}, \vec m_{\cdot,\cdot}, e_{\cdot,\cdot}  ) := \frac{ f_{i+\frac12,j}^x -  f_{i-\frac12,j}^x}{\Delta x} -\frac{f_{i,j+\frac12}^y -  f_{i,j+\frac12}^y}{\Delta y} \label{eq:eulerrhs}
\end{align}

Again, the right-hand side has as many components as there are conserved quantities; the $\rho$-component, for example, shall be denoted by $\rhs^\rho_{ij}$:
\begin{align}
 \rhs_{ij}^\rho( \rho_{\cdot,\cdot}, \vec m_{\cdot,\cdot}, e_{\cdot,\cdot}  ) := \frac{ ( f^\rho)_{i+\frac12,j}^x -  ( f^\rho)_{i-\frac12,j}^x}{\Delta x} -\frac{( f^\rho)_{i,j+\frac12}^y -  ( f^\rho)_{i,j+\frac12}^y}{\Delta y} 
\end{align}

The simplest version of the sequential explicit time integration then reads as follows: update first the momentum equations, and then density and energy using the new value of the momentum. (Other choices concerning the order of equations, or the usage of the new value of $\rho$ in the last equation are possible, but not conceptually different.) Together, this gives

\begin{scheme} \label{scheme:euler}
\begin{align}
 \vec m^{n+1}_{ij} &= \vec m_{ij}^n - \Delta t \,\rhs^{\vec m}_{ij}( \rho^{n}_{\cdot,\cdot}, \vec m^n_{\cdot,\cdot}, e^{n}_{\cdot,\cdot}  ) \label{eq:updateeulermomentum}\\
 \rho^{n+1}_{ij} &= \rho_{ij}^n - \Delta t\, \rhs^\rho_{ij}( \rho^n_{\cdot,\cdot}, \vec m^{n+1}_{\cdot,\cdot}, e^n_{\cdot,\cdot}  ) \label{eq:updateeulerrho}\\
 e^{n+1}_{ij} &= e_{ij}^n - \Delta t \,\rhs^e_{ij}( \rho^n_{\cdot,\cdot}, \vec m^{n+1}_{\cdot,\cdot}, e^n_{\cdot,\cdot}  ) \label{eq:updateeulerene}
\end{align}
with the definition \eqref{eq:eulerrhs} and the fluxes \eqref{eq:eulerfluxxmultid}--\eqref{eq:eulerfluxymultid}.
\end{scheme}

In practice, the computations \eqref{eq:updateeulermomentum} are performed first for all cells in the grid, then the fluxes of the right-hand sides of \eqref{eq:updateeulerrho}--\eqref{eq:updateeulerene} are computed, and finally \eqref{eq:updateeulerrho}--\eqref{eq:updateeulerene} is performed for all grid cells. This is a very similar computational effort as for a usual, fully explicit numerical method.

\begin{theorem}[Low Mach compliance of Scheme \ref{scheme:euler}] \label{thm:lowmacheuler}
 As $\epsilon \to 0$, the limit equations of \eqref{eq:updateeulermomentum}--\eqref{eq:updateeulerene} are solved by
 \begin{align}
  p_{ij} &= \mathrm{const} + \mathcal O(\epsilon^2) \label{eq:limitdiscreulerp} \\
  \mathscr D_{i+\frac12,j+\frac12} &:= \frac{u_{i+1,j} - u_{ij} + u_{i+1,j+1} - u_{i,j+1}}{2 \Delta x} \label{eq:limitdiscreulerdiv} \\ & \phantom{mmmmm} \nonumber + \frac{v_{i+1,j+1} - v_{i+1,j} + v_{i,j+1} - v_{ij}}{2 \Delta y} \in \mathcal O(\epsilon) 
 \end{align}
\end{theorem}
\begin{proof}
 To perform the asymptotic analysis, every quantity is expanded as a power series in $\epsilon$, e.g.
 \begin{align}
  p_{ij} = p_{ij}^{(0)} + \epsilon p_{ij}^{(1)} + \epsilon^2 p_{ij}^{(2)} + \mathcal O(\epsilon^3)
 \end{align}

 As there are no $\epsilon$-factors appearing explicitly in \eqref{eq:eulerfluxxmultid} apart from the $1/\epsilon^2$ scaling of the pressure in the momentum flux, the limit equations read ($\ell = 0, 1$)
 \begin{align}
  \left [ \frac{\frac{ \{\{ \{ p^{(\ell)}\}_{\cdot} \}\}_{j\pm\frac12}   }{8 \Delta x} }{1 + \Delta t \left(\frac{\{\{[u^{(0)}]_{\cdot}\}\}_{j\pm\frac12}}{4\Delta x} + \frac{[\{ v^{(0)} \}_{\cdot}]_{j\pm1}}{4 \Delta y}\right)} \right ]_{i\pm\frac12} &= 0 \label{eq:prooflowmachpx}\\
  \left[ \frac{  \{ \{\{  p^{(\ell)} \}\}_{i+\frac12}  \}  }{1 + \Delta t \left(\frac{ \{ [u^{(0)}]_{i\pm1} \}}{4\Delta x} + \frac{[\{\{ v^{(0)} \}\}_{i\pm\frac12}]}{4 \Delta y}\right)} \right ]_{j\pm\frac12} &= 0 \label{eq:prooflowmachpy}
 \end{align}
 Both divergences vanish when $\mathscr D^{(0)}_{i+\frac12,j+\frac12}$ vanishes identically, for
 \begin{align}
  \frac{\{\{[u^{(0)}]_{i+\frac12}\}\}_{j\pm\frac12}}{4\Delta x} + \frac{[\{ v^{(0)} \}_{i+\frac12}]_{j\pm1}}{4 \Delta y} &= \frac12 \{\mathscr D^{(0)}_{i+\frac12} \}_{j\pm\frac12} = \frac{\mathscr D^{(0)}_{i+\frac12,j+\frac12} + \mathscr D^{(0)}_{i+\frac12,j-\frac12}}{2}\\
  \frac{ \{ [u^{(0)}]_{i\pm1} \}_{j+\frac12}}{4\Delta x} + \frac{[\{\{ v^{(0)} \}\}_{i\pm\frac12}]_{j+\frac12}}{4 \Delta y} &= \frac12 \{ \mathscr D^{(0)} \}_{i\pm\frac12,j+\frac12} = \frac{\mathscr D^{(0)}_{i+\frac12,j+\frac12} + \mathscr D^{(0)}_{i-\frac12,j+\frac12}}{2}
 \end{align}
 Therefore, indeed, \eqref{eq:prooflowmachpx}--\eqref{eq:prooflowmachpy} are fulfilled. The $\mathcal O(1)$ energy equation then reads
 \begin{align}
  (e^{(0)} + p^{(0)}) \left( \frac{\{\{ [ u^{(0)} ]_{i\pm1} \}\}_{j\pm\frac12}}{8 \Delta x} + \frac{[\{\{ v^{(0)} \}\}_{i\pm\frac12} ]_{j\pm1} }{8 \Delta y} \right ) = 0 
 \end{align}
 The diffusion associated with advection vanishes because $e^{(0)}$ is spatially constant as soon as $p^{(0)}$ is. The divergence appearing in the energy equation is
 \begin{align}
   \frac{\{\{ [ u^{(0)} ]_{i\pm1} \}\}_{j\pm\frac12}}{8 \Delta x} + \frac{[\{\{ v^{(0)} \}\}_{i\pm\frac12} ]_{j\pm1} }{8 \Delta y}  = \frac14 \{\{ \mathscr D \}_{i\pm\frac12} \}_{j\pm\frac12}
 \end{align}
 and thus vanishes as well. This completes the proof.
 
\end{proof}

\emph{Note}: Consider discrete data such that two unrelated divergence discretizations vanish simultaneously. By consistency, their difference is a discretization of, in general, second derivatives of the data. If both discrete divergences vanish, so does their difference, and thus the data must be such that, additionally to the divergences, some discrete second derivatives vanish. These data therefore are a discretization of not just a divergenceless vector field, but also satisfy additional constraints. For low Mach number compliance, however, only the divergence constraint should appear, and no other. The appearance of additional constraints which reduce the set of divergenceless vector fields is at the origin of non-low Mach compliant numerical schemes, and leads to artefacts visible in simulations. For further discussion of this the reader is referred to \cite{barsukow20cgk}. Here, it is shown that \eqref{eq:limitdiscreulerp}--\eqref{eq:limitdiscreulerdiv} is a limit solution of the method without any further constraints.

An outline of theoretical results on the stability of this method for $\Delta y = \Delta x$ is given in Section \ref{app:stabilityeuler}. It is shown that one expects a CFL condition
\begin{align}
 \frac{\Delta t}{\Delta x} \leq \frac{1}{|\bar u| + |\bar v|}
\end{align}
for speeds much larger than the speed of sound (advective regime), and
\begin{align}
 \frac{  \Delta t \bar c}{ \epsilon } < \sqrt{\frac{\gamma}{2}} \Delta x
\end{align}
for vanishing velocities (acoustic regime). Numerical studies of the intermediate regime, as detailed in Section \ref{app:stabilityeuler}, suggest therefore a CFL condition of the form
\begin{align}
 \frac{\Delta t}{\Delta x} < \frac{1}{|\bar u| + |\bar v| + \frac{\bar c}{\epsilon} \sqrt{\frac{2}{\gamma}}} \label{eq:cflconditioneuler}
\end{align}

It is reassuring that the linearized method (in case $\bar u = \bar v = 0$) is stationarity preserving, which by results from \cite{barsukow17a} for linear acoustics is equivalent to its low Mach number compliance.

\subsection{Numerical results}

The ideal gas equation of state with $\gamma = 1.4$ is used everywhere. The CFL number mentioned in the texts refers to the time step condition as defined in Equation \ref{eq:cflconditioneuler}.

\subsubsection{Shock tubes}

\begin{figure} 
 \centering
 \includegraphics[width=0.75\textwidth]{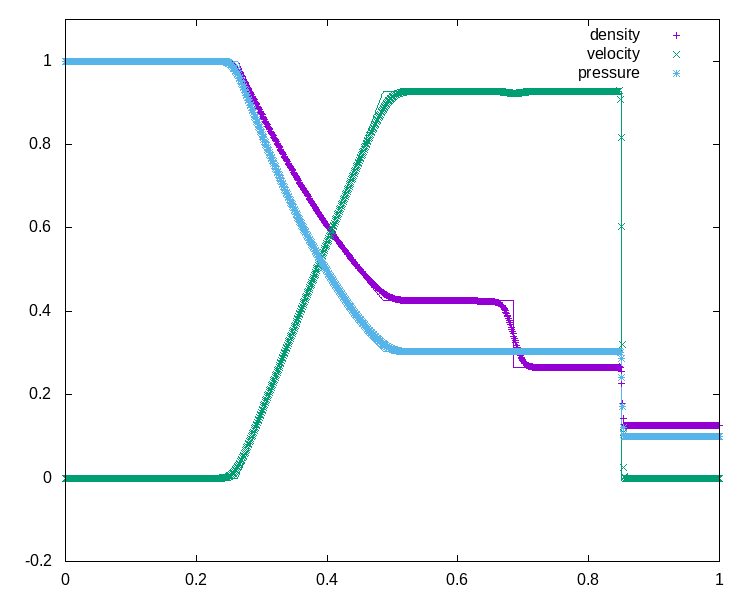}
 \caption{Shock tube test 1 solved with Scheme \ref{scheme:euler} with $\text{CFL} = 0.65$ and $\Delta x = 1/1000$ at time $t=0.2$. The solid lines show the exact solution.}
 \label{fig:sod}
\end{figure}

\begin{figure} 
 \centering
 \includegraphics[width=0.75\textwidth]{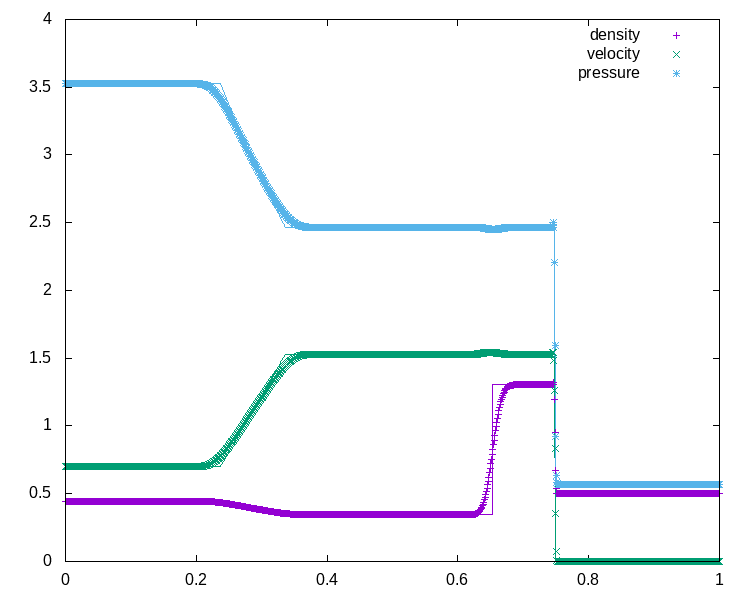}
 \caption{Shock tube test 2 solved with Scheme \ref{scheme:euler} at $t=0.1$. Otherwise, the setup is that of Figure \ref{fig:sod}.}
 \label{fig:lax}
\end{figure}

\begin{figure} 
 \centering
 \includegraphics[width=0.75\textwidth]{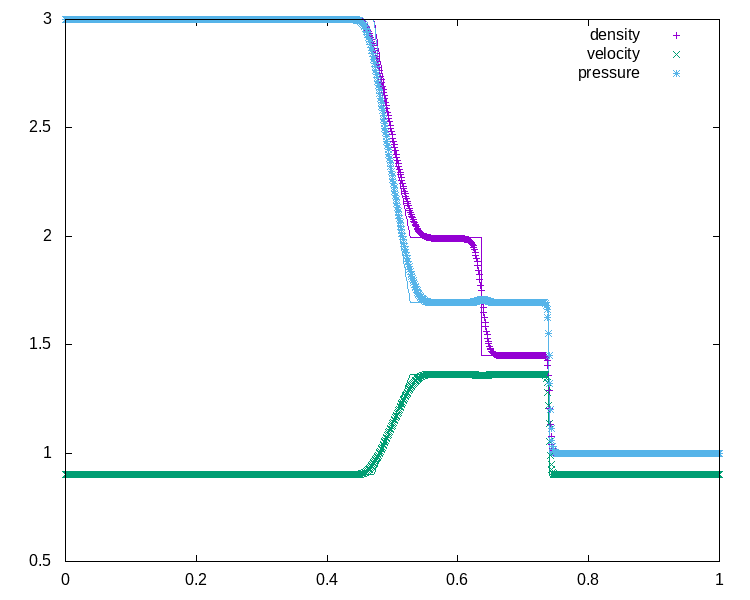}
 \caption{Shock tube test 3 solved with Scheme \ref{scheme:euler} at $t=0.1$. Otherwise, the setup is that of Figure \ref{fig:sod}. The transonic rarefaction does not show any entropy-violating artefacts.}
 \label{fig:leveque}
\end{figure}

Three standard test cases are considered in order to demonstrate the ability of Scheme \ref{scheme:euler} to cope with supersonic phenomena. They are summarized in Table \ref{tab:shocks}. The results are shown in Figures \ref{fig:sod}--\ref{fig:leveque}. The solution is computed on a larger grid in order to exclude the influence of boundaries.

\begin{table}
 \centering
 \begin{tabular}{c|c||c|c|c||c|c|c}
    ID&Name&$p_\text L$ & $\rho_\text L$ & $u_\text L$ & $p_\text R$ & $\rho_\text R$ & $u_\text{R}$ \\
     \hline
    1&Sod test (\cite{sod78}) & 1 & 1 & 0 & 0.1 & 0.125 & 0\\
    2&Lax test (\cite{lax64}) & 3.528 & 0.445 & 0.698 & 0.571 & 0.5 & 0\\
    3&Leveque test (\cite{leveque02}) & 3 & 3 & 0.9 & 1 & 1 & 0.9
 \end{tabular}
 \caption{Overview of parameters of shock tube tests.}
 \label{tab:shocks}
\end{table}

It has been found that for many shock tube setups the time step can be chosen up to $\text{CFL} = 1$, in line with the stability analysis. However, it has been also found that Shock tube test 1 violates positivity unless the CFL number is below 0.7. For better comparison, therefore, all the test cases have been run with a lower CFL number. An in-detail investigation of positivity preservation, however, is subject of future work.

\subsubsection{Incompressible vortex} 

The low Mach number compliance is assessed experimentally using an divergencefree stationary vortex. This flow can be endowed with any Mach number by modifying the background pressure, and thus the speed of sound. The setup is as follows:
\begin{align}
  \rho &= 1  \qquad  v_\phi = \begin{cases}
                         	5 r & r < 0.2  \\ 2-5 r & 0.2 \leq  r < 0.4 \\ 0 &\text{else} \end{cases}
  \\p &=  \begin{cases} p_0 + 12.5 r^2 & r < 0.2 \\ p_0 + 4 \ln(5r) + 4 - 20r + 12.5 r^2 & 0.2 \leq  r < 0.4 \\ p_0 + 4 \ln 2 - 2 & \text{else}
                        \end{cases}
\end{align}
with $p_0 := \frac{1}{\gamma \mathcal M^2} - \frac12$. The definition of $\mathcal M$ here is such that it is equal to the maximum Mach number of the flow.

\begin{figure} 
 \centering
 \includegraphics[width=0.45\textwidth]{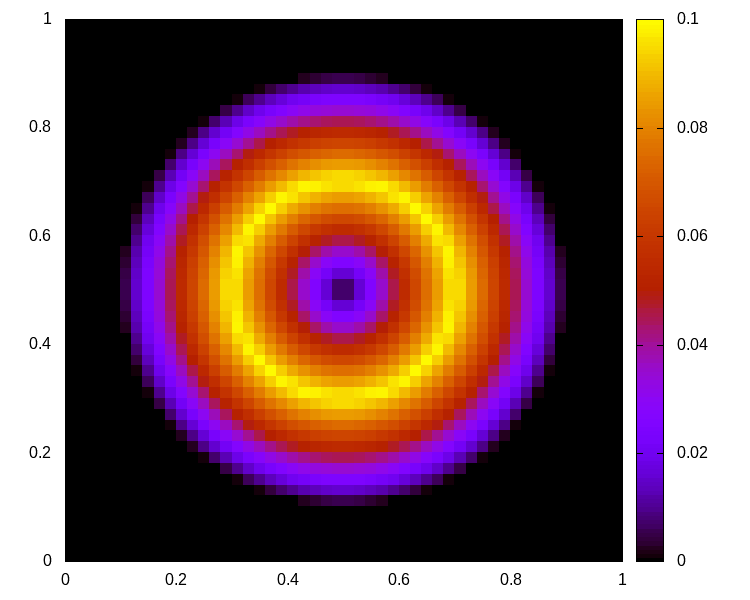} \hfill \includegraphics[width=0.45\textwidth]{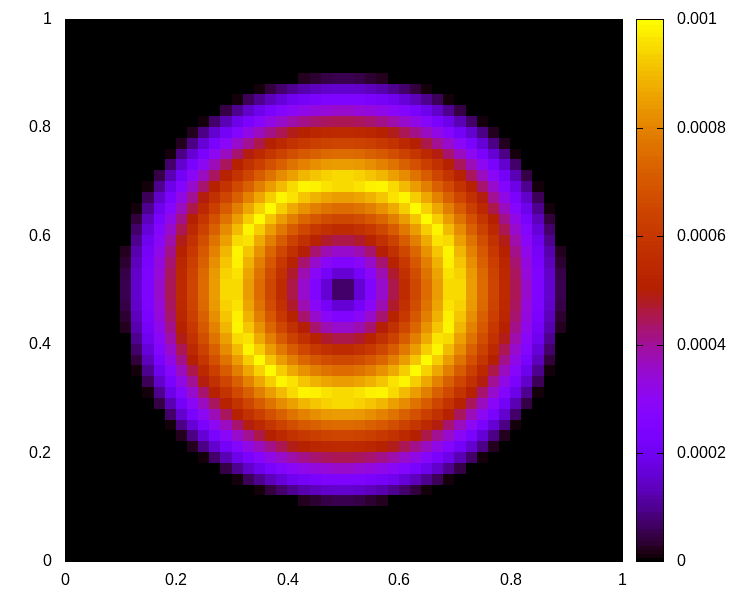}\\
 \includegraphics[width=0.45\textwidth]{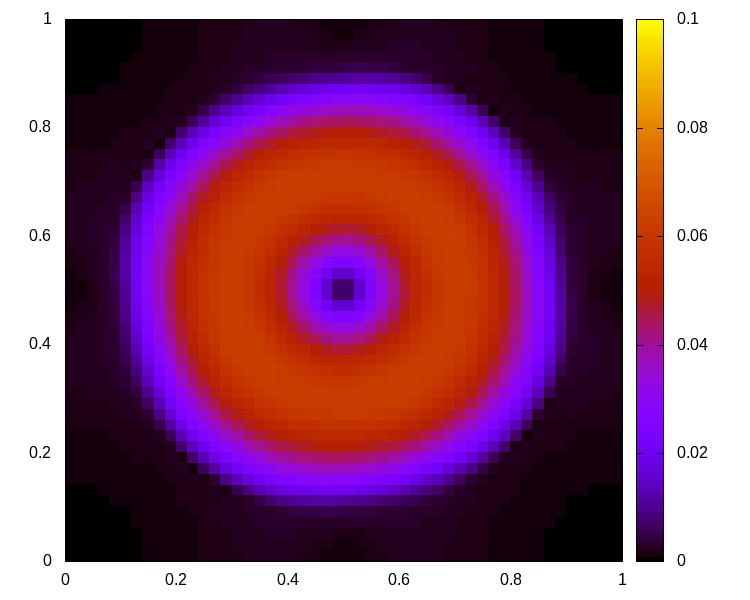} \hfill \includegraphics[width=0.45\textwidth]{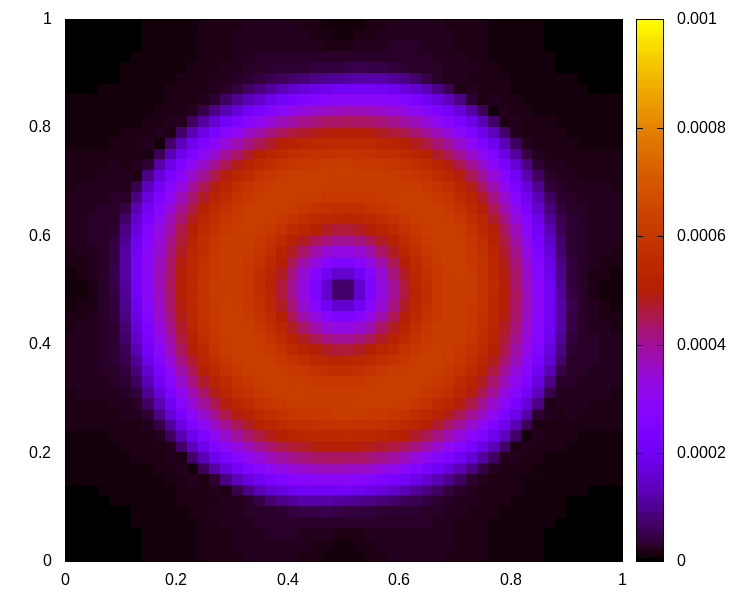}
 \caption{Incompressible stationary vortex solved with Scheme \ref{scheme:euler} with $\text{CFL = 0.9}$ and $\Delta x = \Delta y = 1/50$. Color coded is the Mach number.  \emph{Top}: Initial data, and exact solution for all times. \emph{Bottom}: Numerical results at time $t=1$. \emph{Left}: $\mathcal M = 10^{-1}$. \emph{Right}: $\mathcal M = 10^{-3}$. The presence of numerical diffusion is obvious (Scheme \ref{scheme:euler} is first-order accurate), but its effect is independent of Mach number. This demonstrates experimentally the low Mach compliance of Scheme \ref{scheme:euler}.}
 \label{fig:gresho}
\end{figure}

Experimentally, low Mach compliance manifests itself in a Mach number-independent numerical evolution of the vortex. While the exact solution is stationary, the numerical solution of course displays a certain amount of numerical diffusion. What matters for low Mach compliance is that this diffusion is asymptotically independent of $\epsilon$. This is the case, as shown in Figure \ref{fig:gresho}.

\begin{figure}
 \centering
 \includegraphics[width=0.49\textwidth]{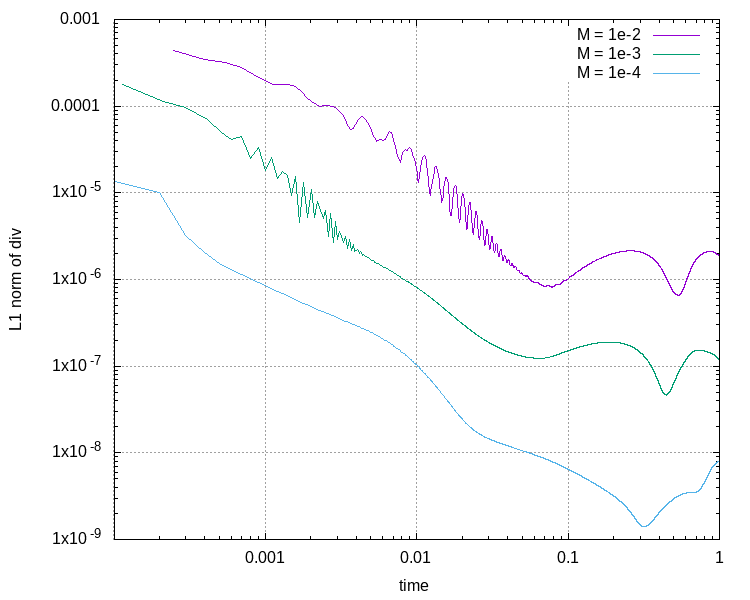} \hfill
 \includegraphics[width=0.49\textwidth]{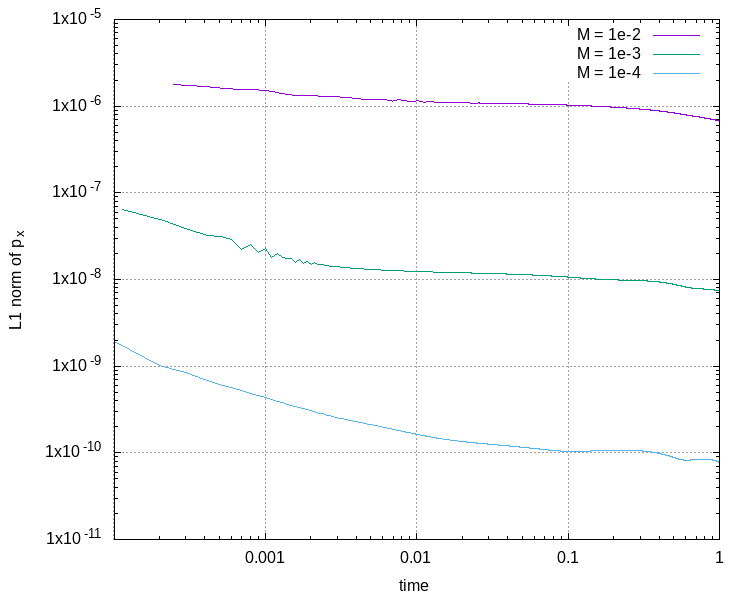}
 \caption{Analysis of the low Mach limit on a $50 \times 50$ grid using CFL = 0.9 for the incompressible vortex. \emph{Left}: The $\ell^1$ norm of the discrete divergence \eqref{eq:limitdiscreulerdiv}, i.e. $\frac{1}{N} \sum_{ij} \Big( \{[u]_{i+\frac12}\}_{j+\frac12} + [\{v\}_{i+\frac12}]_{j+\frac12} \Big)$ is shown as a function of time for different values of $\mathcal M$. \emph{Right}: The $\ell^1$ norm of the discrete $x$-derivative of the pressure, i.e. $\frac{1}{N} \sum_{ij}[\mathcal M^2 p]_{i\pm1}/2$ is shown (the perpendicular derivative is exactly the same and is not shown). Here, $N$ is the number of cells in the grid. One observes the theoretical scaling.}
 \label{fig:lowmachscaling}
\end{figure}

It is also possible to measure the scalings $\nabla p_\epsilon \in \mathcal O(\epsilon^2)$ and $\nabla \cdot \vec v_\epsilon \in \mathcal O(\epsilon)$, where notation has been reintroduce to distinguish the rescaled quantities $p_\epsilon, \vec v_\epsilon$ from their non-rescaled counterparts. There is no difference for $\vec v$, in fact, but $p = \frac{1}{\epsilon^2} p_\epsilon = \frac{1}{\epsilon^2} (p^{(0)} + \epsilon p^{(1)} + \epsilon^2 p^{(2)})$. Thus, the scaling of the non-rescaled variable is $\nabla p_\epsilon \in \mathcal O(1)$. As the numerical implementation uses non-rescaled equations, there is no notion of $\epsilon$ in the code. However, the pressure of the vortex is $p = \frac{\text{const}}{\mathcal M^2} + \mathcal O(1)$. As $M \in \mathcal O(\epsilon)$, for convenience, the rescaled pressure $\mathcal M^2 \nabla p$ is used in the experimental setup, for which again the usual relation $\nabla(\mathcal M^2 \nabla p) \in \mathcal O(\mathcal M^2)$ should hold true. In Figure \ref{fig:lowmachscaling}, the evolution of $\nabla(\mathcal M^2 p)$ and $\nabla \cdot \vec v$ is shown as a function of time for different values of $\mathcal M \sim \epsilon$ for the incompressible vortex. The gradient and divergence are measured using a central extended discretization of the gradient and of the divergence according to Theorem \ref{thm:lowmacheuler}. One observes clearly the theoretical scalings. For comparison with other methods, the reader is referred to Figures 5 and 7 in \cite{barsukow20cgk}. There, among other, the behaviour of a non low Mach compliant method is shown, which deviates from the analytic scalings.

\subsubsection{Convergence study}

\begin{figure}
 \centering
 \includegraphics[width=0.7\textwidth]{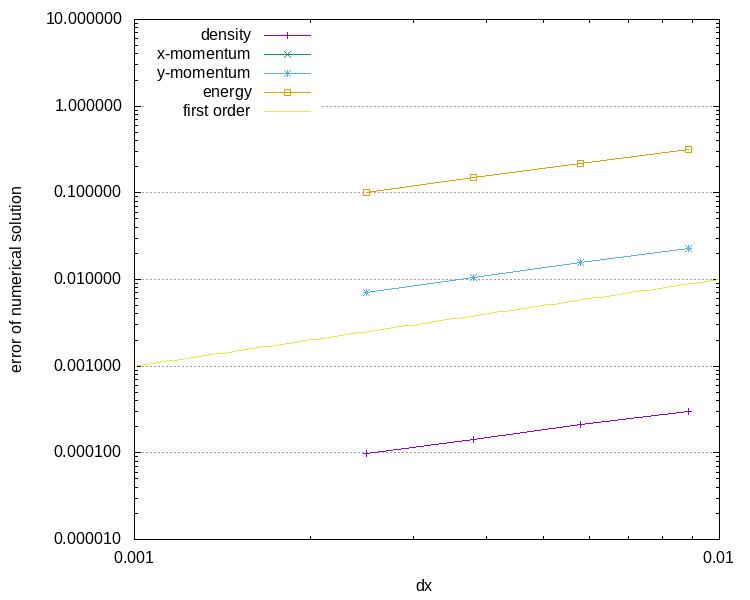}
 \caption{Convergence analysis of Scheme \ref{scheme:euler}. The $\ell^1$ errors of the conserved quantities are shown at time $t=0.05$ and $\mathcal M = 0.3$ as a function of the grid spacing $\Delta x = \Delta y$. The CFL number is 0.9. One observes first order convergence.}
 \label{fig:convergence}
\end{figure}

Scheme \ref{scheme:euler} is formally first-order accurate. To assess this behaviour in the experiment a smooth version of the vortex presented above is used. The initial data are
\begin{align}
  \rho &= 1  \qquad\qquad  v_\phi = v_0 r^2 \exp(-\alpha r)
  \\p &= p_0 + \frac{v_0^2}{8 \alpha^4} \Big(3 + \exp(-2\alpha r) (-3-2 \alpha r (3+\alpha r (3+2 \alpha r))) \Big) 
\end{align}
with $p_0 := \frac{20}{\gamma \mathcal M^2}$, $v_0 = \frac{\alpha^2}{0.13}$ and $\alpha = 20$. The definition of $\mathcal M$ here again corresponds to the maximum Mach number of the flow. This vortex is a stationary solution of the Euler equations. Figure \ref{fig:convergence} shows the experimental errors, which confirm first-order accuracy of the method.

\subsubsection{Kelvin-Helmholtz instability}

\begin{figure}
 \centering
 \includegraphics[width=0.43\textwidth]{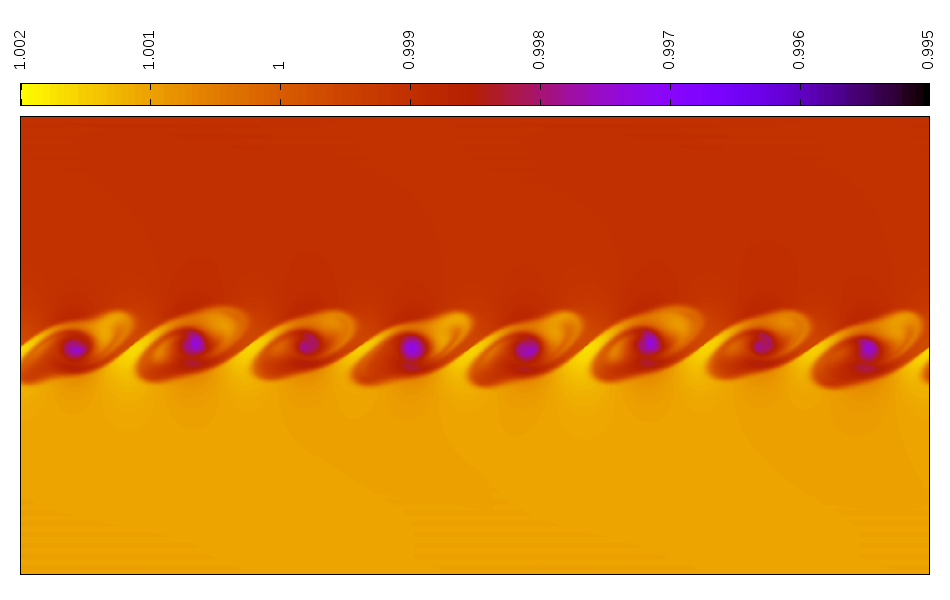} \hfill \includegraphics[width=0.43\textwidth]{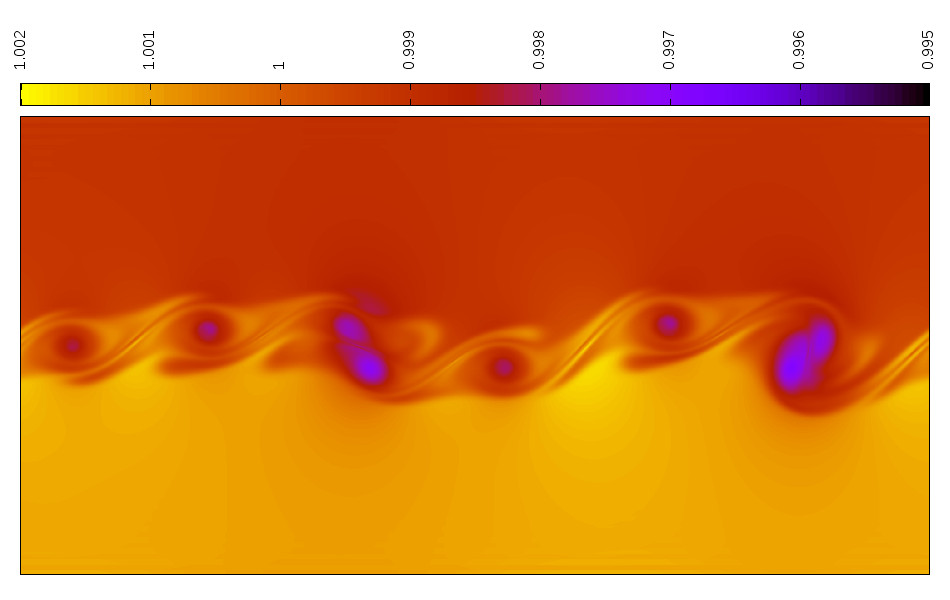}\\
 \includegraphics[width=0.43\textwidth]{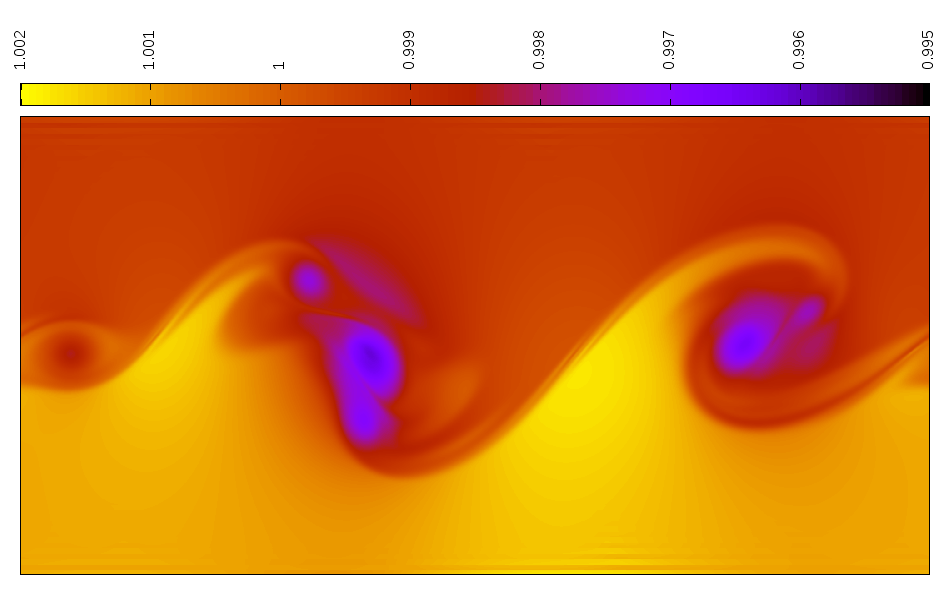} \hfill \includegraphics[width=0.43\textwidth]{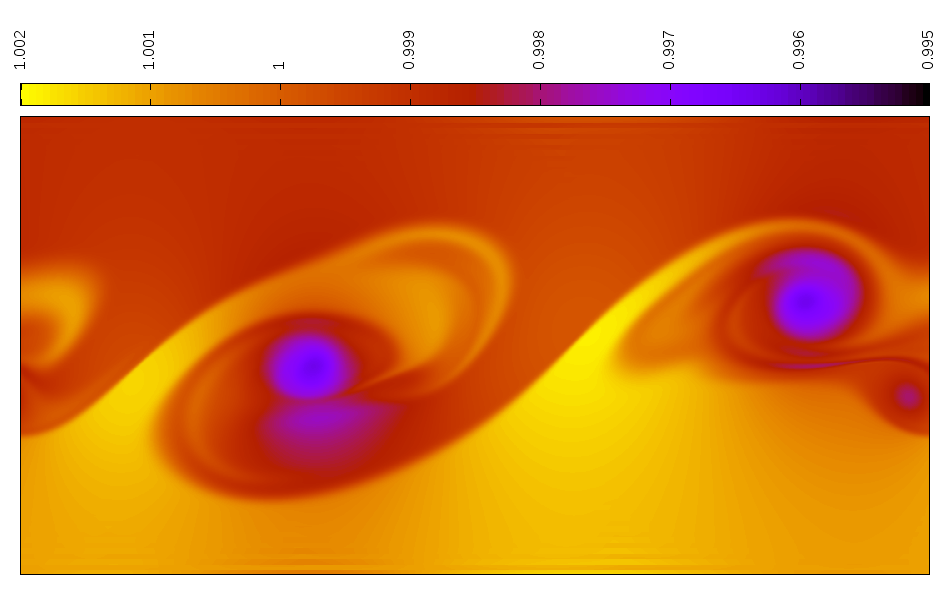}\\
 \includegraphics[width=0.43\textwidth]{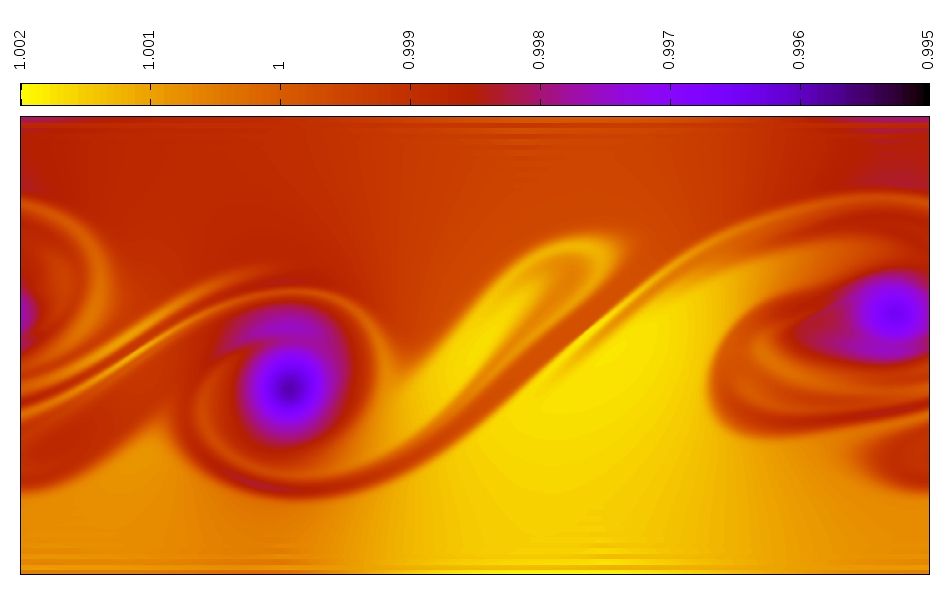} \hfill \includegraphics[width=0.43\textwidth]{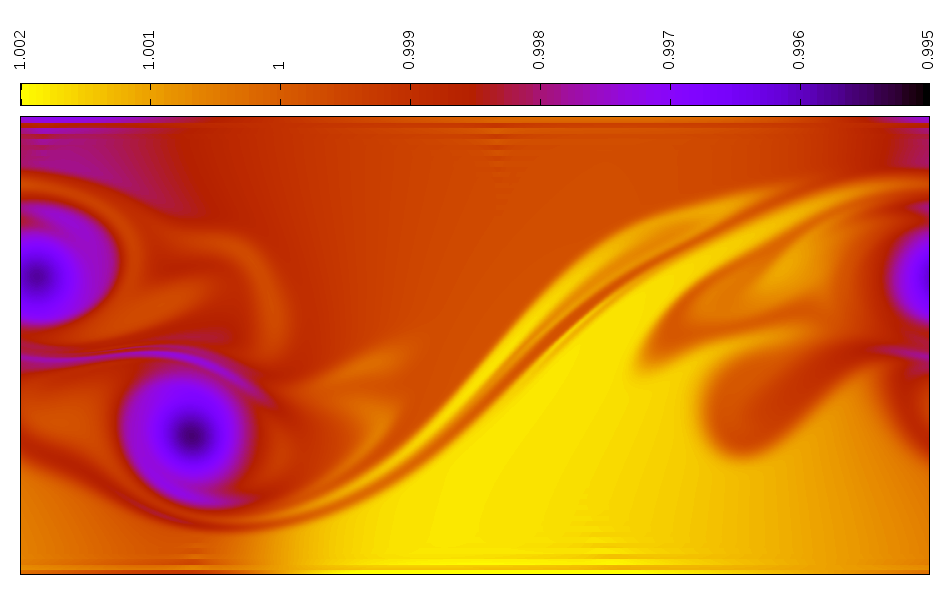}\\
 \includegraphics[width=0.43\textwidth]{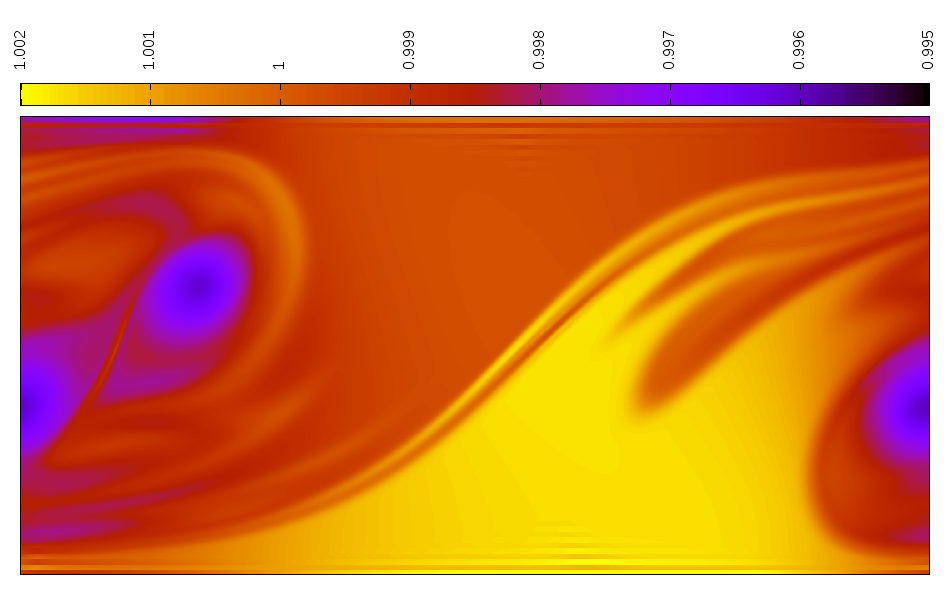} \hfill \includegraphics[width=0.43\textwidth]{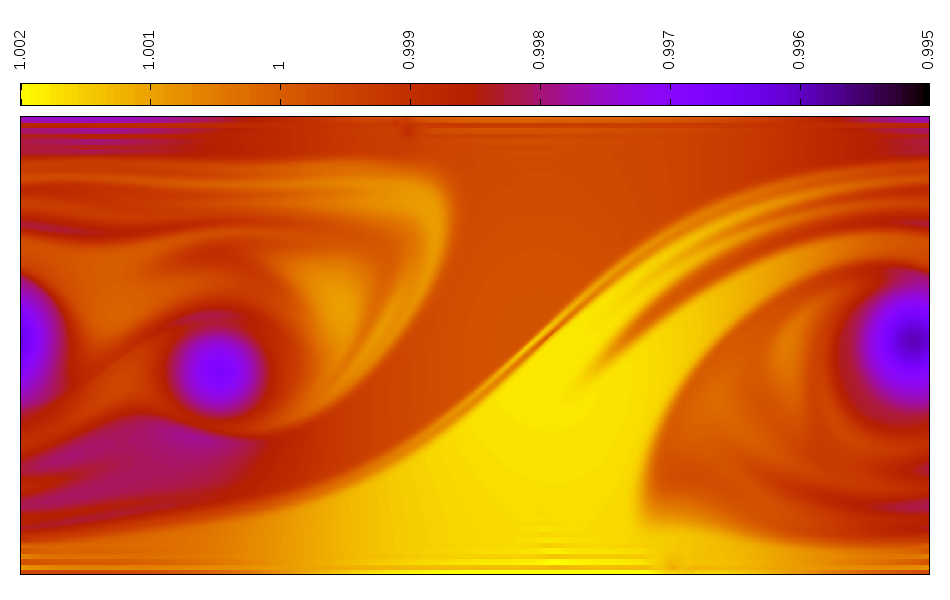}\\
 \includegraphics[width=0.43\textwidth]{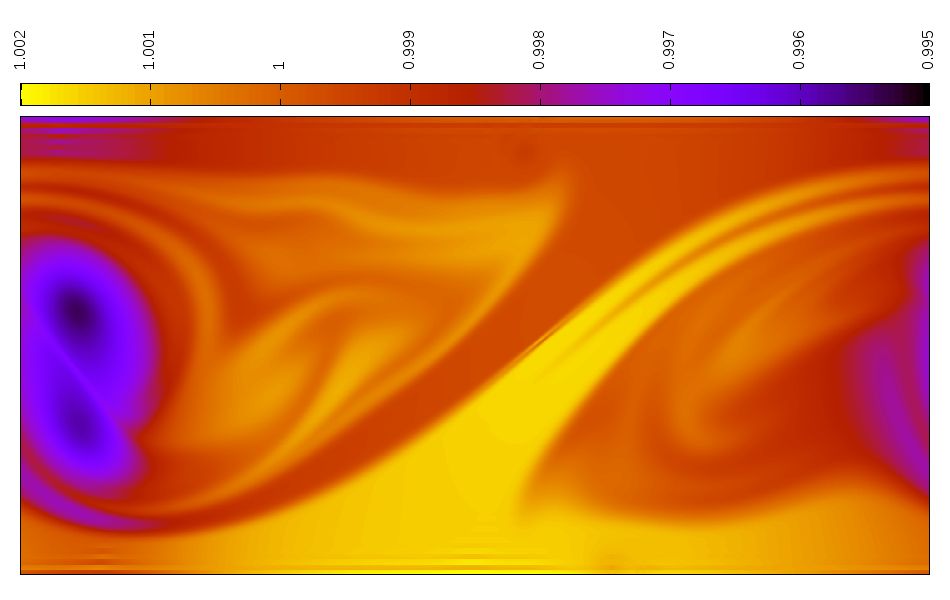} \hfill \includegraphics[width=0.43\textwidth]{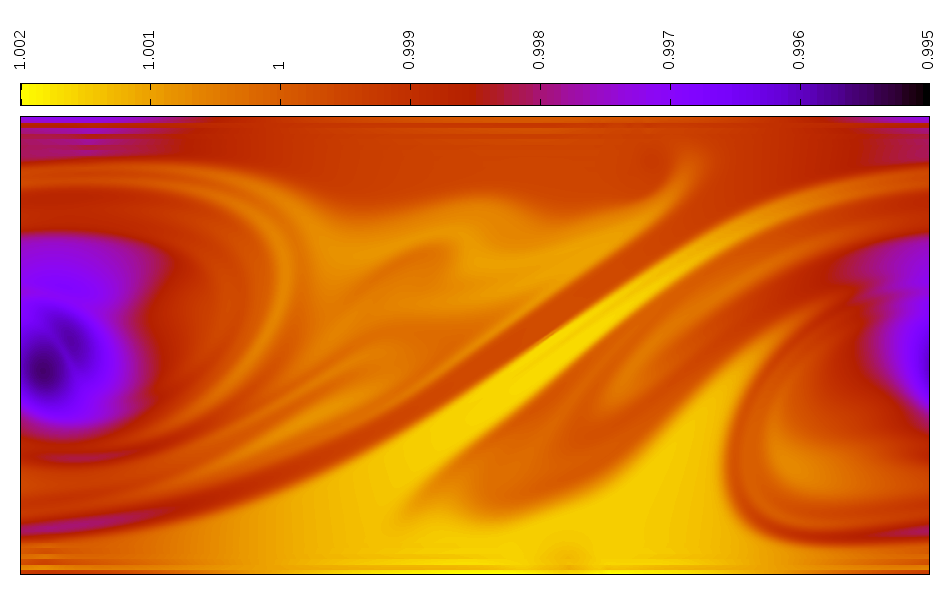}
 \caption{Kelvin-Helmholtz instability solved with Scheme \ref{scheme:euler} with $\text{CFL} = 0.7$ on a grid of $2000 \times 1000$ cells covering $[0,2] \times [0,1]$ with periodic (left and right) and fixed boundaries (top and bottom). Color coded is the density which acts approximately as a passive scalar in this setup. \emph{Left to right, top to bottom}: $t = 5, 10, 15, \ldots, 50$.}
 \label{fig:kh}
\end{figure}

Finally, a Kelvin-Helmholtz instability is studied. The initial setup is that of a contact discontinuity between a fluid with density $1.001$ and $0.999$ streaming in opposite directions with speed $\pm0.1$ in $x$-direction. The perpendicular velocity is $10^{-3} \sin(2 \pi x)$ initially, and acts as a perturbation. The pressure is uniformly equal to $5$. Results of the simulation are shown in Figure \ref{fig:kh}. Despite the method being only first-order accurate, the numerical results display a stunning amount of detail. This is only possible because of low Mach compliance, which reduces the numerical diffusion; in order to resolve these features with a standard method one would require much more resolved grids, entailing significantly longer computational time.

\section{Conclusions and outlook}

The present paper aims at outlining connections between structure-preserving numerical methods on staggered grids and collocated methods. The main result is the insight that a sequential explicit (leap-frog-type) time integrator is the essential ingredient. In use for a long time already in the setting of staggered grids, or Hamiltonian problems, it is applied in the present work to collocated methods without any reference to a Hamiltonian structure. The sequential explicit time integrator is understood simply as a way to stabilize central derivatives while using an essentially time-explicit scheme. Central derivatives are well-known to lead to structure preservation and low Mach compliance and much of the effort of deriving low Mach number compliant schemes can be understood as identifying ways to stabilize a central derivative, i.e. reducing numerical diffusion.

This paper additionally is a contribution to the question how a numerical method for linear acoustics can be extended to include advection, and ultimately handle the nonlinear Euler equations (Scheme \ref{scheme:euler}). First, it is shown that numerical methods for Maxwell's equations can be used to solve linear acoustics, and that these methods easily can be vorticity- and stationarity-preserving upon usage of a sequential explicit time integrator. This relies on these systems of PDEs being off-diagonal, in a way made precise in the text. Advection, however, is a diagonal operator, and in this work a simple explicit time integration is suggested to account for it. It is argued that the compressive terms play an equally important role and that it is not sufficient to discretize them with central derivatives alone. Taking inspiration from Lagrange-Projection methods, a particular discretization of these terms is shown to improve stability, when these three ingredients are combined to a method for the full Euler equations.

Finally, this paper presents a selection of novel staggered-grid and collocated methods for Maxwell's equations based on Yee's method, in particular a new way of staggering the electric and magnetic fields that leads to an improved CFL constraint (Scheme \ref{scheme:yeeextended}). It is a staggered-grid method inspired by collocated methods for acoustics, i.e. this time the transfer of ideas is the other way around.

Contrary to achieving low Mach number compliance with ad hoc fixes, the present work demonstrates a way how central derivatives can be stabilized, while keeping the simplicity of an explicit time integrator. The result is an all-speed method able to deal with both low and high Mach number flow without free parameters. The fact that the method is truly multi-dimensional can be considered a contribution to truly multi-dimensional design principles for numerical methods that avoid excessive diffusion inherited from a dimensionally split approach, which neglects subtle balances between different directions.

Future work will be devoted to extensions of the present method to higher order of accuracy, as well as an analysis of further properties, such as entropy stability or positivity preservation.

\bibliographystyle{alpha}

\begin{thebibliography}{MRKG03}

\bibitem[AG15]{amadori2015}
Debora Amadori and Laurent Gosse.
\newblock {\em Error Estimates for Well-Balanced Schemes on Simple Balance
  Laws: One-Dimensional Position-Dependent Models}.
\newblock BCAM Springer Briefs in Mathematics, Springer, 2015.

\bibitem[AIP19]{abbate19}
Emanuela Abbate, Angelo Iollo, and Gabriella Puppo.
\newblock An asymptotic-preserving all-speed scheme for fluid dynamics and
  nonlinear elasticity.
\newblock {\em SIAM Journal on Scientific Computing}, 41(5):A2850--A2879, 2019.

\bibitem[Bar18]{barsukow18thesis}
Wasilij Barsukow.
\newblock {\em Low {M}ach number finite volume methods for the acoustic and
  {E}uler equations}.
\newblock Doctoral thesis, University of Wuerzburg, 2018.

\bibitem[Bar19]{barsukow17a}
Wasilij Barsukow.
\newblock Stationarity preserving schemes for multi-dimensional linear systems.
\newblock {\em Mathematics of Computation}, 88(318):1621--1645, 2019.

\bibitem[Bar20]{barsukow18hypproceeding}
Wasilij Barsukow.
\newblock Stationary states of finite volume discretizations of
  multi-dimensional linear hyperbolic systems.
\newblock In {\em XVII International Conference on Hyperbolic Problems},
  volume~10, pages 296--303. AIMS Series on Applied Mathematics, 2020.

\bibitem[Bar21a]{barsukow19activeflux}
Wasilij Barsukow.
\newblock The active flux scheme for nonlinear problems.
\newblock {\em Journal of Scientific Computing}, 86(1):1--34, 2021.

\bibitem[Bar21b]{barsukow20cgk}
Wasilij Barsukow.
\newblock Truly multi-dimensional all-speed schemes for the {E}uler equations
  on {C}artesian grids.
\newblock {\em accepted in J. Comp. Phys.}, 2021.

\bibitem[BBT06]{berthon06}
Christophe Berthon, Michael Breu{\ss}, and Marc-Olivier Titeux.
\newblock A relaxation scheme for the approximation of the pressureless euler
  equations.
\newblock {\em Numerical Methods for Partial Differential Equations: An
  International Journal}, 22(2):484--505, 2006.

\bibitem[BDL09]{bouchut09}
Fran{\c{c}}ois Bouchut and T~Morales De~Luna.
\newblock Semi-discrete entropy satisfying approximate {R}iemann solvers. {T}he
  case of the {S}uliciu relaxation approximation.
\newblock {\em Journal of Scientific Computing}, 41(3):483--509, 2009.

\bibitem[BDL{\etalchar{+}}20]{boscheri20}
Walter Boscheri, Giacomo Dimarco, Rapha{\"e}l Loub{\`e}re, Maurizio Tavelli,
  and Marie-H{\'e}l{\`e}ne Vignal.
\newblock A second order all mach number imex finite volume solver for the
  three dimensional euler equations.
\newblock {\em Journal of Computational Physics}, 415:109486, 2020.

\bibitem[BDT21]{boscheri21}
Walter Boscheri, Giacomo Dimarco, and Maurizio Tavelli.
\newblock An efficient second order all {M}ach finite volume solver for the
  compressible {N}avier-{S}tokes equations.
\newblock {\em Computer Methods in Applied Mechanics and Engineering},
  374:113602, 2021.

\bibitem[BEK{\etalchar{+}}17]{barsukow16}
Wasilij Barsukow, Philipp~VF Edelmann, Christian Klingenberg, Fabian Miczek,
  and Friedrich~K R{\"o}pke.
\newblock A numerical scheme for the compressible low-{M}ach number regime of
  ideal fluid dynamics.
\newblock {\em Journal of Scientific Computing}, 72(2):623--646, 2017.

\bibitem[BK22]{barsukow17}
Wasilij Barsukow and Christian Klingenberg.
\newblock Exact solution and a truly multidimensional {G}odunov scheme for the
  acoustic equations.
\newblock {\em ESAIM: M2AN}, 56(1), 2022.

\bibitem[BLMY17]{bispen17}
Georgij Bispen, Maria Lukacova-Medvidova, and Leonid Yelash.
\newblock Asymptotic preserving {IMEX} finite volume schemes for low {M}ach
  number {E}uler equations with gravitation.
\newblock {\em Journal of Computational Physics}, 335:222--248, 2017.

\bibitem[BM05]{birken05}
Philipp Birken and Andreas Meister.
\newblock Stability of preconditioned finite volume schemes at low {M}ach
  numbers.
\newblock {\em BIT Numerical Mathematics}, 45(3):463--480, 2005.

\bibitem[Bou04]{bouchut04}
Fran{\c{c}}ois Bouchut.
\newblock {\em Nonlinear stability of finite Volume Methods for hyperbolic
  conservation laws and Well-Balanced schemes for sources}.
\newblock Springer Science \& Business Media, 2004.

\bibitem[BP21]{boscheri21a}
Walter Boscheri and Lorenzo Pareschi.
\newblock High order pressure-based semi-implicit {IMEX} schemes for the 3{D}
  {N}avier-{S}tokes equations at all {M}ach numbers.
\newblock {\em Journal of Computational Physics}, 434:110206, 2021.

\bibitem[BQRX19]{boscarino19}
Sebastiano Boscarino, Jing-Mei Qiu, Giovanni Russo, and Tao Xiong.
\newblock A high order semi-implicit {IMEX WENO} scheme for the all-{M}ach
  isentropic {E}uler system.
\newblock {\em Journal of Computational Physics}, 392:594--618, 2019.

\bibitem[BW96]{bijl96}
Hester Bijl and Pieter Wesseling.
\newblock A numerical method for the computation of compressible flows with low
  mach number regions.
\newblock {\em Computational Fluid Dynamics}, 96:206--212, 1996.

\bibitem[CCG{\etalchar{+}}10]{chalons10}
Christophe Chalons, Fr{\'e}d{\'e}ric Coquel, Edwige Godlewski, Pierre-Arnaud
  Raviart, and Nicolas Seguin.
\newblock Godunov-type schemes for hyperbolic systems with parameter-dependent
  source: the case of euler system with friction.
\newblock {\em Mathematical Models and Methods in Applied Sciences},
  20(11):2109--2166, 2010.

\bibitem[CCZ14a]{cheng14a}
Yingda Cheng, Andrew~J Christlieb, and Xinghui Zhong.
\newblock Energy-conserving discontinuous {G}alerkin methods for the
  {V}lasov-{A}mpere system.
\newblock {\em Journal of Computational Physics}, 256:630--655, 2014.

\bibitem[CCZ14b]{cheng14}
Yingda Cheng, Andrew~J Christlieb, and Xinghui Zhong.
\newblock Energy-conserving discontinuous {G}alerkin methods for the
  {V}lasov-{M}axwell system.
\newblock {\em Journal of Computational Physics}, 279:145--173, 2014.

\bibitem[CDK12]{cordier12}
Floraine Cordier, Pierre Degond, and Anela Kumbaro.
\newblock An asymptotic-preserving all-speed scheme for the {E}uler and
  {N}avier-{S}tokes equations.
\newblock {\em Journal of Computational Physics}, 231(17):5685--5704, 2012.

\bibitem[CFL28]{courant28}
Richard Courant, Kurt Friedrichs, and Hans Lewy.
\newblock {\"U}ber die partiellen differenzengleichungen der mathematischen
  {P}hysik.
\newblock {\em Mathematische Annalen}, 100(1):32--74, 1928.

\bibitem[CG84]{casulli84}
V~Casulli and D~Greenspan.
\newblock Pressure method for the numerical solution of transient, compressible
  fluid flows.
\newblock {\em International Journal for Numerical Methods in Fluids},
  4(11):1001--1012, 1984.

\bibitem[CGK16]{chalons16}
Christophe Chalons, Mathieu Girardin, and Samuel Kokh.
\newblock An all-regime lagrange-projection like scheme for the gas dynamics
  equations on unstructured meshes.
\newblock {\em Communications in Computational Physics}, 20(1):188--233, 2016.

\bibitem[CGLM14]{cheng14b}
Yingda Cheng, Irene~M Gamba, Fengyan Li, and Philip~J Morrison.
\newblock Discontinuous {G}alerkin methods for the {V}lasov-{M}axwell
  equations.
\newblock {\em SIAM Journal on Numerical Analysis}, 52(2):1017--1049, 2014.

\bibitem[CN47]{crank47}
John Crank and Phyllis Nicolson.
\newblock A practical method for numerical evaluation of solutions of partial
  differential equations of the heat-conduction type.
\newblock In {\em Mathematical Proceedings of the Cambridge Philosophical
  Society}, volume~43, pages 50--67. Cambridge University Press, 1947.

\bibitem[DBTF19]{dumbser19}
Michael Dumbser, Dinshaw~S Balsara, Maurizio Tavelli, and Francesco Fambri.
\newblock A divergence-free semi-implicit finite volume scheme for ideal,
  viscous, and resistive magnetohydrodynamics.
\newblock {\em International Journal for Numerical Methods in Fluids},
  89(1-2):16--42, 2019.

\bibitem[DC16]{dumbser16}
Michael Dumbser and Vincenzo Casulli.
\newblock A conservative, weakly nonlinear semi-implicit finite volume scheme
  for the compressible {N}avier-{S}tokes equations with general equation of
  state.
\newblock {\em Applied Mathematics and Computation}, 272:479--497, 2016.

\bibitem[Del10]{dellacherie10}
St{\'e}phane Dellacherie.
\newblock Analysis of {G}odunov type schemes applied to the compressible
  {E}uler system at low {M}ach number.
\newblock {\em Journal of Computational Physics}, 229(4):978--1016, 2010.

\bibitem[DJOR16]{dellacherie16}
St{\'e}phane Dellacherie, Jonathan Jung, Pascal Omnes, and P-A Raviart.
\newblock Construction of modified {G}odunov-type schemes accurate at any
  {M}ach number for the compressible {E}uler system.
\newblock {\em Mathematical Models and Methods in Applied Sciences},
  26(13):2525--2615, 2016.

\bibitem[DJY07]{degond07}
Pierre Degond, S~Jin, and J~Yuming.
\newblock {M}ach-number uniform asymptotic-preserving gauge schemes for
  compressible flows.
\newblock {\em Bulletin-Institute of Mathematics Academia Sinica}, 2(4):851,
  2007.

\bibitem[DLV17]{dimarco17}
Giacomo Dimarco, Rapha{\"e}l Loub{\`e}re, and Marie-H{\'e}l{\`e}ne Vignal.
\newblock Study of a new asymptotic preserving scheme for the {E}uler system in
  the low {M}ach number limit.
\newblock {\em SIAM Journal on Scientific Computing}, 39(5):A2099--A2128, 2017.

\bibitem[DT11]{degond11}
Pierre Degond and Min Tang.
\newblock All speed scheme for the low mach number limit of the isentropic
  euler equations.
\newblock {\em Communications in Computational Physics}, 10(1):1--31, 2011.

\bibitem[Ebi77]{ebin77}
David~G Ebin.
\newblock The motion of slightly compressible fluids viewed as a motion with
  strong constraining force.
\newblock {\em Annals of mathematics}, pages 141--200, 1977.

\bibitem[ER13]{eymann13}
Timothy~A Eymann and Philip~L Roe.
\newblock Multidimensional active flux schemes.
\newblock In {\em 21st AIAA computational fluid dynamics conference}, 2013.

\bibitem[FLLP05]{fezoui05}
Loula Fezoui, St{\'e}phane Lanteri, St{\'e}phanie Lohrengel, and Serge Piperno.
\newblock Convergence and stability of a discontinuous {G}alerkin time-domain
  method for the 3{D} heterogeneous {M}axwell equations on unstructured meshes.
\newblock {\em ESAIM: Mathematical Modelling and Numerical Analysis},
  39(6):1149--1176, 2005.

\bibitem[GG86a]{guerra86a}
Jaime Guerra and Bertil Gustafsson.
\newblock A numerical method for incompressible and compressible flow problems
  with smooth solutions.
\newblock {\em Journal of Computational Physics}, 63(2):377--397, 1986.

\bibitem[GG86b]{guerra86}
Jaime Guerra and Bertil Gustafsson.
\newblock A semi-implicit method for hyperbolic problems with different
  time-scales.
\newblock {\em SIAM journal on numerical analysis}, 23(4):734--749, 1986.

\bibitem[Gir14]{girardin14}
Mathieu Girardin.
\newblock {\em Asymptotic preserving and all-regime Lagrange-Projection like
  numerical schemes: application to two-phase flows in low Mach regime}.
\newblock PhD thesis, 2014.

\bibitem[Gus87]{gustafsson87}
Bertil Gustafsson.
\newblock Unsymmetric hyperbolic systems and the {E}uler equations at low
  {M}ach numbers.
\newblock {\em Journal of scientific computing}, 2(2):123--136, 1987.

\bibitem[GV99]{guillard99}
Herv{\'e} Guillard and C{\'e}cile Viozat.
\newblock On the behaviour of upwind schemes in the low {M}ach number limit.
\newblock {\em Computers \& fluids}, 28(1):63--86, 1999.

\bibitem[HA71]{harlow71}
Francis~H Harlow and Anthony~A Amsden.
\newblock A numerical fluid dynamics calculation method for all flow speeds.
\newblock {\em Journal of Computational Physics}, 8(2):197--213, 1971.

\bibitem[HJL12]{haack12}
Jeffrey Haack, Shi Jin, and Jian-Guo Liu.
\newblock An all-speed asymptotic-preserving method for the isentropic {E}uler
  and {N}avier-{S}tokes equations.
\newblock {\em Communications in Computational Physics}, 12(4):955--980, 2012.

\bibitem[HW65]{harlow65}
Francis~H Harlow and J~Eddie Welch.
\newblock Numerical calculation of time-dependent viscous incompressible flow
  of fluid with free surface.
\newblock {\em The physics of fluids}, 8(12):2182--2189, 1965.

\bibitem[JT06]{jeltsch06}
Rolf Jeltsch and Manuel Torrilhon.
\newblock On curl-preserving finite volume discretizations for shallow water
  equations.
\newblock {\em BIT Numerical Mathematics}, 46(1):35--53, 2006.

\bibitem[Kim04]{kim04}
Cheolwan Kim.
\newblock Higher-order upwind leapfrog methods for multi-dimensional acoustic
  equations.
\newblock {\em International Journal for Numerical Methods in Fluids},
  44(5):505--523, 2004.

\bibitem[KKL{\etalchar{+}}06]{koh06}
Il-Suek Koh, Hyun Kim, Jung-Mi Lee, Jong-Gwan Yook, and Chang~Sung Pil.
\newblock Novel explicit 2-d {FDTD} scheme with isotropic dispersion and
  enhanced stability.
\newblock {\em IEEE transactions on antennas and propagation},
  54(11):3505--3510, 2006.

\bibitem[Kle95]{klein95}
Rupert Klein.
\newblock Semi-implicit extension of a {G}odunov-type scheme based on low
  {M}ach number asymptotics {I}: {O}ne-dimensional flow.
\newblock {\em Journal of Computational Physics}, 121(2):213--237, 1995.

\bibitem[KM81]{klainerman81}
Sergiu Klainerman and Andrew Majda.
\newblock Singular limits of quasilinear hyperbolic systems with large
  parameters and the incompressible limit of compressible fluids.
\newblock {\em Communications on Pure and Applied Mathematics}, 34(4):481--524,
  1981.

\bibitem[KP89]{karki89}
KC~Karki and SV~Patankar.
\newblock Pressure based calculation procedure for viscous flows at all
  speedsin arbitrary configurations.
\newblock {\em AIAA journal}, 27(9):1167--1174, 1989.

\bibitem[LeV02]{leveque02}
Randall~J LeVeque.
\newblock {\em Finite volume methods for hyperbolic problems}, volume~31.
\newblock {C}ambridge {U}niversity {P}ress, 2002.

\bibitem[LG08]{li08}
Xue-song Li and Chun-wei Gu.
\newblock An all-speed {R}oe-type scheme and its asymptotic analysis of low
  {M}ach number behaviour.
\newblock {\em Journal of Computational Physics}, 227(10):5144--5159, 2008.

\bibitem[LG13]{li13}
Xue-song Li and Chun-wei Gu.
\newblock Mechanism of {R}oe-type schemes for all-speed flows and its
  application.
\newblock {\em Computers \& Fluids}, 86:56--70, 2013.

\bibitem[LR56]{lax56}
Peter~D Lax and Robert~D Richtmyer.
\newblock Survey of the stability of linear finite difference equations.
\newblock {\em Communications on pure and applied mathematics}, 9(2):267--293,
  1956.

\bibitem[LW64]{lax64}
Peter~D Lax and Burton Wendroff.
\newblock Difference schemes for hyperbolic equations with high order of
  acuracy.
\newblock {\em Comm. Pure Appl. Math.}, 17:381--398, 1964.

\bibitem[Mil71]{miller71}
John~JH Miller.
\newblock On the location of zeros of certain classes of polynomials with
  applications to numerical analysis.
\newblock {\em IMA Journal of Applied Mathematics}, 8(3):397--406, 1971.

\bibitem[MR01]{morton01}
Keith~William Morton and Philip~L Roe.
\newblock Vorticity-preserving {L}ax-{W}endroff-type schemes for the system
  wave equation.
\newblock {\em SIAM Journal on Scientific Computing}, 23(1):170--192, 2001.

\bibitem[MRE15]{miczek15}
F~Miczek, FK~R{\"o}pke, and PVF Edelmann.
\newblock New numerical solver for flows at various {M}ach numbers.
\newblock {\em Astronomy \& Astrophysics}, 576:A50, 2015.

\bibitem[MRKG03]{munz03}
C-D Munz, Sabine Roller, Rupert Klein, and Karl~J Geratz.
\newblock The extension of incompressible flow solvers to the weakly
  compressible regime.
\newblock {\em Computers \& Fluids}, 32(2):173--196, 2003.

\bibitem[MS01]{metivier01}
Guy M{\'e}tivier and Steve Schochet.
\newblock The incompressible limit of the non-isentropic {E}uler equations.
\newblock {\em Archive for rational mechanics and analysis}, 158(1):61--90,
  2001.

\bibitem[MT09]{mishra09preprint}
Siddhartha Mishra and Eitan Tadmor.
\newblock Constraint preserving schemes using potential-based fluxes {II}.
  genuinely multi-dimensional central schemes for systems of conservation laws.
\newblock {\em ETH preprint}, (2009-32), 2009.

\bibitem[OSB{\etalchar{+}}16]{birken16}
Kai O{\ss}wald, Alexander Siegmund, Philipp Birken, Volker Hannemann, and
  Andreas Meister.
\newblock L2roe: a low dissipation version of {R}oe's approximate {R}iemann
  solver for low {M}ach numbers.
\newblock {\em International Journal for Numerical Methods in Fluids},
  81(2):71--86, 2016.

\bibitem[PJSS14]{pinto14}
Martin~Campos Pinto, S{\'e}bastien Jund, St{\'e}phanie Salmon, and Eric
  Sonnendr{\"u}cker.
\newblock Charge-conserving fem--pic schemes on general grids.
\newblock {\em Comptes Rendus Mecanique}, 342(10-11):570--582, 2014.

\bibitem[PM05]{park05}
JH~Park and C-D Munz.
\newblock Multiple pressure variables methods for fluid flow at all {M}ach
  numbers.
\newblock {\em International journal for numerical methods in fluids},
  49(8):905--931, 2005.

\bibitem[PRF02]{piperno02}
Serge Piperno, Malika Remaki, and Loula Fezoui.
\newblock A nondiffusive finite volume scheme for the three-dimensional
  {M}axwell's equations on unstructured meshes.
\newblock {\em SIAM Journal on Numerical Analysis}, 39(6):2089--2108, 2002.

\bibitem[Rem00]{remaki99}
Malika Remaki.
\newblock A new finite volume scheme for solving {M}axwell’s system.
\newblock {\em COMPEL-The international journal for computation and mathematics
  in electrical and electronic engineering}, 2000.

\bibitem[Rie11]{rieper11}
Felix Rieper.
\newblock A low-{M}ach number fix for {R}oe’s approximate {R}iemann solver.
\newblock {\em Journal of Computational Physics}, 230(13):5263--5287, 2011.

\bibitem[RLM15]{roe15}
Philip~L Roe, Tyler Lung, and Jungyeoul Maeng.
\newblock New approaches to limiting.
\newblock In {\em 22nd AIAA Computational Fluid Dynamics Conference}, page
  2913, 2015.

\bibitem[RM00]{roller00}
Sabine Roller and Claus-Dieter Munz.
\newblock A low mach number scheme based on multi-scale asymptotics.
\newblock {\em Computing and Visualization in Science}, 3(1):85--91, 2000.

\bibitem[Roe98]{roe98}
Philip Roe.
\newblock Linear bicharacteristic schemes without dissipation.
\newblock {\em SIAM Journal on Scientific Computing}, 19(5):1405--1427, 1998.

\bibitem[Roe17]{roe17}
Philip Roe.
\newblock Multidimensional upwinding.
\newblock {\em Handbook of Numerical Analysis}, 18:53--80, 2017.

\bibitem[RW01]{rodrigue01}
Garry Rodrigue and Daniel White.
\newblock A vector finite element time-domain method for solving {M}axwell's
  equations on unstructured hexahedral grids.
\newblock {\em SIAM Journal on Scientific Computing}, 23(3):683--706, 2001.

\bibitem[Sch17]{schur17}
Issai Schur.
\newblock {\"U}ber {P}otenzreihen, die im {I}nnern des {E}inheitskreises
  beschr{\"a}nkt sind.
\newblock {\em Journal f{\"u}r die reine und angewandte {M}athematik},
  147:205--232, 1917.

\bibitem[Sch18]{schur18}
Issai Schur.
\newblock {\"U}ber {P}otenzreihen, die im {I}nnern des {E}inheitskreises
  beschr{\"a}nkt sind. part ii.
\newblock {\em J. Reine Angew. Math}, 148:122--145, 1918.

\bibitem[SCS92]{shyy92}
Wei Shyy, Ming-Hsiung Chen, and Chia-Sheng Sun.
\newblock Pressure-based multigrid algorithm for flow at all speeds.
\newblock {\em AIAA journal}, 30(11):2660--2669, 1992.

\bibitem[Sod78]{sod78}
Gary~A Sod.
\newblock A survey of several finite difference methods for systems of
  nonlinear hyperbolic conservation laws.
\newblock {\em Journal of computational physics}, 27(1):1--31, 1978.

\bibitem[TB75]{taflove75}
Allen Taflove and Morris~E Brodwin.
\newblock Numerical solution of steady-state electromagnetic scattering
  problems using the time-dependent {M}axwell's equations.
\newblock {\em IEEE transactions on microwave theory and techniques},
  23(8):623--630, 1975.

\bibitem[TD08]{thornber08}
BJR Thornber and D~Drikakis.
\newblock Numerical dissipation of upwind schemes in low {M}ach flow.
\newblock {\em International journal for numerical methods in fluids},
  56(8):1535--1541, 2008.

\bibitem[TPK20]{thomann20}
Andrea Thomann, Gabriella Puppo, and Christian Klingenberg.
\newblock An all speed second order well-balanced imex relaxation scheme for
  the euler equations with gravity.
\newblock {\em Journal of Computational Physics}, 420:109723, 2020.

\bibitem[TR93]{thomas93}
Jeffrey Thomas and Philip Roe.
\newblock Development of non-dissipative numerical schemes for computational
  aeroacoustics.
\newblock In {\em 11th Computational Fluid Dynamics Conference}, page 3382,
  1993.

\bibitem[Tur87]{turkel87}
Eli Turkel.
\newblock Preconditioned methods for solving the incompressible and low speed
  compressible equations.
\newblock {\em Journal of computational physics}, 72(2):277--298, 1987.

\bibitem[VBW11]{viallet11}
M~Viallet, I~Baraffe, and R~Walder.
\newblock Towards a new generation of multi-dimensional stellar evolution
  models: development of an implicit hydrodynamic code.
\newblock {\em Astronomy \& Astrophysics}, 531:A86, 2011.

\bibitem[WS95]{weiss95}
Jonathan~M Weiss and Wayne~A Smith.
\newblock Preconditioning applied to variable and constant density flows.
\newblock {\em AIAA journal}, 33(11):2050--2057, 1995.

\bibitem[WSW02]{wenneker02}
I~Wenneker, A~Segal, and P~Wesseling.
\newblock A mach-uniform unstructured staggered grid method.
\newblock {\em International Journal for Numerical Methods in Fluids},
  40(9):1209--1235, 2002.

\bibitem[Yee66]{yee66}
Kane Yee.
\newblock Numerical solution of initial boundary value problems involving
  {M}axwell's equations in isotropic media.
\newblock {\em IEEE Transactions on antennas and propagation}, 14(3):302--307,
  1966.

\end{thebibliography}
\newcommand{\etalchar}[1]{$^{#1}$}

\appendix
\section{The discrete Fourier transform of the different extensions of the Yee scheme} \label{app:fourier}

Considering linear numerical methods on Cartesian grids in Fourier space has been found advantageous, not only to study stability, but also involution preservation / stationarity preservation, see \cite{barsukow17a}. On Cartesian grids the discrete Fourier transform reads
\begin{align}
 q_{ij}^n = \hat q^n \exp(\ii k_x \Delta x i + \ii k_y \Delta y j)
\end{align}
$\vec k = (k_x, k_y)$ is the wave number, and $\hat q^n$, or $\hat q^n(\vec k)$ is the discrete (spatial) Fourier transform of $q_{ij}^n$. Shifts $i \mapsto i+1$ are replaced by multiplications with $t_x := \exp(\ii k_x \Delta x)$, and similarly $j \mapsto j+1$ is conveyed by the translation operator $t_y := \exp(\ii k_y \Delta y)$. For example, a finite difference such as
\begin{align}
\frac{q^n_{i+1}-  q^n_{i-1}}{2 \Delta x} 
\end{align}
then translates into
\begin{align}
 \hat q^n \frac{t_x - t_x^{-1}}{2 \Delta x} = \hat q^n \frac{(t_x+1)(t_x-1)}{2 t_x \Delta x}
\end{align}

\subsection{Original scheme and its collocated interpretation}

Renaming of variables does not alter the Fourier transform of the scheme, so that both Schemes \ref{scheme:yeecollocated} and \ref{scheme:yeeoriginal} become
\begin{align}
 \frac{(\widehat B^z)^{n+1} - (\widehat B^z)^{n}}{\Delta t} &= -\left( (\widehat E^y)^n\frac{t_x  - 1}{\Delta x} - (\widehat E^x)^n\frac{t_y  - 1}{\Delta y}  \right ) \label{eq:yeecolloc1fourier}\\
 \frac{(\widehat E^x)^{n+1} - (\widehat E^x)^{n}}{\Delta t} &= (\widehat B^z)^{n+1}\frac{t_y - 1}{t_y\Delta y} \label{eq:yeecolloc2fourier}\\
 \frac{(\widehat E^y)^{n+1} - (\widehat E^y)^{n}}{\Delta t} &= -(\widehat B^z)^{n+1} \frac{t_x - 1}{t_x\Delta x}  \label{eq:yeecolloc3fourier}
\end{align}

In the language of the Fourier transform, this scheme keeps stationary the dimensionally split curl
\begin{align}
 (\widehat E^y)^n\frac{t_x  - 1}{\Delta x} - (\widehat E^x)^n\frac{t_y  - 1}{\Delta y}
\end{align}
as is readily verified upon explicit computation.

The stationarity preserving scheme from \cite{barsukow17a} (when applied to Maxwell's equations instead of linear acoustics) would keep stationary the truly multi-dimensional curl
\begin{align}
 (\widehat E^y)^n\frac{t_x  - 1}{\Delta x} \frac{t_y+1}{2} - (\widehat E^x)^n\frac{t_y  - 1}{\Delta y} \frac{t_x+1}{2} \label{eq:firstchoice}
\end{align}
or possibly
\begin{align}
 (\widehat E^y)^n\frac{t_x  - 1}{\Delta x} \frac{t_y+1}{2t_y} - (\widehat E^x)^n\frac{t_y  - 1}{\Delta y} \frac{t_x+1}{2t_x} \label{eq:secondchoice}
\end{align}

\subsection{Multi-dimensional extension of the Yee scheme} \label{ssec:fourieryeeextension}

The first choice \eqref{eq:firstchoice} can be achieved by modifying the Yee scheme as follows:
\begin{align}
 \frac{(\widehat B^z)^{n+1} - (\widehat B^z)^{n}}{\Delta t} &= -\left( (\widehat E^y)^n\frac{t_x  - 1}{\Delta x} \frac{t_y+1}{2} - (\widehat E^x)^n\frac{t_y  - 1}{\Delta y} \frac{t_x+1}{2}  \right ) \label{eq:yeeextended1fourier}\\
 \frac{(\widehat E^x)^{n+1} - (\widehat E^x)^{n}}{\Delta t} &= (\widehat B^z)^{n+1}\frac{t_y - 1}{t_y\Delta y} \frac{t_x+1}{2t_x} \label{eq:yeeextended2fourier}\\
 \frac{(\widehat E^y)^{n+1} - (\widehat E^y)^{n}}{\Delta t} &= -(\widehat B^z)^{n+1} \frac{t_x - 1}{t_x\Delta x} \frac{t_y+1}{2t_y}  \label{eq:yeeextended3fourier}
\end{align}

This gives \eqref{eq:yeecollocextended1}--\eqref{eq:yeecollocextended3} / Scheme \ref{scheme:yeeextended}.

The second choice \eqref{eq:secondchoice} leads to
\begin{align}
 \frac{(\widehat B^z)^{n+1} - (\widehat B^z)^{n}}{\Delta t} &= -\left( (\widehat E^y)^n\frac{t_x  - 1}{\Delta x} \frac{t_y+1}{2t_y} - (\widehat E^x)^n\frac{t_y  - 1}{\Delta y} \frac{t_x+1}{2t_x}  \right ) \\
 \frac{(\widehat E^x)^{n+1} - (\widehat E^x)^{n}}{\Delta t} &= (\widehat B^z)^{n+1}\frac{t_y - 1}{t_y\Delta y} \frac{t_x+1}{2} \\
 \frac{(\widehat E^y)^{n+1} - (\widehat E^y)^{n}}{\Delta t} &= -(\widehat B^z)^{n+1} \frac{t_x - 1}{t_x\Delta x} \frac{t_y+1}{2}  
\end{align}
i.e.
\begin{align}
 \frac{(B^z)^{n+1}_{ij} - (B^z)^{n}_{ij}}{\Delta t} &= -\left( \frac{(E^y)^n_{i+1,j} - (E^y)^n_{ij}  + (E^y)^n_{i+1,j-1}  - (E^y)^n_{i,j-1}}{2 \Delta x} \right . \\\nonumber &\quad\quad\quad - \left. \frac{(E^x)^n_{i,j+1}  - (E^x)^n_{ij}  + (E^x)^n_{i-1,j+1}  - (E^x)^n_{i-1,j}}{2  \Delta y}   \right ) \\
 \frac{(E^x)^{n+1}_{ij} - (E^x)^{n}_{ij}}{\Delta t} &= \frac{(B^z)^{n+1}_{i+1,j} - (B^z)^{n+1}_{i+1,j-1} + (B^z)^{n+1}_{ij}  - (B^z)^{n+1}_{i,j-1}}{2\Delta y} \\
 \frac{(E^y)^{n+1}_{ij} - (E^y)^{n}_{ij}}{\Delta t} &= -\frac{(B^z)^{n+1}_{i,j+1} - (B^z)^{n+1}_{i-1,j+1} + (B^z)^{n+1}_{ij} - (B^z)^{n+1}_{i-1,j} }{2 \Delta x} 
\end{align}
Here, from the first equation it seems natural to rename $(B^z)^{n}_{ij} \mapsto (B^z)^{n+1}_{i+\frac12,j-\frac12}$, but this still does not give a symmetric scheme:
\begin{align}
 &\frac{(B^z)^{n+1}_{i+\frac12,j-\frac12} - (B^z)^{n}_{i+\frac12,j-\frac12}}{\Delta t} = -\left( \frac{(E^y)^n_{i+1,j} - (E^y)^n_{ij}  + (E^y)^n_{i+1,j-1}  - (E^y)^n_{i,j-1}}{2 \Delta x} \right . \\ \nonumber & \phantom{mmmmmmmmmmmmmm} - \left. \frac{(E^x)^n_{i,j+1}  - (E^x)^n_{ij}  + (E^x)^n_{i-1,j+1}  - (E^x)^n_{i-1,j}}{2  \Delta y}   \right ) \\
 &\frac{(E^x)^{n+1}_{ij} - (E^x)^{n}_{ij}}{\Delta t} = \frac{(B^z)^{n+1}_{i+\frac32,j-\frac12} - (B^z)^{n+1}_{i+\frac32,j-\frac32} + (B^z)^{n+1}_{i+\frac12,j-\frac12}  - (B^z)^{n+1}_{i+\frac12,j-\frac32}}{2\Delta y} \\
 &\frac{(E^y)^{n+1}_{ij} - (E^y)^{n}_{ij}}{\Delta t} = -\frac{(B^z)^{n+1}_{i+\frac12,j+\frac12} - (B^z)^{n+1}_{i-\frac12,j+\frac12} + (B^z)^{n+1}_{i+\frac12,j-\frac12} - (B^z)^{n+1}_{i-\frac12,j-\frac12} }{2 \Delta x} 
\end{align}
Observe how, e.g. in the second equation the right hand-side is centered around $(i+1,j-1)$ instead of $(i,j)$.

\subsection{Collocated sequential explicit schemes with central derivatives}

The Fourier transform of \eqref{eq:centralextendedacoustics1}--\eqref{eq:centralextendedacoustics3} (Scheme \ref{scheme:centralextended}) for linear acoustics reads

\begin{align}
 \frac{\hat p^{n+1} - \hat p^{n}}{\Delta t} &= -c^2\left( \frac{(t_x-1)(t_x+1)}{2t_x \Delta x} \boxed{\frac{(t_y+1)^2}{4 t_y}} \hat u^n - \frac{(t_y-1)(t_y+1)}{2t_y\Delta y} \boxed{\frac{(t_x+1)^2}{4 t_x}} \hat v^n \right ) \label{eq:centralextended1fourier}\\ 
 \frac{\hat u^{n+1} - \hat u^{n}}{\Delta t} &= -\frac{ (t_x - 1)(t_x+1)}{2t_x\Delta x} \boxed{\frac{(t_y+1)^2}{4 t_y}} \hat p^{n+1} \\
 \frac{\hat v^{n+1} - \hat v^{n}}{\Delta t} &= -\frac{(t_y - 1)(t_y+1)}{2 t_y\Delta y} \boxed{\frac{(t_x+1)^2}{4 t_y}} \hat p^{n+1} \label{eq:centralextended3fourier}
\end{align}
where terms responsible for the perpendicular averaging are boxed.

\section{Stability of sequential explicit schemes for the Maxwell equations} \label{app:stability}

Consider the Fourier transform of a sequential explicit scheme for the Maxwell equations:
\begin{align}
 \frac{(B^z)^{n+1} - (B^z)^{n}}{\Delta t} &= -\left( (E^y)^n \mathscr D_x - (E^x)^n \mathscr D_y  \right ) \label{eq:stabilitygeneral1} \\
 \frac{(E^x)^{n+1} - (E^x)^{n}}{\Delta t} &= (B^z)^{n+1} \mathscr D'_y \\
 \frac{(E^y)^{n+1} - (E^y)^{n}}{\Delta t} &= -(B^z)^{n+1} \mathscr D'_x  \label{eq:stabilitygeneral3}
\end{align}
with complex numbers $\mathscr D_x,\mathscr D'_x,\mathscr D_y,\mathscr D'_y$. For example, the truly-multi-dimensional extension \eqref{eq:yeeextended1}--\eqref{eq:yeeextended3} / Scheme \ref{scheme:yeeextended}, with its Fourier transform given by \eqref{eq:yeeextended1fourier}--\eqref{eq:yeeextended3fourier}, yields
\begin{align}
 \mathscr D_x &= \frac{t_x  - 1}{\Delta x} \frac{t_y+1}{2} & \mathscr D_y &= \frac{t_y  - 1}{\Delta y} \frac{t_x+1}{2} \\
 \mathscr D'_x &= \frac{t_x - 1}{t_x\Delta x} \frac{t_y+1}{2t_y} & \mathscr D'_y &= \frac{t_y - 1}{t_y\Delta y} \frac{t_x+1}{2t_x}
\end{align}

Converting \eqref{eq:stabilitygeneral1}--\eqref{eq:stabilitygeneral3} into a fully explicit scheme yields
\begin{align}
 (B^z)^{n+1}  &=  (B^z)^{n}-\Delta t\left( (E^y)^n \mathscr D_x - (E^x)^n \mathscr D_y  \right ) \\
 (E^x)^{n+1}  &= (E^x)^{n} + \Delta t(B^z)^{n}\mathscr D'_y -\Delta t^2\left( (E^y)^n \mathscr D_x\mathscr D'_y - (E^x)^n \mathscr D_y\mathscr D'_y  \right )  \\
 (E^y)^{n+1}   &=(E^y)^{n} - \Delta t(B^z)^{n}\mathscr D'_x +\Delta t^2\left( (E^y)^n \mathscr D_x\mathscr D'_x - (E^x)^n \mathscr D_y \mathscr D'_x \right )   
\end{align}
i.e.
\begin{align}
 \veccc{B^z}{E^x}{E^y}^{n+1} &= \left( \begin{array}{ccc}  
                                         1 & \Delta t  \mathscr D_y &-\Delta t  \mathscr D_x    \\
                                           \Delta t \mathscr D'_y & 1  +\Delta t^2  \mathscr D_y\mathscr D'_y & -\Delta t^2  \mathscr D_x\mathscr D'_y \\
                                          - \Delta t\mathscr D'_x & -\Delta t^2  \mathscr D_y \mathscr D'_x  & 1  +\Delta t^2  \mathscr D_x\mathscr D'_x
                                       \end{array} \right )\veccc{B^z}{E^x}{E^y}^{n} 
\end{align}
The characteristic polynomial of this matrix reads

\begin{align}
 (1-z) \left( (1-z)^2   -z \Delta t^2
 \Big( \mathscr D_x\mathscr D'_x +\mathscr D_y\mathscr D'_y    \Big )      \right ) \label{eq:charpolyfull}
\end{align}

One of the eigenvalues therefore is always 1, which by results of \cite{barsukow17a} yields
\begin{theorem} 
 The numerical method \eqref{eq:stabilitygeneral1}--\eqref{eq:stabilitygeneral3} is stationarity preserving / involution preserving for any choice of $D_x, D'_x, D_y, D'_y$. 
\end{theorem}

For stability, the other eigenvalues are of interest. Consider the following Theorem from \cite{miller71}, which is based on results from \cite{schur17,schur18}:

\begin{theorem} \label{thm:millerstability}
 Given a non-constant complex-valued polynomial $f(z) = \sum_{j=0}^n a_j z^j \in P^n$, and construct
 \begin{align}
  f^*(z) &:= \sum_{j=0}^n  \bar a_{n-j} z^j \in P^{n} &
  f_1(z) &:= \frac{f^*(0) f(z) - f(0) f^*(z)}{z} \in P^{n-1}
 \end{align}
 The complex conjugate of $a \in \mathbb C$ is denoted by $\bar a$, and $f' \in P^{n-1}$ is the derivative of $f$ with respect to $z$. Then, all the zeros of $f$ are contained in the closed unit disc iff either
 \begin{itemize}
  \item $|f^*(0)| > |f(0)|$ and all the zeros of $f_1$ are contained in the closed unit disc (or $f_1$ constant), or
  \item $f_1(z) \equiv 0$ and all the zeros of $f'$ are contained in the closed unit disc (or $f'$ constant).
 \end{itemize}
\end{theorem}
Observe that $f_1$ and $f'$ are of smaller degree than $f$ and thus the algorithm surely terminates.

Here, from \eqref{eq:charpolyfull},
\begin{align}
 f(z) &= 1 - z \left( 2 + \Delta t^2 \Big( \mathscr D_x\mathscr D'_x +\mathscr D_y\mathscr D'_y \right ) \Big) + z^2\\
 f^*(z) &= f(z) \qquad \Rightarrow \qquad f_1(z) \equiv 0\\
 f'(z) &= - \left( 2 + \Delta t^2 \Big( \mathscr D_x\mathscr D'_x +\mathscr D_y\mathscr D'_y \right )\Big) + 2z\\
\end{align}
and performing the same analysis with $f'$ yields as stability condition
\begin{align}
 1 > \left | 1 + \frac12 \Delta t^2 \Big( \mathscr D_x\mathscr D'_x +\mathscr D_y\mathscr D'_y \Big ) \right |
\end{align}

\begin{corollary}\label{cor:yeeextendedstability}
 The multi-dimensional extension of the Yee scheme (Schemes \ref{scheme:yeecollocatedextended}/\ref{scheme:yeeextended}) is stable for CFL $< 1$ if $\Delta y = \Delta x$.
\end{corollary}
\begin{proof}
One finds 
\begin{align}
 \mathscr D_x\mathscr D'_x +\mathscr D_y\mathscr D'_y &= \frac{t_x  - 1}{\Delta x} \frac{t_y+1}{2}  \frac{t_x - 1}{t_x\Delta x} \frac{t_y+1}{2t_y}   + \frac{t_y  - 1}{\Delta y} \frac{t_x+1}{2} \frac{t_y - 1}{t_y\Delta y} \frac{t_x+1}{2t_x} \\
 &= \frac{(t_x  - 1)^2}{\Delta x^2 t_x} \frac{(t_y+1)^2}{4 t_y}  + \frac{(t_y  - 1)^2}{t_y\Delta y^2} \frac{(t_x+1)^2}{4 t_x} \\
 &= \frac{t_x^2 - 2 t_x + 1}{\Delta x^2 t_x} \frac{t_y^2 + 2 t_y +1}{4 t_y}  + \frac{t_y^2 - 2 t_y + 1}{t_y\Delta y^2} \frac{t_x^2 + 2 t_x+1}{4 t_x} \\
\intertext{Insert $t_x = \ee^{\ii \beta_x}$, $t_y = \ee^{\ii \beta_y}$}
 &= \frac{ (\cos \beta_x - 1) (\cos \beta_y + 1) }{\Delta x^2}  + \frac{(\cos \beta_y - 1)(\cos \beta_x + 1)}{\Delta y^2}   \\
 &= 2\frac{\cos \beta_x \cos \beta_y   - 1 }{\Delta x^2}
\end{align}
The stability condition thus reads
\begin{align}
 1 > \left | 1 + \frac{\Delta t^2}{\Delta x^2} (\cos \beta_x \cos \beta_y   - 1) \right |
\end{align}
for all $\beta_x, \beta_y$. Thus
\begin{align}
 1 &> \left | 1 -2 \frac{\Delta t^2}{\Delta x^2}  \right |\\
 1 &> \frac{\Delta t}{\Delta x}
\end{align}

\end{proof}

\begin{corollary}\label{cor:centralstability}
 The collocated sequential explicit scheme with central derivatives (Scheme \ref{scheme:central}) is stable for CFL $< 1$ if $\Delta y = \Delta x$.
\end{corollary}
\begin{proof}
For Scheme \ref{scheme:central}
\begin{align}
 \mathscr D_x = \mathscr D'_x &= \frac{(t_x  - 1)(t_x+1)}{2 t_x\Delta x} & \mathscr D_y = \mathscr D'_y &= \frac{(t_y  - 1)(t_y+1)}{2 t_y\Delta y} 
\end{align}
and thus (again, with $\Delta y = \Delta x$)
\begin{align}
 \mathscr D_x\mathscr D'_x +\mathscr D_y\mathscr D'_y &= \frac{(t_x  - 1)^2(t_x+1)^2}{4 t_x^2\Delta x^2} + \frac{(t_y  - 1)^2(t_y+1)^2}{4 t_y^2\Delta y^2} \\
 &=\frac{(\cos^2 \beta_x - 1 ) + (\cos^2 \beta_y   -1 )}{\Delta x^2}
\end{align}
The stability condition amounts to
\begin{align}
 1 &> \left | 1 - \frac12  \frac{ \Delta t^2}{\Delta x^2} (2 - \cos^2 \beta_x - \cos^2 \beta_y) \right |\\
 -1 &< 1 - \frac{ \Delta t^2}{\Delta x^2}
\end{align}
Thus, $\sqrt{2} > \frac{\Delta t}{\Delta x}$, in agreement with \cite{remaki99}.

\end{proof}

\begin{corollary} \label{cor:centralextendedstability}
 The multi-dimensionally extended version of the collocated sequential explicit scheme with central derivatives (Scheme \ref{scheme:centralextended}) is stable for CFL $< 2$ if $\Delta y = \Delta x$.
\end{corollary}
\begin{proof}
Indeed, for Scheme \ref{scheme:centralextended}
\begin{align}
 \mathscr D_x = \mathscr D'_x &= \frac{(t_x  - 1)(t_x+1)}{2 t_x\Delta x} \frac{(t_y+1)^2}{4t_y} \\
 \mathscr D_y = \mathscr D'_y &= \frac{(t_y  - 1)(t_y+1)}{2 t_y\Delta y} \frac{(t_x+1)^2}{4 t_x}
\end{align}
With $\Delta y = \Delta x$ one finds
\begin{align}
 \mathscr D_x\mathscr D'_x +\mathscr D_y\mathscr D'_y &= \frac{(t_x  - 1)^2(t_x+1)^2}{4 t_x^2\Delta x^2} \frac{(t_y+1)^4}{16 t_y^2} + \frac{(t_y  - 1)^2 (t_y+1)^2}{4 t_y^2\Delta x^2} \frac{(t_x+1)^4}{16 t_x^2}\\
 &=\frac{( \cos^2 \beta_x  - 1 )( \cos \beta_y + 1 )^2 + ( \cos^2 \beta_y -1 ) ( \cos \beta_x + 1 )^2}{4\Delta x^2}\\
 &=( \cos \beta_x + 1 )( \cos \beta_y +1 )\frac{ \cos \beta_x \cos \beta_y  -1}{2\Delta x^2}
\end{align}
The image of $( \cos \beta_x + 1 )( \cos \beta_y +1 )(\cos \beta_x \cos \beta_y  -1)$ is $[-2,0]$. Therefore, the image of
\begin{align}
 1 + \frac12 \Delta t^2 \Big( \mathscr D_x\mathscr D'_x +\mathscr D_y\mathscr D'_y \Big ) = 1 + \frac14 \frac{\Delta t^2}{\Delta x^2} ( \cos \beta_x + 1 )( \cos \beta_y +1 )( \cos \beta_x \cos \beta_y  -1)
\end{align}
is $\displaystyle \left[1 - \frac12 \frac{\Delta t^2}{\Delta x^2} , 1 \right]$ and stability requires
\begin{align}
 -1 &< 1 - \frac12 \frac{\Delta t^2}{\Delta x^2}\\
 2 &>  \frac{\Delta t}{\Delta x}
\end{align}

\end{proof}

This section is concluded by rederiving the well-known stability condition for the original Yee scheme (Scheme \ref{scheme:yeecollocated}/\ref{scheme:yeeoriginal}):

\begin{corollary} \label{cor:yeeoriginalstability}
 The Yee scheme (Scheme \ref{scheme:yeecollocated}/\ref{scheme:yeeoriginal}) in two spatial dimensions is stable for CFL $< 1/\sqrt{2}$ if $\Delta y = \Delta x$.
\end{corollary}
\begin{proof}

\begin{align}
 \mathscr D_x &= \frac{t_x  - 1}{\Delta x} & \mathscr D_y &= \frac{t_y  - 1}{\Delta y}\\
 \mathscr D'_x &= \frac{t_x - 1}{t_x\Delta x} & \mathscr D'_y &= \frac{t_y - 1}{t_y\Delta y}
\end{align}
\begin{align}
 \mathscr D_x\mathscr D'_x +\mathscr D_y\mathscr D'_y  &= \frac{(t_x  - 1)^2}{\Delta x^2 t_x} + \frac{(t_y  - 1)^2}{\Delta y^2 t_y}\\
 &= 2\frac{ \cos \beta_x  -1}{\Delta x^2 } + 2\frac{\cos \beta_y  -1 }{\Delta y^2 }
\end{align}
With $\Delta y = \Delta x$ one thus finds
\begin{align}
 1 > \left | 1 + \frac{\Delta t^2}{\Delta x^2} ( \cos \beta_x  +\cos \beta_y  -2)  \right |
\end{align}
Thus,
\begin{align}
 \frac12 >  \frac{\Delta t^2}{\Delta x^2}
\end{align}

\end{proof}

\section{The extended Yee scheme for the Maxwell equations in three spatial dimensions} \label{app:maxwell3d}

It is easiest to use the Fourier transform to construct the three-dimensional analogue of the extension leading from the original Yee scheme \ref{scheme:yeeoriginal} to its extended version \ref{scheme:yeeextended}, the latter employing staggered grids of cell centers and nodes only. Comparing to \eqref{eq:yeeextended1fourier}, the natural extension is
\begin{align}
 \frac{(\widehat B^x)^{n+1} - (\widehat B^x)^{n}}{\Delta t} &= -\left( (\widehat E^z)^n\frac{t_y  - 1}{\Delta y} \boxed{ \frac{t_x+1}{2}  \frac{t_z+1}{2} } - (\widehat E^y)^n\frac{t_z  - 1}{\Delta z} \boxed{\frac{t_x+1}{2}  \frac{t_y+1}{2} } \right ) \label{eq:yeeextended3d1fourier}\\
 \frac{(\widehat B^y)^{n+1} - (\widehat B^y)^{n}}{\Delta t} &= -\left( (\widehat E^x)^n\frac{t_z  - 1}{\Delta z} \boxed{\frac{t_x+1}{2}  \frac{t_y+1}{2} }- (\widehat E^z)^n\frac{t_x  - 1}{\Delta x} \boxed{ \frac{t_y+1}{2}  \frac{t_z+1}{2} } \right ) \\
 \frac{(\widehat B^z)^{n+1} - (\widehat B^z)^{n}}{\Delta t} &= -\left( (\widehat E^y)^n\frac{t_x  - 1}{\Delta x} \boxed{ \frac{t_y+1}{2}  \frac{t_z+1}{2}} - (\widehat E^x)^n\frac{t_y  - 1}{\Delta y} \boxed{ \frac{t_x+1}{2}  \frac{t_z+1}{2} } \right ) \\
 \frac{(\widehat E^x)^{n+1} - (\widehat E^x)^{n}}{\Delta t} &= (\widehat B^z)^{n+1}\frac{t_y - 1}{t_y\Delta y} \boxed{\frac{t_x+1}{2t_x} \frac{t_z+1}{2t_z}} -  (\widehat B^y)^{n+1}\frac{t_z - 1}{t_z\Delta z} \boxed{\frac{t_x+1}{2t_x} \frac{t_y+1}{2t_y}}\\
 \frac{(\widehat E^y)^{n+1} - (\widehat E^y)^{n}}{\Delta t} &=  (\widehat B^x)^{n+1}\frac{t_z - 1}{t_z\Delta z} \boxed{\frac{t_x+1}{2t_x} \frac{t_y+1}{2t_y}} -(\widehat B^z)^{n+1} \frac{t_x - 1}{t_x\Delta x} \boxed{\frac{t_y+1}{2t_y} \frac{t_z+1}{2t_z}}\\
 \frac{(\widehat E^z)^{n+1} - (\widehat E^z)^{n}}{\Delta t} &=  (\widehat B^y)^{n+1}\frac{t_x - 1}{t_x\Delta x} \boxed{\frac{t_y+1}{2t_y} \frac{t_z+1}{2t_z}} -(\widehat B^x)^{n+1} \frac{t_y - 1}{t_y\Delta y} \boxed{\frac{t_x+1}{2t_x} \frac{t_z+1}{2t_z}} \label{eq:yeeextended3d3fourier}
\end{align}
while, without the boxed terms, this would be the Fourier transform of the original Yee scheme from \cite{yee66} in three dimensions.

Associate now the electric fields $E^x, E^y, E^z$ to the cell centers $(i,j,k)$ and the magnetic fields to the corners/nodes $(i+\frac12,j+\frac12,k+\frac12)$. Upon the inverse Fourier transform, a term such as
\begin{align}
 (\widehat E^z)^n\frac{t_y  - 1}{\Delta y} \boxed{ \frac{t_x+1}{2}  \frac{t_z+1}{2} }
\end{align}
becomes 
\begin{align*}
 \frac{ 1 }{4 \Delta y} &\left ( (E^z)^n_{i+1,j+1,k+1} - (E^z)^n_{i+1,j,k+1} + (E^z)^n_{i,j+1,k+1} - (E^z)^n_{i,j,k+1} \right. \\  &\quad \left. + (E^z)^n_{i+1,j+1,k} - (E^z)^n_{i+1,j,k} + (E^z)^n_{i,j+1,k} - (E^z)^n_{i,j,k} \right )
\end{align*}
Using notation introduced in Section \ref{ssec:allspeedeulermethod}, this finite difference formula can be written more concisely as
\begin{align}
 \frac{   \{[\{ (E^z)^n \}_{i+\frac12} ]_{j+\frac12} \}_{k+\frac12} }{4 \Delta y}
\end{align}
making obvious that this expression is now centered at a node $(i+\frac12,j+\frac12,k+\frac12)$. This fits well together with the fact that this term appears in the update of a magnetic field centered at the same location. The complete scheme is
\begin{scheme}[Yee extended, 3D] \label{scheme:yeeextended3d}
 \begin{align}
 &\frac{( B^x)^{n+1}_{i+\frac12,j+\frac12,k+\frac12} - ( B^x)^{n}_{i+\frac12,j+\frac12,k+\frac12}}{\Delta t} = \label{eq:yeeextended3d1}\\ &\phantom{mmmmmmmmm} \nonumber -\left(   \frac{ [ \{ \{ ( E^z)^n \}_{i+\frac12}  ]_{j+\frac12} \}_{k+\frac12}  }{4\Delta y}  - \frac{\{ [ \{ ( E^y)^n \}_{i+\frac12} ]_{k+\frac12}  \}_{k+\frac12} }{4\Delta z}  \right )  \\
 &\frac{( B^y)^{n+1}_{i+\frac12,j+\frac12,k+\frac12} - ( B^y)^{n}_{i+\frac12,j+\frac12,k+\frac12}}{\Delta t} = \\ &\phantom{mmmmmmmmm} \nonumber -\left( \frac{ [ \{ \{ ( E^x)^n \}_{i+\frac12} \}_{j+\frac12} ]_{k+\frac12} }{4\Delta z} - \frac{  \{\{ [ ( E^z)^n ]_{i+\frac12} \}_{j+\frac12} \}_{k+\frac12}  }{4\Delta x} \right ) \\
 &\frac{( B^z)^{n+1}_{i+\frac12,j+\frac12,k+\frac12} - ( B^z)^{n}_{i+\frac12,j+\frac12,k+\frac12}}{\Delta t} = \\ &\phantom{mmmmmmmmm} \nonumber -\left( \frac{ \{\{[ ( E^y)^n ]_{i+\frac12} \}_{j+\frac12} \}_{k+\frac12} }{4\Delta x}  -  \frac{ \{ [ \{ ( E^x)^n \}_{i+\frac12} ]_{j+\frac12} \}_{k+\frac12} }{4\Delta y}  \right ) \\
 &\frac{( E^x)^{n+1}_{ijk} - ( E^x)^{n}_{ijk}}{\Delta t} = \frac{ \{ [ \{ ( B^z)^{n+1}\}_{i\pm\frac12} ]_{j\pm\frac12} \}_{k\pm\frac12}  }{4\Delta y}  -  \frac{ [\{\{  ( B^y)^{n+1} \}_{i\pm\frac12} \}_{j\pm\frac12} ]_{k\pm\frac12}   }{4\Delta z} \\
 &\frac{( E^y)^{n+1}_{ijk} - ( E^y)^{n}_{ijk}}{\Delta t} =  \frac{ [\{\{  ( B^x)^{n+1} \}_{i\pm\frac12} \}_{j\pm\frac12}  ]_{k\pm\frac12}  }{4\Delta z}  - \frac{ \{\{[  ( B^z)^{n+1} ]_{i\pm\frac12} \}_{j\pm\frac12} \}_{k\pm\frac12}  }{4\Delta x} \\
 &\frac{( E^z)^{n+1}_{ijk} - ( E^z)^{n}_{ijk}}{\Delta t} =  \frac{  \{\{[ ( B^y)^{n+1} ]_{i\pm\frac12} \}_{j\pm\frac12}\}_{k\pm\frac12}  }{4\Delta x}  - \frac{ \{[\{  ( B^x)^{n+1} \}_{i\pm\frac12} ]_{j\pm\frac12} \}_{k\pm\frac12}  }{4\Delta y} \label{eq:yeeextended3d3}
\end{align}
\end{scheme}

The stability of this method can be analyzed along the lines of Section \ref{app:stability}. The characteristic polynomial for the case $\Delta z = \Delta y = \Delta x$ reads
\begin{align}
 \frac{(z-1)^2 \Big (2 \Delta x^2 (z-1)^2 + \mathscr F(\beta_x, \beta_y, \beta_z) \Delta t^2 z \Big)^2}{4 \Delta x^2}
\end{align}
with
\begin{align}
 \mathscr F(\beta_x, \beta_y, \beta_z) &= 3 + \cos \beta_x + \cos \beta_y +\cos \beta_z - \cos \beta_x \cos \beta_y - \cos \beta_x \cos \beta_z - \cos \beta_y \cos \beta_z \nonumber \\ &- 3 \cos \beta_x \cos \beta_y \cos \beta_z
\end{align}
Theorem \ref{thm:millerstability} amounts to studying the polynomial
\begin{align}
  z - 1 + \frac14 \mathscr F(\beta_x, \beta_y, \beta_z) \frac{\Delta t^2}{ \Delta x^2} 
\end{align}
with the stability condition
\begin{align}
 \left| 1 - \frac14 \mathscr F(\beta_x, \beta_y, \beta_z) \frac{\Delta t^2}{ \Delta x^2}  \right | \leq 1
\end{align}
The range of $\mathscr F$ is $[0, 8]$, and thus the range of this term is $ \left[ 1 , 1 - 2 \frac{\Delta t^2}{ \Delta x^2} \right ] $. The CFL condition therefore is
\begin{align}
 \frac{\Delta t}{ \Delta x} < 1
\end{align}

\section{The Lagrange-Projection method} \label{app:lagproj}

For the equation
\begin{align}
 \del_t q + \del_x (U(x) q) &= 0
\end{align}
consider the quantity contained in $[x_{i-\frac12} + \xi_1(t), x_{i+\frac12} + \xi_2(t)]$, with $\xi_1(0)  = \xi_2(0)= 0$:
\begin{align}
 \del_t \left(\int\limits_{x_{i-\frac12} + \xi_1(t)}^{x_{i+\frac12} + \xi_2(t)} \dd x \,q(t, x)\right ) &= q\Big(t, x_{i+\frac12} + \xi_2(t)\Big)\xi_2'(t) - q\Big(t, x_{i-\frac12} + \xi_1(t)\Big) \xi_1'(t) \\&\phantom{mmmmmmmmmm} \nonumber + \int\limits_{x_{i-\frac12} + \xi_1(t)}^{x_{i+\frac12} + \xi_2(t)} \dd x \,\del_t q(t, x)\\
 &\!\!\!\!\!\!\!\!\!\!\!\!\!\!\!\!\!\!\!\!\!\!\!\!\!\!\!\!\!\!\!\!\!\!\!\!\!\!\!\!\!\!\!\!\!\!\!\!\!\!\!\!\!\!\!\!\!\!\!\!\!\!\!\!\!\!\!\!\!\!\!\!= q\Big(t, x_{i+\frac12} + \xi_2(t)\Big) \xi_2'(t) - q\Big(t, x_{i-\frac12} + \xi_1(t)\Big) \xi_1'(t)  \\&\!\!\!\!\!\!\!\!\!\!\!\!\!\!\!\!\!\!\!\!\!\!\!\!\!\!\!\!\!\!\!\!\!\!\!\!\!\!\!\!\!\!\!\!\!\!\!\!\!\!\!\!\!\!\!\!\!\!\!\!\!\!\!\!\!+U\Big(x_{i+\frac12} + \xi_2(t)\Big) q\Big(t, x_{i+\frac12} + \xi_2(t)\Big) - U\Big(x_{i-\frac12} + \xi_1(t)\Big) q\Big(t, x_{i-\frac12} + \xi_1(t)\Big) 
\end{align}
The mass inside this time dependent volume is constant if
\begin{align}
 U\Big(x_{i-\frac12} + \xi_1(t)\Big) &= \xi_1'(t)\\
 U\Big(x_{i+\frac12} + \xi_2(t)\Big) &= \xi_2'(t) 
\end{align}
i.e. if the volume is bounded by the characteristics emerging at $x_{i\pm\frac12}$.

As an approximation, consider $\xi_1(t) = x_{i-\frac12} + U(x_{i-\frac12}) t$ and $\xi_2(t) = x_{i+\frac12} + U(x_{i+\frac12}) t$. Assume $q$ to be constant inside $[\xi_1(t), \xi_2(t)]$, and compute the average of $q$ inside $[x_{i-\frac12}, x_{i+\frac12}]$ at time $\Delta t$. There are several cases, depending on the signs of the two velocities.

\begin{figure}
 \centering
 \includegraphics[width=0.7\textwidth]{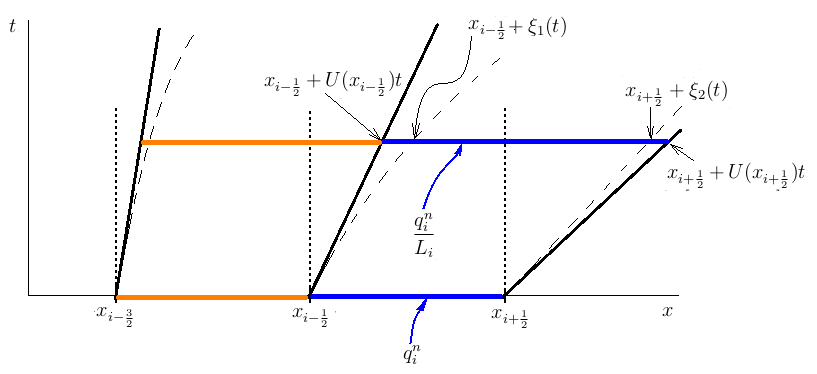}
 \caption{Illustration of the Lagrangian step. The characteristics emerging at the boundaries of two computational cells are shown (dashed lines), together with their approximations (straight lines). At a later time, the average $q_i^n$ in computational cell $[x_{i-\frac12},x_{i+\frac12}]$ has become $q_i^n/ L_i$.}
 \label{fig:lagprojsketch}
\end{figure}

If both velocities are positive($U(x_{i-\frac12}) > 0$, $U(x_{i+\frac12}) > 0$), all the mass $q_i^n \Delta x$ situated initially inside $[x_{i-\frac12}, x_{i+\frac12}]$ is moved to the right and now is distributed over $[x_{i-\frac12} + U(x_{i-\frac12}) \Delta t, x_{i+\frac12} + U(x_{i+\frac12}) \Delta t]$ (see Figure \ref{fig:lagprojsketch}) with a density
 \begin{align}
  \frac{q_i^n \Delta x}{\Delta x + U(x_{i+\frac12}) \Delta t - U(x_{i-\frac12}) \Delta t}
 \end{align}
 The mass still inside $[x_{i-\frac12}, x_{i+\frac12}]$ is 
 \begin{align}
  \frac{q_i^n (\Delta x - U(x_{i-\frac12}) \Delta t) }{1 + \Delta t \frac{U(x_{i+\frac12}) - U(x_{i-\frac12})}{\Delta x}}
 \end{align}
 From the left cell, there is mass flowing into cell $i$. At time $\Delta t$ the total mass flown in from the left is
 \begin{align}
  \frac{q_{i-1}^n U(x_{i-\frac12}) \Delta t }{1 + \Delta t \frac{U(x_{i-\frac12}) - U(x_{i-\frac32})}{\Delta x}}
 \end{align}
 which, on total gives the new average 
 \begin{align}
   q_i^{n+1} &=q_{i-1}^n \frac{ U(x_{i-\frac12}) \frac{\Delta t}{\Delta x} }{L_{i-1}} + q_i^n\frac{ 1 - U(x_{i-\frac12}) \frac{\Delta t}{\Delta x} }{L_i}\\
   &= q_i^n -\frac{\Delta t}{\Delta x} \left( q_i^n\frac{U(x_{i+\frac12})  }{L_i} - q_{i-1}^n \frac{ U(x_{i-\frac12}) }{L_{i-1}} \right ) 
 \end{align}
 having defined 
 \begin{align}
  L_{i} := 1 + \Delta t \frac{U(x_{i+\frac12}) - U(x_{i-\frac12})}{\Delta x}
 \end{align}
 
 The scheme is conservative, with the numerical flux given by
 \begin{align}
  f_{i+\frac12} &= \begin{cases} \displaystyle U(x_{i+\frac12}) \frac{q_i^n}{L_i} & U(x_{i+\frac12}) > 0 \\\\ \displaystyle U(x_{i+\frac12}) \frac{q_{i+1}^n}{L_{i+1}} & U(x_{i+\frac12}) < 0 \end{cases}
  \\ &= \frac{U(x_{i+\frac12})}{2} \left( \frac{q_i^n}{L_i} + \frac{q_{i+1}^n}{L_{i+1}} \right ) - \frac{|U(x_{i+\frac12})|}{2} \left(  \frac{q_{i+1}^n}{L_{i+1}} - \frac{q_i^n}{L_i} \right ) 
 \end{align}

 Observe that in order to account for the compressive terms, the denominators contain an approximation to the cell-centered divergence.

\section{Stability analysis of Scheme \ref{scheme:euler}} \label{app:stabilityeuler}

The choice of including the denominator in \eqref{eq:eulerfluxx} is dictated by the stability requirement of the scheme. Omitting the denominator greatly simplifies the method, and also the proof of the low Mach number compliance (in which case there is only one discrete divergence). However, this comes at the cost of decreased stability as is described at the end of this Section. 

First, a detailed stability analysis of the scheme \ref{scheme:euler} is presented. Consider therefore the linearization of Scheme \ref{scheme:euler} around a constant background $\bar q$. Write $q  = \bar q + \delta q$ in \eqref{eq:eulerfluxxmultid} and neglect terms quadratic in $\delta q$ (again, indices referring to the time discretization are suppressed because they are different for different components of $q$ and $ f_{i+\frac12,j}^x$):
\begin{align}
  f^x_{i+\frac12,j} &\simeq 
 \left( f^x(\bar q) + J^x(\bar q) \frac{ \{\{ \{\delta q\}_{i+\frac12}  \}\}_{j\pm\frac12} }{8} - \frac12 |\bar u | [\delta q]_{i+\frac12,j}  \right ) \times \\ &\phantom{mmmmmm}\nonumber\left( 1 - \Delta t \left(\frac{\{\{ [\delta u]_{i+\frac12}  \}\}_{j\pm\frac12}}{4\Delta x} + \frac{   [\{ \delta v \}_{i+\frac12}]_{j\pm1}  }{4 \Delta y}\right) \right )\\
 &\simeq f^x(\bar q) + J^x(\bar q) \frac{ \{\{ \{\delta q\}_{i+\frac12}  \}\}_{j\pm\frac12} }{8} - \frac12 |\bar u | [\delta q]_{i+\frac12,j} \\
  &\phantom{mmmmmm}\nonumber-
 f^x(\bar q) \Delta t \left(\frac{\{\{ [\delta u]_{i+\frac12}  \}\}_{j\pm\frac12}}{4\Delta x} + \frac{   [\{ \delta v \}_{i+\frac12}]_{j\pm1}  }{4 \Delta y}\right)
\end{align}

Here, $J^x = \nabla_q f^x(q)$. Observe also that $\displaystyle \delta u = \delta \left(\frac{\rho u}{\rho}\right) = \left(-\frac{\bar u}{\bar \rho}, \frac{1}{\bar \rho}, 0, 0\right) \delta q$ and $\displaystyle \delta v = \left(-\frac{\bar v}{\bar \rho}, 0, \frac{1}{\bar \rho}, 0\right) \delta q$.

Writing $q$ instead of $\delta q$ and omitting constant terms which would drop out anyway, the linearization of $f^x_{i+\frac12,j}$ becomes
\begin{align}
  J^x(\bar q) \frac{ \{\{ \{q\}_{i+\frac12}  \}\}_{j\pm\frac12} }{8} - \frac12 |\bar u | [q]_{i+\frac12,j} &- 
 \Delta t f^x(\bar q)  \left(-\frac{\bar u}{\bar \rho}, \frac{1}{\bar \rho}, 0, 0\right)  \frac{\{\{ [q]_{i+\frac12}  \}\}_{j\pm\frac12}}{4\Delta x} \\\nonumber &- \Delta t f^x(\bar q) \left(-\frac{\bar v}{\bar \rho}, 0, \frac{1}{\bar \rho}, 0\right)   \frac{   [\{ q \}_{i+\frac12}]_{j\pm1}  }{4 \Delta y}
\end{align}

The linearization of the $y$-flux $f^y_{i,j+\frac12}$ is
\begin{align}
  J^y(\bar q) \frac{ \{\{ \{q\}\}_{i\pm\frac12}  \}_{j+\frac12} }{8} - \frac12 |\bar v | [q]_{i,j+\frac12} &-
 \Delta t f^y(\bar q)  \left(-\frac{\bar u}{\bar \rho}, \frac{1}{\bar \rho}, 0, 0\right)  \frac{\{ [q]_{i\pm1}  \}_{j+\frac12}}{4\Delta x} \\\nonumber &- \Delta t f^y(\bar q) \left(-\frac{\bar v}{\bar \rho}, 0, \frac{1}{\bar \rho}, 0\right)   \frac{   [\{\{ q \}\}_{i\pm\frac12}]_{j+\frac12}  }{4 \Delta y}
\end{align}
 
The full scheme, upon applying the Fourier transform as in Section \ref{app:fourier} can be written as
\begin{align}
 \hat q^{n+1} = \mathcal A \hat q^n
\end{align}
The characteristic polynomial of the amplification matrix $\mathcal A$ for the case $\bar p = 0$ (i.e. for the purely advective regime) is 
\begin{align}
 (z+z_0)^4
\end{align}
with
\begin{align}
 z_0 &= \frac{2 (|\bar u|+|\bar v|) \Delta t-2 \Delta x-2 \Delta t (|\bar u| \cos \beta_x+|\bar v| \cos \beta_y) }{2 \Delta x} \\
 &+ \nonumber\ii \Delta t \frac{\bar u (1 + \cos \beta_y) \sin \beta_x +  \bar v (1+\cos\beta_x) \sin\beta_y}{2 \Delta x}
\end{align}

\begin{align}
 |z_0|^2 &= \frac{\Big(2 (|\bar u|+|\bar v|) \Delta t-2 \Delta x-2 \Delta t (|\bar u| \cos \beta_x+|\bar v| \cos \beta_y) \Big)^2}{4 \Delta x^2}  \\&+\nonumber \Delta t^2 \frac{\Big( \bar u (1 + \cos \beta_y) \sin \beta_x + \bar v (1+\cos\beta_x) \sin\beta_y\Big )^2}{4 \Delta x^2} 
\end{align}

Its gradient with respect to $\beta_x,\beta_y$ 
\begin{align}
 \frac{2 (|\bar u|+|\bar v|) \Delta t-2 \Delta x-2 \Delta t (|\bar u| \cos \beta_x+|\bar v| \cos \beta_y) }{4 \Delta x^2} \vecc{2 \Delta t |\bar u| \sin \beta_x}{2 \Delta t |\bar v| \sin \beta_y} \\ 
 +\nonumber \Delta t^2 \frac{\Big( \bar u (1 + \cos \beta_y) \sin \beta_x + \bar v (1+\cos\beta_x) \sin\beta_y\Big )}{4 \Delta x^2} \vecc{\bar u (1 + \cos \beta_y) \cos \beta_x - \bar v \sin\beta_x \sin\beta_y}{-\bar u \sin \beta_y \sin \beta_x + \bar v (1+\cos\beta_x) \cos\beta_y}
\end{align}
vanishes wherever $\sin \beta_x = \sin \beta_y = 0$. It can be easily seen that these extrema are maxima if $\beta_x = \pm \pi, \beta_y = \pm \pi$. Based on numerical studies of the function, the conjecture that these maxima are the ones relevant for the stability bound seems reasonable. Then, at these maxima
\begin{align}
 |z_0|^2 = \frac{\Big(4 (|\bar u|+|\bar v|) \Delta t-2 \Delta x  \Big)^2 }{4 \Delta x^2} \leq 1
\end{align}
and the stability bound is
\begin{align}
 (|\bar u|+|\bar v|) \Delta t \leq  \Delta x
\end{align}

On the other hand, the characteristic polynomial of the amplification matrix for $\bar u = \bar v = 0$ is $(z-1)^2$ times
\begin{align}
 z^2 &-z \left( 2 - \frac{\Delta t^2 \bar p  (1-\cos \beta_x \cos \beta_y) (4+\gamma+\gamma \cos\beta_y+\gamma \cos\beta_x +\gamma \cos \beta_x \cos \beta_y)}{2 \Delta x^2 \epsilon^2 \bar \rho}  \right ) \\&\nonumber + \left (1- \frac{ 2 \Delta t^2 \bar p (1-\cos\beta_x \cos\beta_y)}{\Delta x^2 \epsilon^2 \bar \rho} \right )
\end{align}
The prefactor $(z-1)^2$ immediately proves stationarity preservation (with nontrivial stationary states given either by constant pressure and velocity (contact) or by constant pressure and density and divergencefree velocity).
 
Applying Theorem \ref{thm:millerstability} yields first the condition
\begin{align}
 1 \geq \left|  1- \frac{ 2 \Delta t^2 \bar c^2}{\Delta x^2 \epsilon^2 \gamma} (1-\cos\beta_x \cos\beta_y)  \right |
\end{align}
having defined $\bar c^2 = \frac{\gamma \bar p}{\bar \rho}$. If
\begin{align}
 \frac{  \Delta t^2 \bar c^2}{ \epsilon^2 } < \frac{\gamma}{2} \Delta x^2 \label{eq:appaccfl}
\end{align}
then this condition is always true because $1 - \cos \beta_x \cos \beta_y \in [0,2]$. 

The application of Theorem \ref{thm:millerstability} produces next the condition
\begin{align}
 \left| 2- \frac{ 2 \Delta t^2 \bar c^2}{\Delta x^2 \epsilon^2 \gamma} (1-\cos\beta_x \cos\beta_y) \right| > \phantom{mmmmmmmmmmmmmmmmmmmmm}\\\nonumber \left| 2 - \frac{2\Delta t^2 \bar c^2}{\Delta x^2 \epsilon^2 \gamma} \frac{ (1-\cos \beta_x \cos \beta_y) (4+\gamma+\gamma \cos\beta_y+\gamma \cos\beta_x +\gamma \cos \beta_x \cos \beta_y)}{4 }  \right |
\end{align}
On the left side, the absolute value is taken of a positive number, provided \eqref{eq:appaccfl} is fulfilled. The right side is the absolute value of a number between $1-\frac{\gamma}{2} - \frac{1}{2\gamma} < 0$ and 2, because
\begin{align}
 (1-\cos \beta_x \cos \beta_y) (4+\gamma+\gamma \cos\beta_y+\gamma \cos\beta_x +\gamma \cos \beta_x \cos \beta_y) \in \left[0, 4 + 2\gamma + \frac{2}{\gamma} \right]
\end{align}
A maximum is located, for example, at $\beta_x = 0, \cos \beta_y = -\frac{1}{\gamma}$. Therefore the absolute value on the right side is important. In all cases when the absolute value is taken of a positive number, one has
\begin{align}
   1  \leq \frac{  4+\gamma+\gamma \cos\beta_y+\gamma \cos\beta_x +\gamma \cos \beta_x \cos \beta_y}{4 }  
\end{align}
which is always true. In case, $\beta_x$ and $\beta_y$ are such that the absolute value on the right side is taken of a negative number, the condition can only be proven for $\gamma < 2$. It is conjectured though, that this is not sharp, as in the numerical evaluation of these conditions the result seems to hold true for all $\gamma > 1$. Thus, in the purely acoustic case, the maximum CFL number is $\sqrt{\gamma/2}$, i.e. 0.84 for $\gamma = 1.4$. In future one might investigate further the possibility for modification of Scheme \ref{scheme:euler} in order to obtain a $\gamma$-independent CFL number.

\begin{figure}
 \centering
 \includegraphics[width=0.45\textwidth]{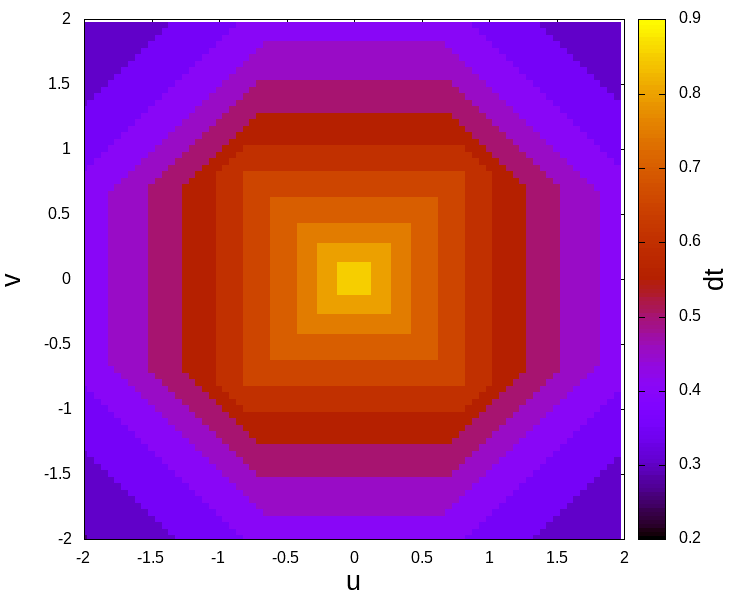} \hfill \includegraphics[width=0.45\textwidth]{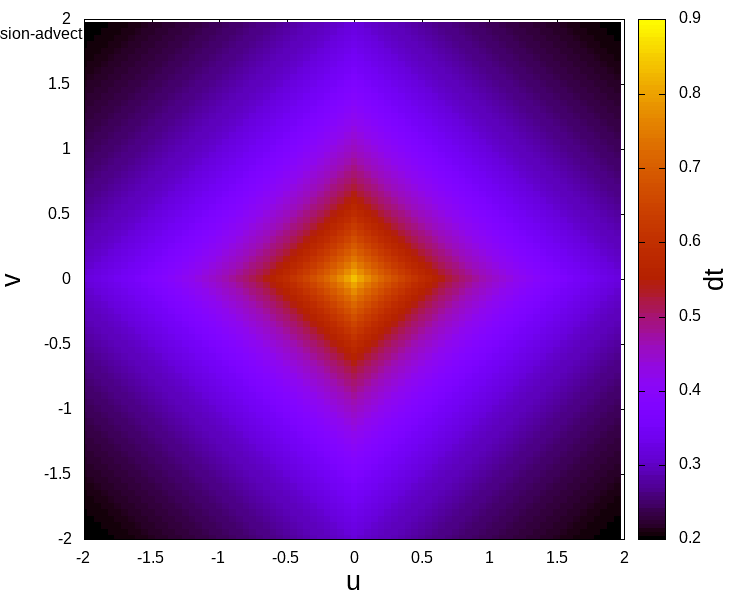}
 \caption{\emph{Left}: Maximal $\frac{\Delta t}{\Delta x}$ according to a numerical evaluation of the criterion from Theorem \ref{thm:millerstability}. The sampling sizes are: for $\beta:x, \beta_y$ 0.01, for $\bar u$, $\bar v$ 0.05, for $\Delta t$ 0.05. $\Delta x$ was set to 1, and $\gamma = 1.4$. \emph{Right}: Maximal $\frac{\Delta t}{\Delta x}$ according to the proposed formula \ref{eq:stabilityconditionproosed}}.
 \label{fig:stabilityeuler}
\end{figure}

For the complete Scheme \ref{scheme:euler} it thus seems reasonable to use a CFL condition of the form
\begin{align}
 \frac{\Delta t}{\Delta x} < \frac{1}{|\bar u| + |\bar v| + \frac{\bar c}{\epsilon} \sqrt{\frac{2}{\gamma}}} \label{eq:stabilityconditionproosed}
\end{align}
This choice is confirmed by numerical studies of the zeros of the characteristic polynomial (see Figure \ref{fig:stabilityeuler}).

\begin{figure}
 \centering
 \includegraphics[width=0.7\textwidth]{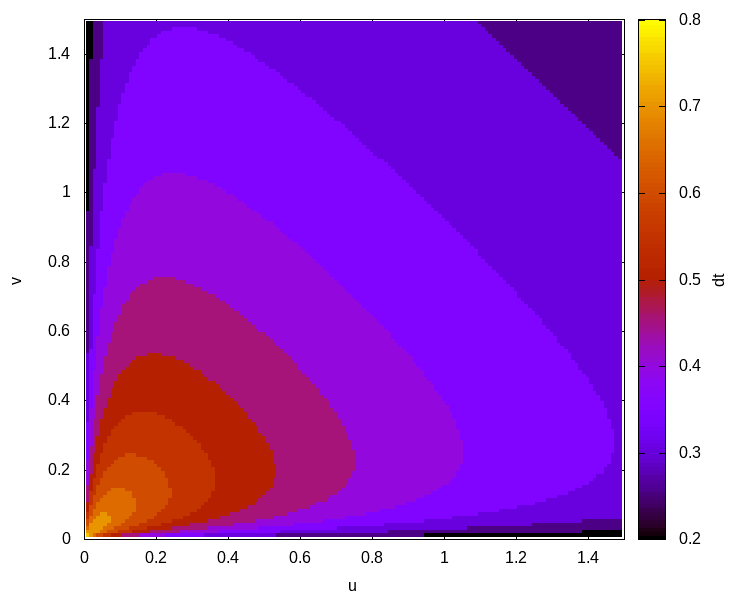}   
 \caption{Maximum $\Delta t$ (color-coded) as a function of $\bar u$, $\bar v$ for $\gamma = 1.4$, $\Delta x = 1$, $\bar c = 1$ for the linearization of a numerical method based on the flux \eqref{eq:eulerfluxx}, but without the denominator. The results were obtained by a numerical evaluation of Theorem \ref{thm:millerstability}. One observes a decrease of stability when the velocity is aligned with the grid axes.}
 \label{fig:stabilitygridaligned}
\end{figure}

Without the denominator, the numerical flux \eqref{eq:eulerfluxxmultid}, and even more so \eqref{eq:eulerfluxx}, would be a simple centered flux with an advective diffusion. Its low Mach compliance is also easier to show. Unfortunately, upon linearization, or application to the equations of linear acoustics and advection, this scheme shows a stability condition that is impractical. An investigation of the maximum allowed time step (through a numerical evaluation of Theorem \ref{thm:millerstability}) is shown in Figure \ref{fig:stabilitygridaligned}. The maximum time step for the two-dimensional method drops significantly when the velocity is aligned with one of the axes. Interestingly, the one-dimensional method does not suffer from a low CFL number. These results are also confirmed in numerical experiments performed with an implementation of this scheme for the full Euler equations.

\end{document}